%% file: main.tex
\documentclass{article}

\usepackage[final]{neurips_2022}

\usepackage[utf8]{inputenc} 
\usepackage[T1]{fontenc}    
\usepackage{hyperref}       
\usepackage{url}            
\usepackage{booktabs}       
\usepackage{amsfonts}       
\usepackage{nicefrac}       
\usepackage{microtype}      
\usepackage{bbold}
\usepackage{capt-of}

\usepackage[center]{caption}

\usepackage{amsmath, amsthm, amssymb, graphicx, cite, epsfig,epstopdf,color,soul}
\usepackage{tabularx}
\usepackage[ruled,vlined]{algorithm2e}
\allowdisplaybreaks
\usepackage{color}
\usepackage{enumerate}
\usepackage[shortlabels]{enumitem}
\usepackage{multicol}
\usepackage{empheq}
\usepackage{caption} 

\usepackage[dvipsnames]{xcolor}
\definecolor{red}{gray}{0}

\input{header}

\input{defs-mlmath.tex}

\usepackage[toc,page,header]{appendix}
\usepackage{minitoc}

\title{Regularized Gradient Descent Ascent for Two-Player Zero-Sum Markov Games}

\author{
  Sihan Zeng \\
  Dept. of Electrical and Computer Engineering\\
  Georgia Institute of Technology\\
  Atlanta, GA 30318 \\
  \texttt{szeng30@gatech.edu} \\
  \And
  Thinh Doan \\
  Dept. of Electrical and Computer Engineering \\
  Virginia Tech \\
  Blacksburg, VA 24061 \\
  \texttt{thinhdoan@vt.edu} \\
  \And
  Justin Romberg \\
  Dept. of Electrical and Computer Engineering \\
  Georgia Institute of Technology \\
  Atlanta, GA 30318 \\
  \texttt{jrom@ece.gatech.edu} \\
}

\begin{document}

\doparttoc 
\faketableofcontents 


\maketitle

\begin{abstract}

We study the problem of finding the Nash equilibrium in a two-player zero-sum Markov game. Due to its formulation as a minimax optimization program, a natural approach to solve the problem is to perform gradient descent/ascent with respect to each player in an alternating fashion. However, due to the non-convexity/non-concavity of the underlying objective function, theoretical understandings of this method are limited.
In our paper, we consider solving an entropy-regularized variant of the Markov game. The regularization introduces structure into the optimization landscape that make the solutions more identifiable and allow the problem to be solved more efficiently. Our main contribution is to show that under proper choices of the regularization parameter, the gradient descent ascent algorithm converges to the Nash equilibrium of the original unregularized problem. We explicitly characterize the finite-time performance of the last iterate of our algorithm, which vastly improves over the existing convergence bound of the gradient descent ascent algorithm without regularization. Finally, we complement the analysis with numerical simulations that illustrate the accelerated convergence of the algorithm.\looseness=-1

\end{abstract}

\input{Introduction}

\input{Preliminaries}

\input{Analysis_PG_Constanttau}

\input{Analysis_PG_PiecewiseConstant}

\input{Analysis_PG_Diminishingtau}

\input{Experiments}

\input{Conclusion}

\section*{Acknowledgement}
Sihan Zeng and Justin Romberg were supported in part by ARL DCIST CRA W911NF-17-2-0181.
The work of Thinh T. Doan was supported in part by the Commonwealth Cyber Initiative.


\bibliographystyle{plainnat}
\bibliography{references}

\newpage
\appendix

    

\vbox{%
    \hsize\textwidth
    \linewidth\hsize
    \vskip 0.1in
  \hrule height 4pt
  \vskip 0.25in
  \vskip -5.5pt%
  \centering
    {\LARGE\bf{ Appendix} \par}
      \vskip 0.29in
  \vskip -5.5pt
  \hrule height 1pt
  \vskip 0.09in%
    
  \vskip 0.2in

  }
  
\setcounter{page}{1}
For convenience, we include a table of contents for the supplementary material below. 
\vspace{-30pt}
\addcontentsline{toc}{section}{Appendix} 
\part{} 
\parttoc 

\input{Proof_Theorem}

\input{Proof_Lemma}

\input{Remark_InitialCondition}

\input{ExperimentDetails}

\end{document}

%% file: header.tex
\newcommand{\Acal}{{\cal A}}
\newcommand{\Bcal}{{\cal B}}

\newcommand{\Fcal}{{\cal F}}

\newcommand{\Hcal}{{\cal H}}

\newcommand{\Mcal}{{\cal M}}

\newcommand{\Ocal}{{\cal O}}
\newcommand{\Pcal}{{\cal P}}

\newcommand{\Scal}{{\cal S}}

\newcommand{\1}{{\mathbf{1}}}

\newcommand{\argmin}{\mathop{\rm argmin}}
\newcommand{\argmax}{\mathop{\rm argmax}}

\newtheorem{lem}{Lemma}
\newtheorem{thm}{Theorem}
\newtheorem{cor}{Corollary}

\newtheorem{assump}{Assumption}
\newtheorem{remark}{Remark}


%% file: defs-mlmath.tex
\usepackage{tikz}

%% file: Introduction.tex
\section{Introduction}
The two-player zero-sum Markov game is a special case of competitive multi-agent reinforcement learning where two agents driven by opposite reward functions jointly determine the state transition in an environment. Usually cast as a non-convex non-concave minimax optimization program, this framework finds applications in many practical problems including game playing \citep{lanctot2019openspiel,vinyals2019grandmaster}, robotics \citep{riedmiller2007experiences,shalev2016safe}, and robust policy optimization \citep{pinto2017robust}.

A convenient class of methods frequently used to solve multi-agent reinforcement learning problems is the independent learning approach. Independent learning algorithms proceed iteratively with each player taking turns to optimize its own objective while pretending that the other players' policies are fixed to their current iterates. In the context of two-player zero-sum Markov games, the independent learning algorithm performs gradient descent ascent (GDA), which alternates between the gradient updates of the two agents that seek to maximize and minimize the same value function.
Despite the popularity of such algorithms in practice, their theoretical understandings are sparse and do not follow from those in the single-agent case as the environment is not stationary from the eye of any agent. \citep{daskalakis2017training} shows that iterates of GDA can possibly diverge or be trapped in limit cycles even in the simplest single-state case when the two players learn with the same rate.

It may be tempting to analyze the two-player zero-sum Markov game by applying the existing theoretical results on minimax optimization. However, as the objective function in a Markov game is not convex or concave, current analytical tools in minimax optimization that require the objective function to be convex/concave at least on one side are inapplicable.
Fortunately, the Markov game has its own structure: it exhibits a ``gradient domination'' condition with respect to each player, which essentially guarantees that every stationary point of the value function is globally optimal. Exploiting this property, \citet{daskalakis2020independent} builds on the theory of \citet{lin2020gradient} and shows that a two-time-scale GDA algorithm converges to the Nash equilibrium of the Markov game with a complexity that depends polynomially on the specified precision.
However, deriving an explicit finite-time convergence rate is still an open problem. In addition, the analysis in \citet{daskalakis2020independent} does not guarantee the convergence of the last iterate; convergence is shown on the average of all past iterates. 

In this paper, we show that introducing an entropy regularizer into the value function significantly accelerates the convergence of GDA to the Nash equilibrium.  By dynamicially adjusting the regularization weight towards zero, we are able to give a finite-time last-iterate convergence guarantee to the Nash equilibrium of the original Markov game.


\textbf{Main Contributions}\vspace{0.1cm}\\
$\bullet$ We show that the entropy-regularized Markov game is highly structured; in particular, it obeys a condition similar to the well-known Polyak-\L ojasiewicz condition, which allows linear convergence of GDA to the (unique) equilibrium point of the regularized game with fixed regularization weight.  We also show that the distance of the equilibrium point of the regularized game to the equilibrium point of the original game can be bounded in terms of the regularizing weight.

$\bullet$ 
We show that by dynamically driving the regularization weight towards zero, we can solve the original Markov game.
We propose two approaches to reduce the regularization weight and study their finite-time convergence. The first approach uses a piecewise constant weight that decays geometrically fast, and its analysis follows as a straightforward consequence of our analysis for the case of fixed regularization weight. 
To reach a Nash equilibrium of the Markov game up to error $\epsilon$, we find that this approach requires at most $\Ocal(\epsilon^{-3})$ gradient updates, where $\Ocal$ only hides structural constants.
The second approach reduces the regularization weight online along with the gradient updates. Through a multi-time-scale analysis, we optimize the regularization weight sequence along with the step size as polynomial functions of $k$, where $k$ is the iteration index. We show that the last iterate of the GDA algorithm converges to the Nash equilibrium of the original Markov game at a rate of $\Ocal(k^{-1/3})$. 
Compared with the state-of-the-art analysis of the GDA algorithm without regularization which shows that the convergence rate of the averaged iterates is polynomial in the desired precision and all related parameters, our algorithms enjoy faster last-iterate convergence guarantees.

\input{Literature}

%% file: Literature.tex
\subsection{Related Work}
A Markov game reduces to a standard Markov Decision Process (MDP) with respect to one player if the policy of the other player is fixed. This is an important observation that allows our work to exploit the recent advances in the analysis of policy gradient methods for MDPs \citep{nachum2017bridging,neu2017unified,agarwal2020optimality,mei2020global,lan2022policy}. Various entropy-based regularizers are introduced in these works that inspire the regularization of this paper. Our particular regularization is also considered by \citet{cen2021fast2}, but we discuss and leverage structure in the regularized Markov game that was previously unknown. 

As the two-player zero-sum Markov game can be formulated a minimax optimization problem, our work relates to the vast volume of literature in this domain. Minimax optimization has been extensively studied in the case where the objective function is convex/concave with respect to at least one variable \citep{lin2020gradient,lin2020near,wang2020improved,ostrovskii2021efficient}. 
In the general non-convex non-concave setting, the problem becomes much more challenging as even the notion of stationarity is unclear \citep{jin2020local}.
{\color{red}
In \citet{nouiehed2019solving}, non-convex non-concave objective functions obeying a one--sided P\L~condition are considered, which the authors use to show the convergence of GDA. \citet{yang2020global} analyzes GDA under a two-sided P\L~condition and has a tight connection to our work as the value function of our regularized Markov game also has structure that is similar to, but weaker than, the P\L~condition on two sides.}

By exploiting the gradient domination condition of a Markov game with respect to each player, \citet{daskalakis2020independent} is the first to show that the GDA algorithm provably converges to a Nash equilibrium of a Markov game. 
A finite-time complexity is not derived in \citet{daskalakis2020independent}, but their analysis and choice of step sizes indicate that the convergence rate is at least worse than $\Ocal(k^{-1/10.5})$.
%
%
Additionally, \citet{daskalakis2020independent} does not guarantee the convergence of the last iterate, but rather analyzes the average of all iterates. In contrast, our work provides a finite-time convergence analysis on the last iterate of the GDA algorithm.

{\color{red}
While our work treats the Markov game purely from the optimization perspective, we would like to point out another related line of works that consider value-based methods \citep{perolat2015approximate,bai2020provable,xie2020learning,cen2021fast2,sayin2022fictitious}. In particular, \citet{perolat2015approximate} is among the first works to extend value-based methods from single-agent MDP to two-player Markov games. Since then, the basic techniques for analyzing value-based methods for Markov games are relatively well-known. \citet{bai2020provable} considers a value iteration algorithm with confidence bounds. In \citet{cen2021fast2}, a nested-loop algorithm is designed where the outer loop employs value iteration and the inner loop runs a gradient-descent-ascent-flavored algorithm to solve a regularized bimatrix game. In comparison, pure policy optimization algorithms are much less understood for Markov games, but this is an important subject to study due to their wide use in practice. In single-agent MDPs, value-based methods and policy optimization methods enjoy comparable convergence guarantees today, and our work aims to narrow the gap between the understanding of these two classes of algorithms in two-player Markov games.
}

Finally, we note the recent surge of interest in solving two-player games and minimax optimization programs with extragradient or optimistic gradient methods in the cases where vanilla gradient algorithms often cannot be shown to converge \citep{chavdarova2019reducing,mokhtari2020unified,li2022convergence,wei2021last,zhao2021provably,cen2021fast2,chen2021sample}.
These methods typically require multiple gradient evaluations at each iteration and are more complicated to implement. 
Most related to our work, \citet{cen2021fast2} shows the linear convergence of an extragradient algorithm for solving regularized bilinear matrix games. 
They also show that a regularized Markov game can be decomposed into a series of regularized matrix games and present a nested-loop extragradient algorithm which solves these games successively and eventually converges to the Nash equilibrium of the regularized Markov game.
The regularization weight can then be selected based on the desired precision of the unregularized problem.
Although our overall goal of finding the Nash equilibrium of a general Markov game is the same, the manner in which we decompose and analyze the problem is different.  
Our analysis here is based on GDA applied directly to a general regularized Markov game.
We show that for a fixed regularization parameter for a general Markov game, GDA has linear convergence to the modified equilibrium point.
We also give a scheduling scheme for adjusting the regularization parameter as the GDA iterations proceed, making them converge to the solution to the original problem.

%% file: Preliminaries.tex
\section{Preliminaries}\label{sec:Preliminaries}

We consider a two-player Markov game characterized by $\Mcal=(\Scal,\Acal,\Bcal,\Pcal,\gamma,r)$. Here, $\Scal$ is the finite state space, $\Acal$ and $\Bcal$ are the finite action spaces of the two players, $\gamma\in(0,1)$ is the discount factor, and $r:\Scal\times\Acal\times\Bcal\rightarrow[0,1]$ is the reward function. Let $\Delta_{\Fcal}$ denote the probability simplex over a set $\Fcal$, and $\Pcal:\Scal\times\Acal\times\Bcal\rightarrow\Delta_{\Scal}$ be the transition probability kernel, with $\Pcal(s'\mid s,a,b)$ specifying the probability of the game transitioning from state $s$ to $s'$ when the first player selects action $a\in\Acal$ and the second player selects $b\in\Bcal$.
The policies of the two players are denoted by $\pi\in\Delta_{\Acal}^{\Scal}$ and $\phi\in\Delta_{\Bcal}^{\Scal}$, with $\pi(a\mid s)$, $\phi(b\mid s)$ denoting the probability of selecting action $a$, $b$ in state $s$ according to $\pi$, $\phi$.
Given a policy pair $(\pi,\phi)$, we measure its performance in state $s\in\Scal$ by the value function
\begin{align*}
    V^{\pi, \phi}(s)=\mathbb{E}_{a_k\sim\pi(\cdot\mid s_k),b_k\sim\phi(\cdot\mid s_k),s_{k+1}\sim\Pcal(\cdot\mid s_k,a_k,b_k)}\Big[\sum\nolimits_{k=0}^{\infty} \gamma^k r\left(s_k, a_k, b_k\right) \mid s_0=s\Big].
\end{align*}
Under a fixed initial distribution $\rho\in\Delta_{\Scal}$, we define the discounted cumulative reward under $(\pi,\phi)$
\begin{align*}
    J(\pi, \phi)\triangleq\mathbb{E}_{s_0\sim\rho}[V^{\pi, \phi}(s_0)],
\end{align*}
where the dependence on $\rho$ is dropped for simplicity.
It is known that the Nash equilibrium always exists in two-player zero-sum Markov games \citep{shapley1953stochastic}, i.e. there exists an optimal policy pair $(\pi^{\star},\phi^{\star})$ such that
\begin{align}
    \max_{\pi\in\Delta_{\Acal}^{\Scal}}\min_{\phi\in\Delta_{\Bcal}^{\Scal}}J(\pi, \phi)=\min_{\phi\in\Delta_{\Bcal}^{\Scal}}\max_{\pi\in\Delta_{\Acal}^{\Scal}}J(\pi, \phi)=J(\pi^{\star},\phi^{\star}).\label{eq:obj}
\end{align}
However, as $J$ is generally non-concave with respect to the policy of the first player and non-convex with respect to that of the second player, direct GDA updates may not find $(\pi^{\star},\phi^{\star})$ and usually exhibit an oscillation behavior, which we illustrate through numerical simulations in Section~\ref{sec:experiment}. Our approach to address this issue is to enhance the structure of the Markov game through regularization.

\subsection{Entropy-Regularized Two-Player Zero-Sum Markov Games}
In this section we define the entropy regularization and discuss structure of the regularized objective function and its connection to the original problem.
Let the regularizers be
\begin{align*}
    \Hcal_{\pi}(s,\pi,\phi)&\triangleq\mathbb{E}_{a_k\sim\pi(\cdot\mid s_k),b_k\sim\phi(\cdot\mid s_k),s_{k+1}\sim\Pcal(\cdot\mid s_k,a_k, b_k)}\Big[\sum\nolimits_{k=0}^{\infty}-\gamma^{k} \log \pi\left(a_k \mid s_k\right)\mid s_0=s\Big],\notag\\
    \Hcal_{\phi}(s,\pi,\phi)&\triangleq\mathbb{E}_{a_k\sim\pi(\cdot\mid s_k),b_k\sim\phi(\cdot\mid s_k),s_{k+1}\sim\Pcal(\cdot\mid s_k,a_k, b_k)}\Big[\sum\nolimits_{k=0}^{\infty}-\gamma^{k} \log \phi\left(b_k \mid s_k\right)\mid s_0=s\Big].\notag\\
\end{align*}
We define the regularized value function
\begin{align*}
    V_{\tau}^{\pi,\phi}(s)&\triangleq V^{\pi,\phi}(s)+\tau\Hcal_{\pi}(s,\pi,\phi)-\tau\Hcal_{\phi}(s,\pi,\phi)\notag\\
    &=\mathbb{E}_{\pi,\phi,\Pcal}\Big[\sum\nolimits_{k=0}^{\infty} \gamma^k \Big(r\left(s_k, a_k, b_k\right)-\tau\log\pi(a_k\mid s_k)+\tau\log\phi(b_k\mid s_k)\Big) \mid s_0=s\Big],
\end{align*}
where $\tau\geq0$ is a weight parameter. Again under a fixed initial distribution $\rho\in\Delta_{\Scal}$ we denote $J_{\tau}(\pi,\phi)\triangleq\mathbb{E}_{s\sim\rho}[V_{\tau}^{\pi,\phi}(s)]$. The regularized advantage function is
\begin{align*}
    A_{\tau}^{\pi,\phi}(s,a,b)&\triangleq r(s,a,b)-\tau\log\pi(a\mid s)+\tau\log\phi(b\mid s)+\gamma\mathbb{E}_{s'\sim\Pcal(\cdot\mid s,a,b)}\left[V_{\tau}^{\pi,\phi}(s')\right]-V_{\tau}^{\pi,\phi}(s),
\end{align*}
which later helps us to express the policy gradient.

We use $d_{\rho}^{\pi,\phi}\in\Delta_{\Scal}$ to denote the discounted visitation distribution under any policy pair $(\pi,\phi)$ and the initial state distribution $\rho$
\[d_{\rho}^{\pi,\phi}(s)\triangleq(1-\gamma)\mathbb{E}_{\pi,\phi,\Pcal}\Big[\sum\nolimits_{k=0}^{\infty} \gamma^k \1(s_k=s) \mid s_0\sim\rho\Big]\]
For sufficient state visitation, we assume that the initial state distribution is bounded away from zero. This is a standard assumption in the entropy-regularized MDP literature \citep{mei2020global, ying2022dual}.
\begin{assump}\label{assump:positive_rho}
The initial state distribution $\rho$ is strictly positive for any state, and we denote $\rho_{\min}=\min_{s\in\Scal}\rho(s)>0$.
\end{assump}



When the policy of the first player is fixed to $\pi\in\Delta_{\Acal}^{\Scal}$, the Markov game reduces to an MDP for the second player with state transition probability $\widetilde{\Pcal}_{\phi}(s'\mid s,b)=\sum_{a\in\Acal}\Pcal(s'\mid s,a,b)\pi(a\mid s)$ and reward function $\widetilde{r}_{\phi}(s,b)=\sum_{a\in\Acal}r(s,a,b)\pi(a\mid s)$. A similar argument holds for the first player if the second player's policy is fixed.
To denote the operators that map one player's policy to the best response of the other player and the corresponding value function, we define
\begin{gather}
    \pi_{\tau}(\phi)\triangleq\argmax_{\pi\in\Delta_{\Acal}^{\Scal}} J_{\tau}(\pi,\phi),\quad \phi_{\tau}(\pi)\triangleq\argmin_{\phi\in\Delta_{\Bcal}^{\Scal}} J_{\tau}(\pi,\phi),\notag\\
    g_{\tau}(\pi)\triangleq\min_{\phi\in\Delta_{\Bcal}^{\Scal}}J_{\tau}(\pi,\phi)=J_{\tau}(\pi,\phi_{\tau}(\pi)).\label{eq:def_g}
\end{gather}

For any $\tau>0$, the following lemma bounds the performance difference between optimal and sub-optimal policies and establishes the uniqueness of $\pi_{\tau}(\phi)$ and $\phi_{\tau}(\pi)$. When $\tau=0$, we use $\pi_{0}(\phi)$ and $\phi_{0}(\pi)$ to denote one of the maximizers and minimizers since they may not be unique.

\begin{lem}[Performance Difference]\label{lem:quadratic_growth}
Under Assumption~\ref{assump:positive_rho} and given $\tau>0$, $\pi_{\tau}(\phi)$ is unique for any $\phi\in\Delta_{\Bcal}^{\Scal}$, and $\phi_{\tau}(\pi)$ is unique for any $\pi\in\Delta_{\Acal}^{\Scal}$. 
Given any min player policy $\phi\in\Delta_{\Bcal}^{\Scal}$,
\begin{align}
    J_{\tau}(\pi_{\tau}(\phi),\phi) -J_{\tau}(\pi,\phi)\geq\frac{\tau\rho_{\min}}{2\log(2)}\|\pi_{\tau}(\phi)-\pi\|^2,\quad\forall \pi\in\Delta_{\Acal}^{\Scal}.
    \label{lem:quadratic_growth:eq1}
\end{align}

Given any max player policy $\pi\in\Delta_{\Acal}^{\Scal}$,
\begin{align}
    J_{\tau}(\pi,\phi_{\tau}(\pi)) -J_{\tau}(\pi,\phi)\leq-\frac{\tau\rho_{\min}}{2\log(2)}\|\phi_{\tau}(\pi)-\phi\|^2,\quad\forall \phi\in\Delta_{\Bcal}^{\Scal}.
    \label{lem:quadratic_growth:eq2}
\end{align}
\end{lem}
The Nash equilibrium of the regularized problem is sometimes referred to as the quantal response equilibrium \citep{mckelvey1995quantal} and is known to exist under any $\tau$.
Leveraging Lemma~\ref{lem:quadratic_growth}, we formally state the conditions guaranteeing its existence and affirm that it is unique.

\begin{lem}[Minimax Theorem for Entropy-Regularized Markov Game]\label{lem:minimax_regMG}
Under Assumption~\ref{assump:positive_rho}, for any regularization weight $\tau>0$, there exists a unique Nash equilibrium policy pair $(\pi_{\tau}^{\star},\phi_{\tau}^{\star})$ such that\looseness=-1
\begin{align}
    \max_{\pi\in\Delta_{\Acal}^{\Scal}}\min_{\phi\in\Delta_{\Bcal}^{\Scal}}J_{\tau}(\pi,\phi)=\min_{\phi\in\Delta_{\Bcal}^{\Scal}}\max_{\pi\in\Delta_{\Acal}^{\Scal}}J_{\tau}(\pi,\phi)=J_{\tau}(\pi_{\tau}^{\star},\phi_{\tau}^{\star}).\label{lem:minimax_regMG:eq0}
\end{align}
\end{lem}

We are only interested in the solution of the regularized Markov game if it gives us knowledge of the original problem in \eqref{eq:obj}. In the following lemma, we show that the distance between the Nash equilibrium of the regularized game and that of the original one is bounded by the regularization weight. This is an important condition guaranteeing that we can find an approximate solution to the original Markov game by solving the regularized problem.
In addition, this lemma also shows that the same policy pair produces value functions with bounded distance under two regularization weights.

\begin{lem}\label{lem:V_tau_diff}
For any $\tau\geq\tau'\geq 0$ and policy $\pi$,
\begin{gather}
    -(\tau-\tau')\log|\Bcal|\leq J_{\tau}(\pi_{\tau}^{\star},\phi_{\tau}^{\star})-J_{\tau'}(\pi_{\tau'}^{\star},\phi_{\tau'}^{\star})\leq (\tau-\tau')\log|\Acal|.\label{lem:V_tau_diff:eq3}\\
    -(\tau-\tau')\log|\Bcal|\leq g_{\tau}(\pi) - g_{\tau'}(\pi)=J_{\tau}(\pi,\phi_{\tau}(\pi))-J_{\tau'}(\pi,\phi_{\tau'}(\pi))\leq(\tau-\tau')\log|\Acal|.\label{lem:V_tau_diff:eq4}\\
    -\frac{\tau-\tau'}{1-\gamma}\log|\Bcal|\leq J_{\tau}(\pi,\phi)-J_{\tau'}(\pi,\phi)\leq \frac{\tau-\tau'}{1-\gamma}\log|\Acal|.\label{lem:V_tau_diff:eq1}
\end{gather}

\end{lem}

\subsection{Softmax Parameterization}\label{sec:softmax}

In this work we use a tabular softmax policy parameterization and maintain two tables $\theta\in\mathbb{R}^{\Scal\times\Acal}$, $\psi\in\mathbb{R}^{\Scal\times\Bcal}$ that parameterize the policies of the two players according to
\begin{align*}
    \pi_{\theta}(a \mid s)=\frac{\exp \left(\theta(s, a)\right)}{\sum_{a' \in \Acal} \exp \left(\theta(s, a')\right)},\quad \text {and} \quad \phi_{\psi}(b \mid s)=\frac{\exp \left(\psi(s, b)\right)}{\sum_{b' \in \Acal} \exp \left(\psi(s, b')\right)}.
\end{align*}
The gradients of the regularized value function with respect to the policy parameters admit closed-form expressions
\begin{align*}
    &\frac{\partial J_{\tau}(\pi_{\theta},\phi_{\psi})}{\partial \theta(s,a)}=\frac{1}{1-\gamma} d_{\rho}^{\pi_{\theta},\phi_{\psi}}(s)\pi_{\theta}(a\mid s) \sum\nolimits_{b\in\Bcal}\phi_{\psi}(b\mid s)A_{\tau}^{\pi_{\theta},\phi_{\psi}}(s, a, b),\\
    &\frac{\partial J_{\tau}(\pi_{\theta},\phi_{\psi})}{\partial \psi(s,b)}=\frac{1}{1-\gamma} d_{\rho}^{\pi_{\theta},\phi_{\psi}}(s)\phi_{\psi}(b\mid s) \sum\nolimits_{a\in\Acal}\pi_{\theta}(a\mid s)A_{\tau}^{\pi_{\theta},\phi_{\psi}}(s, a, b),
\end{align*}
and computing them exactly requires knowledge of the dynamics of the environment. Note that the gradients of value function and the regularizer are Lipschitz with respect to the policy parameters with constants $L_V=\frac{8}{(1-\gamma)^3}$ and $L_{\Hcal}=\frac{4+8\log|\Acal|}{(1-\gamma)^3}$. This property is more formally stated and proved in Lemmas~\ref{lem:value_LipschitzGrad} and \ref{lem:reg_LipschitzGrad} of the appendix.

We next present an important property that we will later exploit to study the convergence of the GDA updates to the solution of the regularized Markov game.
Under the softmax parameterization, the regularized value function enjoys a gradient domination condition with respect to the policy parameter that resembles the P\L~condition.
\begin{lem}[PL-Type Condition]\label{lem:nonuniform_PL}
Under Assumption \ref{assump:positive_rho}, we have for any $\theta\in\mathbb{R}^{\Scal\times\Acal}$ and $\psi\in\mathbb{R}^{\Scal\times\Bcal}$
\begin{align*}
    &\|\nabla_{\theta} J_{\tau}(\pi_{\theta},\phi_{\psi})\|^2\geq\frac{2(1-\gamma)\tau\rho_{\min}^2}{|\Scal|}\left(\min_{s,a}\pi_{\theta}(a\mid s)\right)^2\left(J_{\tau}(\pi_{\tau}(\phi_{\psi}),\phi_{\psi})-J_{\tau}(\pi_{\theta},\phi_{\psi})\right),\\
    &\|\nabla_{\psi} J_{\tau}(\pi_{\theta},\phi_{\psi})\|^2\geq\frac{2(1-\gamma)\tau\rho_{\min}^2}{|\Scal|}\left(\min_{s,b}\phi_{\psi}(b\mid s)\right)^2\left(J_{\tau}(\pi_{\theta},\phi_{\psi})-J_{\tau}(\pi_{\theta},\phi_{\tau}(\pi_{\theta}))\right).
\end{align*}
\end{lem}
The P\L~condition is a tool commonly used in the optimization community to show the linear convergence of the gradient descent algorithm \citep{karimi2016linear,yu2019computation,khaled2020better,zeng2021two}. 
The condition in Lemma~\ref{lem:nonuniform_PL} is weaker than the common P\L~condition in two aspects. First, our P\L~coefficient is a function of the smallest policy entry. When we seek to bound the gradient of the iterates $\|\nabla_{\theta} J_{\tau}(\pi_{\theta_k},\phi_{\psi_k})\|^2$ and $\|\nabla_{\psi} J_{\tau}(\pi_{\theta_k},\phi_{\psi_k})\|^2$ later in the analysis, the P\L~coefficients will depend on $\min_{s,a}\pi_{\theta_k}(a\mid s)$ and $\min_{s,b}\phi_{\psi_k}(b\mid s)$, which may not be lower bounded by any positive constant. Second, the coefficients involve $\tau$, which is not a constant but needs to be carefully chosen to control the error between the regularized problem and the original one.

%% file: Analysis_PG_Constanttau.tex
\section{Solving Regularized Markov Games}

Leveraging the structure introduced in Section~\ref{sec:Preliminaries}, our first aim is to establish the finite-time convergence of the GDA algorithm to the Nash equilibrium of the regularized Markov game under a fixed regularization weight $\tau>0$. 
The GDA algorithm executes the updates
\begin{align}
    \theta_{k+1}&=\theta_k+\alpha_k\nabla_{\theta}J_{\tau}(\pi_{\theta_k},\phi_{\psi_k}),\qquad\psi_{k+1}=\psi_k-\beta_k \nabla_{\psi}J_{\tau}(\pi_{\theta_{k+1}},\phi_{\psi_k}).
    \label{update:GDA_constanttau}
\end{align}
The convergence bound we will derive reflects a trade-off for the regularization weight $\tau$: when $\tau$ is large, we get faster convergence to the Nash equilibrium of the regularized problem, but it is farther away from the Nash equilibrium of the original one. The result in this section will inspire the $\tau$ adjustment schemes designed later in the paper to achieve the best possible convergence to the Nash equilibrium of the original unregularized Markov game. 
%

It can be shown that the Nash equilibrium of the regularized Markov game is a pair of completely mixed policies, i.e. $\forall\tau\hspace{-2pt}>\hspace{-2pt}0$ there exists $c_{\tau}\hspace{-2pt}>\hspace{-2pt}0$ such that $\min_{s,a}\hspace{-1pt}\pi_{\tau}^{\star}(a\mid s)\hspace{-2pt}\geq\hspace{-2pt} c_{\tau}$, and $\min_{s,b}\hspace{-1pt}\phi_{\tau}^{\star}(b\mid s)\hspace{-2pt}\geq\hspace{-2pt} c_{\tau}$ \citep{nachum2017bridging}.
In this work, we further assume the existence of a uniform lower bound on the entries of $(\pi_{\tau}^{\star},\phi_{\tau}^{\star})$ across $\tau$. 
We provide more explanation of the assumption in Remark~\ref{remark:completelymixed}.

\begin{assump}\label{assump:NE_completelymixed}
There exists a positive constant $c$ (independent of $\tau$) such that for any $\tau>0$
\begin{align*}
    \min_{s,a}\pi_{\tau}^{\star}(a\mid s)\geq c,\quad\min_{s,b}\phi_{\tau}^{\star}(b\mid s)\geq c.
\end{align*}
\end{assump}
To measure the convergence of the iterates to the Nash equilibrium of the regularized Markov game, we recall the definition of $g_{\tau}$ in \eqref{eq:def_g} and define
\begin{align}
\delta^{\pi}_k =J_{\tau}(\pi_{\tau}^{\star},\phi_{\tau}^{\star})-g_{\tau}(\pi_{\theta_k}),\quad\delta^{\phi}_k = J_{\tau}(\pi_{\theta_k},\phi_{\psi_k})-g_{\tau}(\pi_{\theta_k}).\label{eq:conv_metric_constanttau}
\end{align}
The convergence metric is asymmetric for two players: the first player is quantified by its performance when the second player takes the most adversarial policy, while the second player is evaluated under the current policy iterate of the first player. We note that $\delta^{\pi}_k$ and $\delta^{\phi}_k$ are non-negative, and $\delta^{\pi}_k=\delta^{\phi}_k=0$ implies that $(\pi_{\theta_k},\phi_{\psi_k})$ is the Nash equilibrium.
Under this convergence metric, the following theorem states that the GDA updates in \eqref{update:GDA_constanttau} solve the regularized Markov game linearly fast. The proofs of the theoretical results of this paper are presented in Section~\ref{sec:proof:thm} of the appendix.


    

\begin{thm}\label{thm:main_constanttau}
We define $L=3 L_{\Hcal}\max\{\tau,1\}$, $C_1=\frac{\rho_{\min}c^2}{64\log(2)}$, and $C_2=\frac{2\sqrt{|\Scal|}}{\sqrt{(1-\gamma)\rho_{\min}}c}$, and choose the initial policy parameters to be $\theta_0=0\in\mathbb{R}^{|\Scal|\times|\Acal|}$ and $\psi_0=0\in\mathbb{R}^{|\Scal|\times|\Bcal|}$ (the initial policies $\pi_{\theta_0}$ and $\phi_{\psi_0}$ are uniform).
Let the step sizes of \eqref{update:GDA_constanttau} be
\begin{align*}
    \alpha_k=\alpha,\quad \beta_k=\beta,
\end{align*}
with $\alpha$, $\beta$ satisfying
\begin{align*}
    &\max\{\alpha,\beta\}\leq\frac{1}{L},\,\, \frac{\alpha}{\beta}\leq\min\{\frac{ (1-\gamma)\rho_{\min}^3 c^2 \tau^2}{152\log(2)|\Scal|L^2},8\},\,\, \alpha\leq\min\{(L+\frac{C_2 L^2}{\tau})^{-1},\frac{16|\Scal|}{(1-\gamma)\rho_{\min}^2 c^2 \tau}\}.
\end{align*}
If Assumption~\ref{assump:positive_rho} holds and
\begin{align}
    3\delta^{\pi}_{0}+\delta^{\phi}_{0}\leq C_1\tau,
    \label{thm:main_constanttau:eq1}
\end{align}
then the iterates of \eqref{update:GDA_constanttau} satisfy for all $k\geq0$
\begin{align*}
    3\delta_k^{\pi}+\delta_k^{\phi}\leq (1-\frac{(1-\gamma)\alpha \tau \rho_{\min}^2 c^2}{32|\Scal|})^{k}(3\delta_0^{\pi}+\delta_0^{\phi}).
\end{align*}
\end{thm}
Theorem \ref{thm:main_constanttau} establishes the linear convergence of the iterates of \eqref{update:GDA_constanttau} to the Nash equilibrium of \eqref{lem:minimax_regMG:eq0}, provided that the initial condition \eqref{thm:main_constanttau:eq1} is satisfied. 
The convergence is faster when $\tau$ is large and slower when $\tau$ is small.
Choosing $\tau$ to be large enough guarantees the initial condition (see Section~\ref{remark:large_tau0} of the appendix for more discussion) but causes the Nash equilibrium of the regularized Markov game to be distant from that of the original Markov game. This motivates us to make the regularization weight a decaying sequence that starts off large enough to meet the initial condition and becomes smaller over time to narrow the gap between the regularized Markov game and the original one. We discuss two such schemes of reducing the regularization weight in the next section.


%% file: Analysis_PG_PiecewiseConstant.tex
\section{Main Results - Solving the Original Markov Game}

This section presents two approaches to adjust the regularization weight that allow the GDA algorithm to converge to the Nash equilibrium of the original Markov game. The first approach uses a piecewise constant weight and results in the nested-loop updates stated in Algorithm \ref{Alg:GDA_piecewiseconstanttau}. In the inner loop the regularization weight and step sizes are fixed, and the two players update their policy iterates towards the Nash equilibrium of the regularized Markov game. The outer loop iteration reduces the regularization weight to make the regularized Markov game approach the original one. The regularization weight decays geometrically in the outer loop, i.e. $\tau_{t+1}=\eta\tau_t$, where $\eta\in(0,1)$ must be carefully balanced. On the one hand, recalling the definition of $g_{\tau}$ in \eqref{eq:def_g} and defining 
\begin{align*}
    \delta^{\pi}_{t,k} =J_{\tau_t}(\pi_{\tau_t}^{\star},\phi_{\tau_t}^{\star})-g_{\tau_t}(\pi_{\theta_{t,k}}),\quad\delta^{\phi}_{t,k} = J_{\tau_t}(\pi_{\theta_{t,k}},\phi_{\psi_{t,k}})-g_{\tau_t}(\pi_{\theta_{t,k}}),
\end{align*}
we need $\eta$ to be large enough that if $\theta_{t,0}$ and $\psi_{t,0}$ observe the initial condition
$3\delta^{\pi}_{t,0}+\delta^{\phi}_{t,0}\leq C_1\tau_{t}$,
then so do $\theta_{t+1,0}$ and $\psi_{t+1,0}$ in the worst case.
On the other hand, an $\eta$ selected excessively large makes the reduction of $\tau_t$ too slow to achieve the best possible convergence rate. Our next theoretical result, as a corollary of Theorem~\ref{thm:main_constanttau}, properly chooses $\eta$ and $K_t$ and establishes the convergence of Algorithm~\ref{Alg:GDA_piecewiseconstanttau} to the Nash equilibrium of the original original problem.

\begin{algorithm}[!h]
\SetAlgoLined
\textbf{Initialize:} Policy parameters $\theta_{0,0}=0\in\mathbb{R}^{\Scal\times\Acal}$ and $\psi_{0,0}=0\in\mathbb{R}^{\Scal\times\Bcal}$, step size sequences $\{\alpha_t\}$ and $\{\beta_t\}$, an initial regularization parameter $\tau_0$

\For{$t=0,1,\cdots,T$}{

 \For{$k=0,1,\cdots,K_t-1$}{
    1) Max player update:
    \begin{align*}
        \theta_{t,k+1}=\theta_{t,k}+\alpha_t\nabla_{\theta}J_{\tau}(\pi_{\theta_{t,k}},\phi_{\psi_{t,k}})
    \end{align*}

    2) Min player update: 
    \begin{align*}
        \psi_{t,k+1}=\psi_{t,k}-\beta_t \nabla_{\psi}J_{\tau}(\pi_{\theta_{t,k+1}},\phi_{\psi_{t,k}})
    \end{align*}
 }
 
 Set initial policies for next outer loop iteration $\theta_{t+1,0}=\theta_{t,K_t}$, $\psi_{t+1,0}=\psi_{t,K_t}$
 
 Reduce regularization weight $\tau_{t+1}=\eta\tau_t$ and properly adjust $\alpha_t,\beta_t$
}
\caption{Nested-Loop Policy Gradient Descent Ascent Algorithm with Piecewise Constant Regularization Weight}
\label{Alg:GDA_piecewiseconstanttau}
\end{algorithm}

\begin{cor}\label{cor:GDA_piecewiseconstanttau}
Suppose that Assumption~\ref{assump:positive_rho}-\ref{assump:NE_completelymixed} hold and $\tau_0$ is chosen such that $3\delta_{0,0}^{\pi}+\delta_{0,0}^{\phi}\leq C_1\tau_0$\footnote{This inequality is guaranteed to hold with a large enough $\tau_0$ if $\pi_{\theta_0}$ and $\phi_{\psi_0}$ are initialized to be uniform. See Section \ref{remark:large_tau0} of the appendix for more discussion.}.
We choose $\eta=\frac{C_1+2L_{\delta}}{2C_1+2L_{\delta}}$, where $L_{\delta}=4\log|\Acal|+3\log|\Bcal|+\frac{\log|\Bcal|}{1-\gamma}$  and $C_1$ is defined in Theorem~\ref{thm:main_constanttau}.
Then, under proper choices of $\alpha_t$ and $\beta_t$, the iterates of Algorithm \ref{Alg:GDA_piecewiseconstanttau} converge to a point such that
\begin{align}
    J(\pi^{\star},\phi^{\star}) - g_0(\pi_{\theta_{T,0}})\leq\epsilon\quad\text{and}\quad J(\pi_{\theta_{T,0}},\phi_{\psi_{T,0}})-g_0(\pi_{\theta_{T,0}})\leq\epsilon
    \label{cor:GDA_piecewiseconstanttau:eq1}
\end{align}
in at most $T=\Ocal(\log(\epsilon^{-1}))$ outer loop iterations.
The total number of gradient updates required is $\sum_{t=0}^{T}K_t=\Ocal(\epsilon^{-3})$.

\end{cor}

Corollary~\ref{cor:GDA_piecewiseconstanttau} guarantees that $(\pi_{\theta_T},\phi_{\psi_T})$ converge to an $\epsilon$-approximate Nash equilibrium of the original Markov game in $T=\Ocal(\epsilon^{-3})$ gradient steps. In order to achieve this rate, $K_t$ has to be adjusted along with $\tau_t$: we need $K_t=\Ocal(\tau_t^{-3})$ when $\tau_t$ becomes smaller than 1. The varying number of inner loop iterations may cause inconvenience for practical implementation. To address this issue, we next propose another scheme of adjusting the regularization weight that is carried out online along with the update of the policy iterates.

%% file: Analysis_PG_Diminishingtau.tex
\begin{algorithm}[!h]
\SetAlgoLined
\textbf{Initialize:} Policy parameters $\theta_0=0\in\mathbb{R}^{\Scal\times\Acal}$ and $\psi_0=0\in\mathbb{R}^{\Scal\times\Bcal}$, step size sequences $\{\alpha_k\}$ and $\{\beta_k\}$, regularization parameter sequence $\{\tau_k\}$

 \For{$k=0,1,\cdots,K$}{
    1) Max player update:
    \begin{align*}
        \theta_{k+1}=\theta_k+\alpha_k\nabla_{\theta}J_{\tau_k}(\pi_{\theta_k},\phi_{\psi_k})
    \end{align*}

    2) Min player update: 
    \begin{align*}
        \psi_{k+1}=\psi_k-\beta_k \nabla_{\psi}J_{\tau_k}(\pi_{\theta_{k+1}},\phi_{\psi_k})
    \end{align*}

 }
\caption{Policy Gradient Descent Ascent Algorithm with Diminishing Regularization Weight}
\label{Alg:GDA}
\end{algorithm}

Presented in Algorithm~\ref{Alg:GDA}, the second approach is a single-loop algorithm that reduces the regularization weight as a polynomial function of the iteration $k$.
We define the auxiliary convergence metrics
\begin{align*}
    \delta^{\pi}_k =J_{\tau_k}(\pi_{\tau_k}^{\star},\phi_{\tau_k}^{\star})-g_{\tau_k}(\pi_{\theta_k}),\quad \delta^{\phi}_k = J_{\tau_k}(\pi_{\theta_k},\phi_{\psi_k})-g_{\tau_k}(\pi_{\theta_k}),
\end{align*}
which measure the convergence of $(\pi_{\theta_k},\phi_{\psi_k})$ to the Nash equilibrium of the Markov game regularized with weight $\tau_k$.
To judge the performance of the iterates in the original Markov game, we are ultimately interested in bounding $J(\pi^{\star},\phi^{\star})-g_{0}(\pi_{\theta_k})$ and $J(\pi_{\theta_k},\phi_{\psi_k})-g_{0}(\pi_{\theta_k})$. Thanks to Lemma~\ref{lem:V_tau_diff}, we can quantify how fast $\delta^{\pi}_k$ and $\delta^{\phi}_k$ approach these desired quantities as $\tau_k$ decays to 0. Under an initial condition on $\delta^{\pi}_k$ and $\delta^{\phi}_k$, we now establish the convergence rate of Algorithm~\ref{Alg:GDA} to $(\pi^{\star},\phi^{\star})$ of \eqref{eq:obj} through a multi-time-scale analysis.

\begin{thm}\label{thm:main}
Let the step sizes and regularization parameter be
\begin{align*}
    \alpha_k=\frac{\alpha_0}{(k+h)^{2/3}},\quad \beta_k=\beta_0,\quad \tau_k=\frac{\tau_0}{(k+h)^{1/3}},
\end{align*}
with $\alpha_0$, $\beta_0$, $\tau_0$, and $h\geq 1$ satisfying a system of inequalities discussed in details in the analysis.
Under Assumption~\ref{assump:positive_rho}-\ref{assump:NE_completelymixed}, the iterates of Algorithm~\ref{Alg:GDA} satisfy for all $k\geq0$
\begin{align}
    J(\pi^{\star},\phi^{\star}) - g_0(\pi_{\theta_k})&\leq\frac{C_1\tau_0+3(\log|\Acal|+\log|\Bcal|)\tau_0}{3{(k+h)^{1/3}}},\label{thm:main:eq1}\\
    J(\pi_{\theta_k},\phi_{\psi_k})- g_0(\pi_{\theta_k})&\leq \frac{(1-\gamma)C_1\tau_0+(\log|\Acal|+\log|\Bcal|)\tau_0}{(1-\gamma){(k+h)^{1/3}}}\label{thm:main:eq2},
\end{align}
where the constant $C_1$ is defined in Theorem~\ref{thm:main_constanttau}.
\end{thm}

Theorem~\ref{thm:main} states that the last iterate of Algorithm~\ref{Alg:GDA} converges to an $\Ocal(k^{-1/3})$-approximate Nash equilibrium of the original Markov game in $k$ iterations. This translates to the same sample complexity as Algorithm~\ref{Alg:GDA_piecewiseconstanttau} derived in Corollary~\ref{cor:GDA_piecewiseconstanttau}. Compared with Algorithm~\ref{Alg:GDA_piecewiseconstanttau}, reducing $\tau_k$ online along with the gradient updates in a single loop simplifies the algorithm and makes tracking the regularization weight, step sizes, and policy iterates simpler and more convenient. 
We note that the techniques in \citet{daskalakis2020independent} may be used to analyze the finite-time performance of GDA for Markov games and lead to a convergence rate at least worse than $\Ocal(k^{-1/10.5})$, which we improve over.

\begin{remark}\label{remark:completelymixed}
{\color{red}
Assumption~\ref{assump:NE_completelymixed} is a restrictive assumption that does not seem necessary but rather arises as an artifact of the current analysis. When we apply the weaker PL-type condition (Lemma~\ref{lem:nonuniform_PL}) in the analysis, the entries of the iterates $\pi_{\theta_k},\phi_{\psi_k}$ need to be uniformly lower bounded, which is difficult to establish using the game structure. We come up with an innovative induction approach to quantify the connection between $\min_{s,a}\pi_{\theta_k}(a\mid s),\min_{s,b}\phi_{\psi_k}(b\mid s)$ and the optimal gap $\delta_{k}^{\pi}, \delta_{k}^{\phi}$. This approach allows us to transform the uniform lower bound requirement on $\pi_{\theta_k},\phi_{\psi_k}$ to that on the Nash equilibrium, leading to Assumption~\ref{assump:NE_completelymixed}. It is a future work to remove/relax this assumption.
}

A Markov game is said to be completely mixed if every Nash equilibrium of the game consists of a pair of completely mixed policies, i.e. $\min_{s,a}\pi^{\star}(a\mid s)>0,\min_{s,b}\phi^{\star}(b\mid s)>0$ for any Nash equilibrium $(\pi^{\star},\phi^{\star})$ of the Markov game (if more than one exists). Assumption~\ref{assump:NE_completelymixed} intuitively seems no stronger than requiring the original Markov game to be completely mixed. If the original Markov game has at least one completely mixed Nash equilibrium, the Nash equilibrium of the regularized Markov game should also be completely mixed even when the regularization weight is small, since the regularization encourages the solution to be more uniform. The reward function that results in completely mixed Markov games is well studied in \citet{raghavan1978completely,kaplansky1995contribution,das2017completely}. 
\end{remark}

%% file: Experiments.tex
\section{Numerical Simulations}\label{sec:experiment}

In this section, we numerically verify the convergence of Algorithm~\ref{Alg:GDA} on small-scale synthetic Markov games. 
Our aim is to confirm that the algorithm indeed converges rather than to visualize the exact convergence rate, as achieving the theoretical rate derived in Theorem~\ref{thm:main} requires very careful selection of all involved parameters. Considering an environment with $|\Scal|=2$ and $|\Acal|=|\Bcal|=2$, we first choose the reward and transition probability kernel such that the Markov game is completely mixed\footnote{To create a completely mixed game with $|\Acal|=|\Bcal|=2$, we simply need to choose the reward function such that $r(s,\cdot,\cdot)$ as a 2x2 matrix is diagonal dominant or sub-diagonal dominant for any state $s\in\Scal$, and we can use an arbitrary transition probability kernel. The exact choice of the reward function and transition kernel as well as the Nash equilibrium of this Markov game are presented in Section~\ref{sec:experimentdeatils} of the appendix.}.

\begin{figure}[h]
  \centering
  \includegraphics[width=\linewidth]{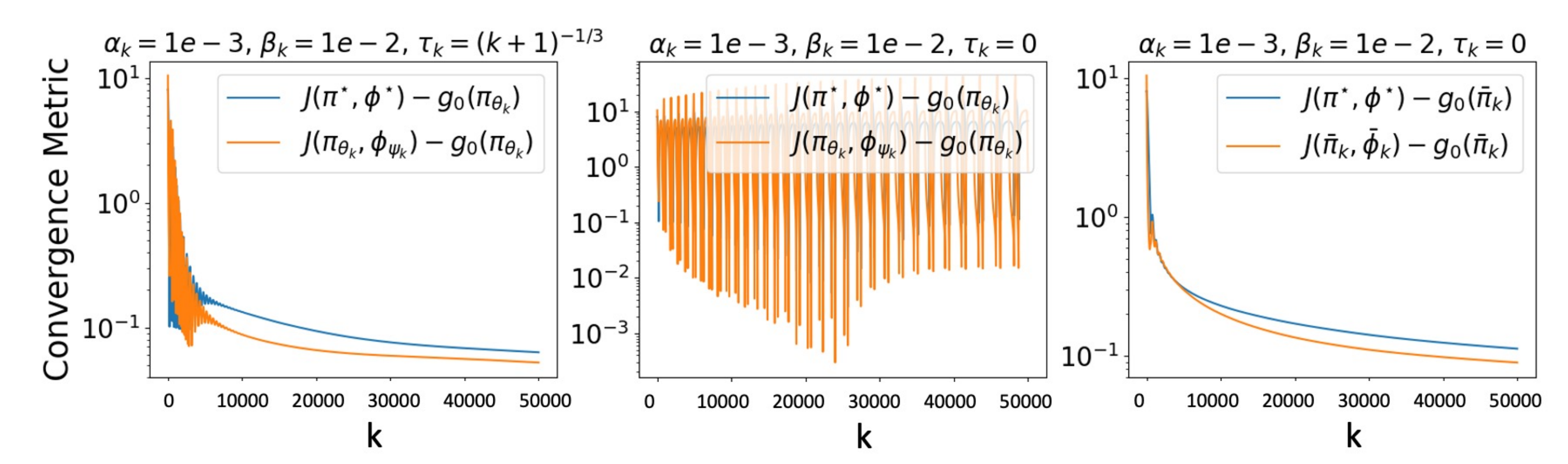}
  \caption{Convergence of GDA for a Completely Mixed Markov game}
  \label{fig:MG_mixed}
\end{figure}
We run Algorithm~\ref{Alg:GDA} for $50000$ iterations with $\alpha_k=10^{-3}$, $\beta_k=10^{-2}$, $\tau_k=(k+1)^{-1/3}$, and measure the convergence of $\pi_k$ and $\phi_k$ by metrics considered in \eqref{thm:main:eq1} and \eqref{thm:main:eq2} of Theorem~\ref{thm:main}. As shown in the first plot of Figure~\ref{fig:MG_mixed}, the last iterate exhibits an initial oscillation behavior but converge smoothly after 10000 iterations. In comparison, we visualize the convergence of the last iterate and averaged iterate of the GDA algorithm without any regularization (second and third plots of Figure~\ref{fig:MG_mixed}), where the average is computed with equal weights as $\bar{\pi}_k=\frac{1}{k+1}\sum\nolimits_{t=0}^{k}\pi_{\theta_t}, \bar{\phi}_k=\frac{1}{k+1}\sum\nolimits_{t=0}^{k}\phi_{\psi_t}$.
The existing theoretical results in this case guarantee the convergence of the averaged iterate but not the last iterate \citep{daskalakis2020independent}. According to our simulations, the last iterate indeed does not converge, while the averaged iterate does, but at a slower rate than the convergence of the last iterate of the GDA algorithm under the decaying regularization.

The theoretical results derived in this paper rely on Assumption~\ref{assump:NE_completelymixed}. To investigate whether this assumption is truly necessary, we also apply Algorithm~\ref{Alg:GDA} to a Markov game that has a deterministic Nash equilibrium and does not observe Assumption~\ref{assump:NE_completelymixed}\footnote{The detailed description of the game is again deferred to Section~\ref{sec:experimentdeatils} of the appendix.}. As illustrated in Figure~\ref{fig:MG_deterministic}, the experiment shows that Algorithm~\ref{Alg:GDA} still converges correctly to $(\pi^{\star},\phi^{\star})$ of \eqref{eq:obj}. This observation suggests that Assumption~\ref{assump:NE_completelymixed} may be an artifact of the current analysis and motivates for us to investigate ways to remove/relax this assumption in the future.
We note that the pure GDA approach without regularization also has a last-iterate convergence and does not exhibit the oscillation behavior observed in Figure~\ref{fig:MG_mixed}, since the gradients of both players never change signs regardless of the policy of the opponent in this Markov game.

\begin{figure}[h]
  \centering
  \includegraphics[width=\linewidth]{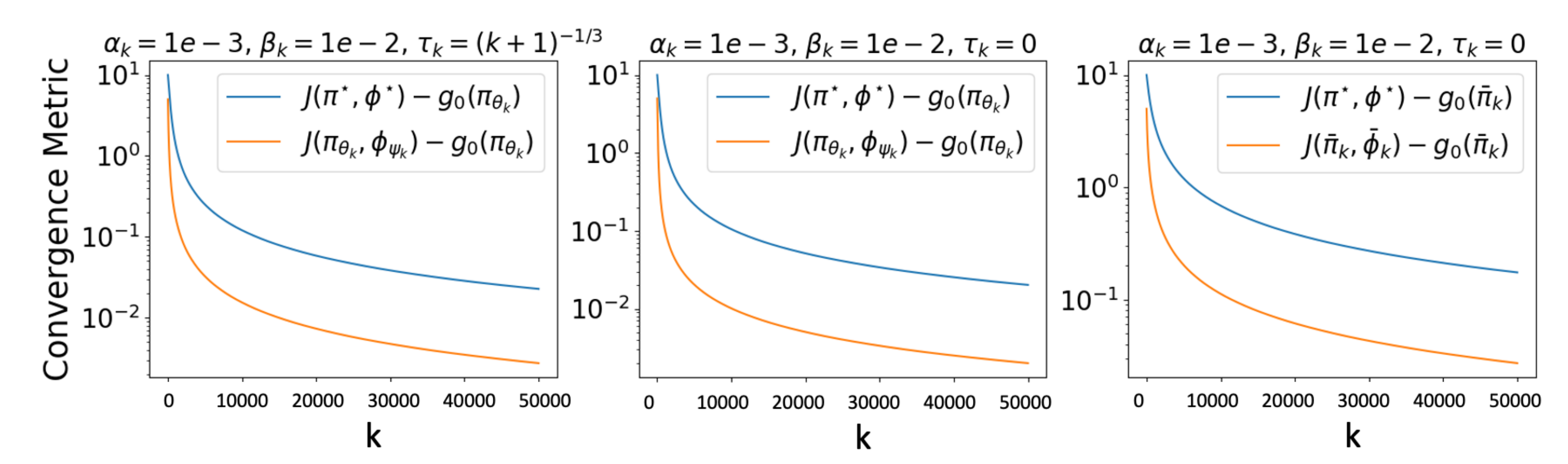}
  \caption{Convergence of GDA for a Deterministic Markov game}
  \label{fig:MG_deterministic}
\end{figure}

%% file: Conclusion.tex
\section{Conclusion \& Future Work}

In this paper, we present the finite-time analysis of two GDA algorithms that provably find the Nash equilibrium of a Markov game with the help of a structured entropy regularization. Future directions of this work include formalizing the link between Assumption~\ref{assump:NE_completelymixed} and completely mixed Markov games, investigating the possibility of relaxing this assumption, and characterizing the convergence of the stochastic GDA algorithm where the players do not have knowledge of the environment dynamics and can only take samples to estimate the gradients.

%% file: Proof_Theorem.tex
\section{Proof of Theorems and Corollaries}\label{sec:proof:thm}

We frequently use the following inequalities which hold for all $\tau\geq 0$, $\pi\in\Delta_{\Scal}^{\Acal}$, and $\phi\in\Delta_{\Scal}^{\Bcal}$,
\begin{align*}
    J_{\tau}(\pi,\phi_{\tau}(\pi))\leq J_{\tau}(\pi,\phi),\quad J_{\tau}(\pi_{\tau}(\phi),\phi)\geq J_{\tau}(\pi,\phi).
\end{align*}

We use $H(\cdot)$ to denote the entropy of a distribution. For example,
\begin{align}
    H(\pi(\cdot\mid s))=-\sum_{a}\pi(a\mid s)\log\pi(a\mid s),\quad H(\phi(\cdot\mid s))=-\sum_{b}\phi(b\mid s)\log\phi(b\mid s).\label{eq:H_def}
\end{align}

Due to the uniqueness of $\phi_{\tau}(\cdot)$, Danskin's Theorem guarantees that
$g_{\tau}(\pi_{\theta})$ defined in \eqref{eq:def_g} is differentiable with respect to $\theta$ \citep{bernhard1995theorem}
\begin{align}
    \nabla_{\theta} g_{\tau}(\pi_{\theta}) = \nabla_{\theta} J_{\tau}(\pi_{\theta},\phi),\quad\phi=\phi_{\tau}(\pi_{\theta}), \quad\forall \theta\in\mathbb{R}^{|\Scal|\times|\Acal|}.
    \label{eq:Danskin}
\end{align}

We also introduce a few lemmas that will be applied regularly in the rest of the paper.

\begin{lem}\label{lem:value_LipschitzGrad}
Let $L_V=\frac{8}{(1-\gamma)^3}$. The value function $J$ is $L_V$-Lipschitz continuous and has $L_V$-Lipschitz gradients, i.e. we have for all $\theta_1,\theta_2\in\mathbb{R}^{|\Scal|\times|\Acal|}$ and $\psi_1,\psi_2\in\mathbb{R}^{|\Scal|\times|\Bcal|}$
\begin{align*}
    &\|\nabla_{\theta}J(\pi_{\theta_1},\phi_{\psi_1})-\nabla_{\theta}J(\pi_{\theta_2},\phi_{\psi_2})\|\leq L_V(\|\theta_1-\theta_2\|+\|\psi_1-\psi_2\|),\\
    &\|\nabla_{\psi}J(\pi_{\theta_1},\phi_{\psi_1})-\nabla_{\psi}J(\pi_{\theta_2},\phi_{\psi_2})\|\leq L_V(\|\theta_1-\theta_2\|+\|\psi_1-\psi_2\|),\\
    &\|J(\pi_{\theta_1},\phi_{\psi_1})-J(\pi_{\theta_2},\phi_{\psi_2})\|\leq L_V(\|\theta_1-\theta_2\|+\|\psi_1-\psi_2\|).
\end{align*}
\end{lem}

\begin{lem}\label{lem:reg_LipschitzGrad}
Let $L_{\Hcal}=\frac{4+8\log|\Acal|}{(1-\gamma)^3}$. The regularization functions $\Hcal_{\pi}$ and $\Hcal_{\phi}$ are $L_{\Hcal}$-Lipschitz continuous and has $L_{\Hcal}$-Lipschitz gradients.
\end{lem}

Lemmas \ref{lem:value_LipschitzGrad} and \ref{lem:reg_LipschitzGrad} imply that $\forall\tau\geq0$, $\nabla_{\theta}J_{\tau}$ is Lipschitz continuous, i.e. for any $\theta_1,\theta_2\in\mathbb{R}^{|\Scal|\times|\Acal|}$, $\psi_1,\psi_2\in\mathbb{R}^{|\Scal|\times|\Bcal|}$
\begin{align}
    \|\nabla_{\theta}J_{\tau}(\pi_{\theta_1},\phi_{\psi_1})-\nabla_{\theta}J_{\tau}(\pi_{\theta_2},\phi_{\psi_2})\|
    &\leq \|\nabla_{\theta}J(\pi_{\theta_1},\phi_{\psi_1})-\nabla_{\theta}J(\pi_{\theta_2},\phi_{\psi_2})\|\notag\\
    &\hspace{20pt}+\tau\|\nabla_{\theta}\Hcal_{\pi}(s,\pi_{\theta_1},\phi_{\psi_1})-\nabla_{\theta}\Hcal_{\pi}(s,\pi_{\theta_2},\phi_{\psi_2})\|\notag\\
    &\hspace{20pt}+\tau\|\nabla_{\theta}\Hcal_{\phi}(s,\pi_{\theta_1},\phi_{\psi_1})-\nabla_{\theta}\Hcal_{\phi}(s,\pi_{\theta_2},\phi_{\psi_2})\|\notag\\
    &\leq (L_V+2\tau L_{\Hcal})(\|\theta_1-\theta_2\|+\|\psi_1-\psi_2\|).\label{eq:Lipschitz_Jtau}
\end{align}

\begin{lem}\label{lem:stepsize_diff}
For any $0\leq a \leq 1$ and integer $k>0$, we have
\begin{align*}
    \frac{1}{(k+h)^a}-\frac{1}{(k+1+h)^a}\leq\frac{8}{3(k+h)^{a+1}}.
\end{align*}
\end{lem}

\subsection{Proof of Theorem \ref{thm:main_constanttau}}

The definition of the constant $L$ and Eq.\eqref{eq:Lipschitz_Jtau} imply for any $\theta_1,\theta_2\in\mathbb{R}^{|\Scal|\times|\Acal|}$, $\psi_1,\psi_2\in\mathbb{R}^{|\Scal|\times|\Bcal|}$
\begin{align}
    \|\nabla_{\theta}J_{\tau}(\pi_{\theta_1},\phi_{\psi_1})-\nabla_{\theta}J_{\tau}(\pi_{\theta_2},\phi_{\psi_2})\|\leq L(\|\theta_1-\theta_2\|+\|\psi_1-\psi_2\|).\label{thm:main_constanttau:proof_eq2}
\end{align}

We will use an induction argument to prove the convergence of $3\delta_k^{\pi}+\delta_k^{\phi}$. 
The base case is $3\delta_0^{\pi}+\delta_0^{\phi}\leq3\delta_0^{\pi}+\delta_0^{\phi}$, which obviously holds.
Now, suppose
\begin{align}
    3\delta_k^{\pi}+\delta_k^{\phi}\leq (1-\frac{\alpha (1-\gamma) \tau \rho_{\min}^2 c^2}{32|\Scal|})^{k}(3\delta_0^{\pi}+\delta_0^{\phi})\label{thm:main_constanttau:proof_eq1}
\end{align}
holds. We aim to show
\begin{align*}
    3\delta_{k+1}^{\pi}+\delta_{k+1}^{\phi}\leq (1-\frac{\alpha (1-\gamma) \tau \rho_{\min}^2 c^2}{32|\Scal|})^{k+1}(3\delta_0^{\pi}+\delta_0^{\phi}).
\end{align*}

We introduce the following technical lemmas.

\begin{lem}\label{lem:min_policy_iterates_constanttau}
Suppose \eqref{thm:main_constanttau:proof_eq1} holds. Then, we have
\begin{align}
    -\left(\min_{s,a}\pi_{\theta_k}(a\mid s)\right)^2 &\leq-\frac{3c^2}{8},
    \label{lem:min_policy_iterates_constanttau:eq1}\\
    -\left(\min_{s,b}\phi_{\psi_k}(b\mid s)\right)^2&\leq -\frac{3c^2}{8}. \label{lem:min_policy_iterates_constanttau:eq2}
\end{align}
\end{lem}

\begin{lem}\label{lem:g_smooth_constanttau}
Suppose \eqref{thm:main_constanttau:proof_eq1} holds.
Under Assumption \ref{assump:positive_rho} and the step size $\alpha_k\leq(L+\frac{2\sqrt{|\Scal|}L^2}{\sqrt{(1-\gamma)\rho_{\min}}\tau c})^{-1}$, we have
\begin{align*}
    &g_{\tau}(\theta_k)-g_{\tau}(\theta_{k+1})\notag\\
    &=J_{\tau}(\pi_{\theta_k},\phi_{\tau}(\pi_{\theta_k}))-J_{\tau}(\pi_{\theta_{k+1}},\phi_{\tau}(\pi_{\theta_{k+1}}))\notag\\
    &= \frac{\alpha_k}{2}\left(\|\nabla_{\theta} J_{\tau}(\pi_{\theta_{k}},\phi_{\tau}(\pi_{\theta_k}))-\nabla_{\theta}J_{\tau}(\pi_{\theta_k},\phi_{\psi_k})\|^2-\|\nabla_{\theta} J_{\tau}(\pi_{\theta_{k}},\phi_{\tau}(\pi_{\theta_k}))\|^2\right).
\end{align*}
\end{lem}

By the lemma above, we have
\begin{align}
    &\delta^{\pi}_{k+1}-\delta^{\pi}_{k}\notag\\
    &=J_{\tau}(\pi_{\theta_k},\phi_{\tau}(\pi_{\theta_k}))-J_{\tau}(\pi_{\theta_{k+1}},\phi_{\tau}(\pi_{\theta_{k+1}}))\notag\\
    &\leq  \frac{\alpha_k}{2}\left(\|\nabla_{\theta} J_{\tau}(\pi_{\theta_{k}},\phi_{\tau}(\pi_{\theta_k}))-\nabla_{\theta}J_{\tau}(\pi_{\theta_k},\phi_{\psi_k})\|^2-\|\nabla_{\theta} J_{\tau}(\pi_{\theta_{k}},\phi_{\tau}(\pi_{\theta_k}))\|^2\right).
    \label{thm:pg_constanttau:eq1}
\end{align}

Similarly, we consider the decay of $\delta^{\phi}_{k}$.
\begin{align}
    \delta^{\phi}_{k+1}-\delta^{\phi}_{k}&=J_{\tau}(\pi_{\theta_{k+1}},\phi_{\psi_{k+1}})-g_{\tau}(\pi_{\theta_{k+1}})-J_{\tau}(\pi_{\theta_k},\phi_{\psi_k})+g_{\tau}(\pi_{\theta_k})\notag\\
    &=\left(J_{\tau}(\pi_{\theta_{k+1}},\phi_{\psi_{k+1}})-J_{\tau}(\pi_{\theta_{k+1}},\phi_{\psi_{k}})\right)\notag\\
    &\hspace{20pt}+\left(J_{\tau}(\pi_{\theta_{k+1}},\phi_{\psi_{k}})-J_{\tau}(\pi_{\theta_k},\phi_{\psi_k})\right)+\left(g_{\tau}(\pi_{\theta_k})-g_{\tau}(\pi_{\theta_{k+1}})\right).
    \label{thm:pg_constanttau:eq2}
\end{align}

Using the $L$-smoothness of the value function derived in \eqref{thm:main_constanttau:proof_eq2}
\begin{align*}
    &J_{\tau}(\pi_{\theta_{k+1}},\phi_{\psi_{k+1}})-J_{\tau}(\pi_{\theta_{k+1}},\phi_{\psi_{k}})\\
    &\leq\langle \nabla_{\psi} J_{\tau}(\pi_{\theta_{k+1}},\phi_{\psi_{k}}),\psi_{k+1}-\psi_{k}\rangle+\frac{L}{2}\|\psi_{k+1}-\psi_k\|^2\notag\\
    &= -\beta_k\|\nabla_{\psi} J_{\tau}(\pi_{\theta_{k+1}},\phi_{\psi_{k}})\|^2+\frac{L\beta_k^2}{2}\|\nabla_{\psi} J_{\tau}(\pi_{\theta_{k+1}},\phi_{\psi_{k}})\|^2\notag\\
    &\leq -\frac{\beta_k}{2}\|\nabla_{\psi} J_{\tau}(\pi_{\theta_{k+1}},\phi_{\psi_{k}})\|^2\notag\\
    &\leq -\frac{(1-\gamma)\beta_k\tau\rho_{\min}^2}{|\Scal|}\left(\min_{s,b}\phi_{\psi_k}(b\mid s)\right)^2\left(J_{\tau}(\pi_{\theta_k},\phi_{\psi_k})-J_{\tau}(\pi_{\theta_k},\phi_{\tau}(\pi_{\theta_k}))\right)\notag\\
    &= -\frac{(1-\gamma)\beta_k\tau\rho_{\min}^2}{|\Scal|}\left(\min_{s,b}\phi_{\psi_k}(b\mid s)\right)^2\delta_k^{\phi},
\end{align*}
where the second inequality uses $\beta_k\leq\frac{1}{L}$ and the third inequality follows from Lemma \ref{lem:nonuniform_PL} and the fact that $d_{\rho}^{\pi,\phi}(s)\leq1$ for all $s\in\Scal$ and policies $\pi,\phi$.

Using Eq.~\eqref{lem:min_policy_iterates_constanttau:eq2} of Lemma \ref{lem:min_policy_iterates_constanttau} to further simplify this inequality,
\begin{align}
    J_{\tau}(\pi_{\theta_{k+1}},\phi_{\psi_{k+1}})-J_{\tau}(\pi_{\theta_{k+1}},\phi_{\psi_{k}})&\leq -\frac{3(1-\gamma)\beta_k\tau\rho_{\min}^2 c^2}{8|\Scal|}\delta_k^{\phi}\label{thm:pg_constanttau:eq3}.
\end{align}

For the second term of \eqref{thm:pg_constanttau:eq2}, we have from the $L$-smoothness of the value function derived in \eqref{thm:main_constanttau:proof_eq2}
\begin{align}
    J_{\tau}(\pi_{\theta_{k+1}},\phi_{\psi_{k}})-J_{\tau}(\pi_{\theta_k},\phi_{\psi_k})&\leq\langle \nabla_{\theta} J_{\tau}(\pi_{\theta_k},\phi_{\psi_k}),\theta_{k+1}-\theta_k\rangle+\frac{L}{2}\|\theta_{k+1}-\theta_k\|^2\notag\\
    &= \alpha_k\|\nabla_{\theta} J_{\tau}(\pi_{\theta_k},\phi_{\psi_k})\|^2+\frac{L\alpha_k^2}{2}\|\nabla_{\theta} J_{\tau}(\pi_{\theta_k},\phi_{\psi_k})\|^2\notag\\
    &\leq \frac{3\alpha_k}{2}\|\nabla_{\theta} J_{\tau}(\pi_{\theta_k},\phi_{\psi_k})\|^2,
    \label{thm:pg_constanttau:eq4}
\end{align}
where in the last inequality we use $\alpha_k L\leq 1$.

Similarly to \eqref{thm:pg_constanttau:eq1}, the last term of \eqref{thm:pg_constanttau:eq2} is bounded as
\begin{align}
    g_{\tau}(\pi_{\theta_k})\hspace{-2pt}-\hspace{-2pt}g_{\tau}(\pi_{\theta_{k+1}}) &=g_{\tau}(\pi_{\theta_k})-g_{\tau}(\pi_{\theta_{k+1}})+g_{\tau}(\pi_{\theta_{k+1}})-g_{\tau}(\pi_{\theta_{k+1}})\notag\\
    &\leq \frac{\alpha_k}{2}\hspace{-2pt}\left(\|\nabla_{\theta} J_{\tau}(\pi_{\theta_{k}},\phi_{\tau}(\pi_{\theta_k}))\hspace{-2pt}-\hspace{-2pt}\nabla_{\theta}J_{\tau}(\pi_{\theta_k},\phi_{\psi_k})\|^2\hspace{-2pt}-\hspace{-2pt}\|\nabla_{\theta} J_{\tau}(\pi_{\theta_{k}},\phi_{\tau}(\pi_{\theta_k}))\|^2\right)
    \label{thm:pg_constanttau:eq5}
\end{align}

Using \eqref{thm:pg_constanttau:eq3}-\eqref{thm:pg_constanttau:eq5} in \eqref{thm:pg_constanttau:eq2}, we have
\begin{align}
    \delta^{\phi}_{k+1} &=\left(J_{\tau}(\pi_{\theta_{k+1}},\phi_{\psi_{k+1}})-J_{\tau}(\pi_{\theta_{k+1}},\phi_{\psi_{k}})\right)\notag\\
    &\hspace{20pt}+\left(J_{\tau}(\pi_{\theta_{k+1}},\phi_{\psi_{k}})-J_{\tau}(\pi_{\theta_k},\phi_{\psi_k})\right)+\left(g_{\tau}(\pi_{\theta_k})-g_{\tau}(\pi_{\theta_{k+1}})\right)\notag\\
    &\leq (1-\frac{3(1-\gamma)\beta_k\tau \rho_{\min}^2 c^2}{8|\Scal|})\delta^{\phi}_{k}+\frac{3\alpha_k}{2}\|\nabla_{\theta} J_{\tau}(\pi_{\theta_k},\phi_{\psi_k})\|^2\notag\\
    &\hspace{20pt}+\frac{\alpha_k}{2}\left(\|\nabla_{\theta} J_{\tau}(\pi_{\theta_{k}},\phi_{\tau}(\pi_{\theta_k}))-\nabla_{\theta}J_{\tau}(\pi_{\theta_k},\phi_{\psi_k})\|^2-\|\nabla_{\theta} J_{\tau}(\pi_{\theta_{k}},\phi_{\tau}(\pi_{\theta_k}))\|^2\right).
    \label{thm:pg_constanttau:eq6}
\end{align}

Combining \eqref{thm:pg_constanttau:eq1} and \eqref{thm:pg_constanttau:eq6},
\begin{align*}
    3\delta^{\pi}_{k+1}\hspace{-2pt}+\hspace{-2pt}\delta^{\phi}_{k+1}\hspace{-2pt}&\leq 3\delta^{\pi}_{k}\hspace{-2pt}+\hspace{-2pt}\frac{3\alpha_k}{2}\left(\|\nabla_{\theta} J_{\tau}(\pi_{\theta_{k}},\phi_{\tau}(\pi_{\theta_k}))\hspace{-2pt}-\hspace{-2pt}\nabla_{\theta}J_{\tau}(\pi_{\theta_k},\phi_{\psi_k})\hspace{-1pt}\|^2\hspace{-2pt}-\hspace{-2pt}\|\nabla_{\theta} J_{\tau}(\pi_{\theta_{k}},\phi_{\tau}(\pi_{\theta_k}))\|^2\right)\notag\\
    &\hspace{20pt}+(1-\frac{3(1-\gamma)\beta_k\tau \rho_{\min}^2 c^2}{8|\Scal|})\delta^{\phi}_{k}+2\delta_k^{\phi})\delta_k^{\phi}+\frac{3\alpha_k}{2}\|\nabla_{\theta} J_{\tau}(\pi_{\theta_k},\phi_{\psi_k})\|^2\notag\\
    &\hspace{20pt}+\frac{\alpha_k}{2}\left(\|\nabla_{\theta} J_{\tau}(\pi_{\theta_{k}},\phi_{\tau}(\pi_{\theta_k}))-\nabla_{\theta}J_{\tau}(\pi_{\theta_k},\phi_{\psi_k})\|^2-\|\nabla_{\theta} J_{\tau}(\pi_{\theta_{k}},\phi_{\tau}(\pi_{\theta_k}))\|^2\right)\notag\\
    &\leq 3\delta^{\pi}_{k}+(1-\frac{3(1-\gamma)\beta_k\tau \rho_{\min}^2 c^2}{8|\Scal|})\delta^{\phi}_{k}+\frac{3\alpha_k}{2}\|\nabla_{\theta} J_{\tau}(\pi_{\theta_k},\phi_{\psi_k})\|^2\notag\\
    &\hspace{20pt}+2\alpha_k\|\nabla_{\theta} J_{\tau}(\pi_{\theta_{k}},\phi_{\tau}(\pi_{\theta_k}))\hspace{-2pt}-\hspace{-2pt}\nabla_{\theta}J_{\tau}(\pi_{\theta_k},\phi_{\psi_k})\|^2\hspace{-2pt}-\hspace{-2pt}2\alpha_k\|\nabla_{\theta} J_{\tau}(\pi_{\theta_{k}},\phi_{\tau}(\pi_{\theta_k}))\|^2.
\end{align*}
Simplifying this inequality with 
\begin{align*}
    \|\nabla_{\theta} J_{\tau}(\pi_{\theta_k},\phi_{\psi_k})\|^2&=\|\nabla_{\theta} J_{\tau}(\pi_{\theta_{k}},\phi_{\tau}(\pi_{\theta_k}))-\left(\nabla_{\theta} J_{\tau}(\pi_{\theta_{k}},\phi_{\tau}(\pi_{\theta_k}))-\nabla_{\theta}J_{\tau}(\pi_{\theta_k},\phi_{\psi_k})\right)\|^2\\
    &\leq \|\nabla_{\theta} J_{\tau}(\pi_{\theta_{k}},\phi_{\tau}(\pi_{\theta_k}))\|^2+\|\nabla_{\theta} J_{\tau}(\pi_{\theta_{k}},\phi_{\tau}(\pi_{\theta_k}))-\nabla_{\theta}J_{\tau}(\pi_{\theta_k},\phi_{\psi_k})\|^2\\
    &\hspace{20pt}+2\langle\nabla_{\theta} J_{\tau}(\pi_{\theta_{k}},\phi_{\tau}(\pi_{\theta_k})),\nabla_{\theta} J_{\tau}(\pi_{\theta_{k}},\phi_{\tau}(\pi_{\theta_k}))-\nabla_{\theta}J_{\tau}(\pi_{\theta_k},\phi_{\psi_k})\rangle\\
    &\leq \frac{5}{4}\|\nabla_{\theta} J_{\tau}(\pi_{\theta_{k}},\phi_{\tau}(\pi_{\theta_k}))\|^2+5\|\nabla_{\theta} J_{\tau}(\pi_{\theta_{k}},\phi_{\tau}(\pi_{\theta_k}))-\nabla_{\theta}J_{\tau}(\pi_{\theta_k},\phi_{\psi_k})\|^2,
\end{align*}
we have
\begin{align}
    3\delta^{\pi}_{k+1}+\delta^{\phi}_{k+1}&\leq 3\delta^{\pi}_{k}+(1-\frac{3(1-\gamma)\beta_k\tau \rho_{\min}^2 c^2}{8|\Scal|})\delta^{\phi}_{k}-\frac{\alpha_k}{8}\|\nabla_{\theta} J_{\tau}(\pi_{\theta_{k}},\phi_{\tau}(\pi_{\theta_k}))\|^2\notag\\
    &\hspace{20pt}+\frac{19\alpha_k}{2}\|\nabla_{\theta} J_{\tau}(\pi_{\theta_{k}},\phi_{\tau}(\pi_{\theta_k}))-\nabla_{\theta}J_{\tau}(\pi_{\theta_k},\phi_{\psi_k})\|^2.
    \label{thm:pg_constanttau:eq7}
\end{align}

Using Lemma \ref{lem:nonuniform_PL} to bound $-\|\nabla_{\theta} J_{\tau}(\pi_{\theta_{k}},\phi_{\tau}(\pi_{\theta_k}))\|^2$, 
\begin{align}
    &-\|\nabla_{\theta} J_{\tau}(\pi_{\theta_{k}},\phi_{\tau}(\pi_{\theta_k}))\|^2\notag  \\
    &\leq-\frac{2(1-\gamma)\tau\rho_{\min}^2}{|\Scal|}\left(\min_{s,a}\pi_{\theta_k}(a\mid s)\right)^2\left(J_{\tau}(\pi_{\tau}(\phi_{\tau}(\pi_{\theta_k})),\phi_{\tau}(\pi_{\theta_k}))-J_{\tau}(\pi_{\theta_k},\phi_{\tau}(\pi_{\theta_k}))\right)\notag\\
    &\leq-\frac{2(1-\gamma)\tau\rho_{\min}^2}{|\Scal|}\left(\min_{s,a}\pi_{\theta_k}(a\mid s)\right)^2\left(J_{\tau}(\pi_{\tau}^{\star},\phi_{\tau}^{\star})-J_{\tau}(\pi_{\theta_k},\phi_{\tau}(\pi_{\theta_k}))\right),
    \label{thm:pg_constanttau:eq7.5}
\end{align}
where the second inequality follows from
\begin{align*}
    J_{\tau}(\pi_{\tau}(\phi_{\tau}(\pi_{\theta_k})),\phi_{\tau}(\pi_{\theta_k}))=\max_{\pi}J_{\tau}(\pi,\phi_{\tau}(\pi_{\theta_k}))\geq \max_{\pi}\min_{\phi}J_{\tau}(\pi,\phi)=J_{\tau}(\pi_{\tau}^{\star},\phi_{\tau}^{\star}).
\end{align*}

From Lemma \ref{lem:min_policy_iterates_constanttau} Eq.~\eqref{lem:min_policy_iterates_constanttau:eq1}, $-\left(\min_{s,a}\pi_{\theta_k}(a\mid s)\right)^2\leq-\frac{3c^2}{8}$, which further simplifies \eqref{thm:pg_constanttau:eq7.5}
\begin{align*}
    -\|\nabla_{\theta} J_{\tau}(\pi_{\theta_{k}},\phi_{\tau}(\pi_{\theta_k}))\|^2&\leq-\frac{2(1-\gamma)\tau\rho_{\min}^2}{|\Scal|}\hspace{-2pt}\left(\min_{s,a}\pi_{\theta_k}(a\mid s)\right)^2\hspace{-2pt}\left(J_{\tau}(\pi_{\tau}^{\star},\phi_{\tau}^{\star})\hspace{-2pt}-\hspace{-2pt}J_{\tau}(\pi_{\theta_k},\phi_{\tau}(\pi_{\theta_k}))\right)\\
    &=-\frac{2(1-\gamma)\tau\rho_{\min}^2}{|\Scal|}\left(\min_{s,a}\pi_{\theta_k}(a\mid s)\right)^2\delta_k^{\pi}\leq-\frac{3(1-\gamma)\tau\rho_{\min}^2 c^2}{4|\Scal|}\delta_k^{\pi}.
\end{align*}

For $\|\nabla_{\theta} J_{\tau}(\pi_{\theta_{k}},\phi_{\tau}(\pi_{\theta_k}))-\nabla_{\theta}J_{\tau}(\pi_{\theta_k},\phi_{\psi_k})\|^2$, we have from the $L$-smoothness of the value function derived in \eqref{thm:main_constanttau:proof_eq2}
\begin{align*}
    \|\nabla_{\theta} J_{\tau}(\pi_{\theta_{k}},\phi_{\tau}(\pi_{\theta_k}))-\nabla_{\theta}J_{\tau}(\pi_{\theta_k},\phi_{\psi_k})\|^2&\leq L^2\|\phi_{\tau}(\pi_{\theta_k})-\phi_{\psi_k}\|^2\notag\\
    &\leq \frac{2\log(2) L^2}{\tau\rho_{\min}}\left(J_{\tau}(\pi_{\theta_k},\phi_{\psi_k})-J_{\tau}(\pi_{\theta_{k}},\phi_{\tau}(\pi_{\theta_k}))\right)\notag\\
    &=\frac{2\log(2) L^2}{\tau\rho_{\min}}\delta_k^{\phi}
\end{align*}

Using the bound on $-\|\nabla_{\theta} J_{\tau}(\pi_{\theta_{k}},\phi_{\tau}(\pi_{\theta_k}))\|^2$ and $\|\nabla_{\theta} J_{\tau}(\pi_{\theta_{k}},\phi_{\tau}(\pi_{\theta_k}))-\nabla_{\theta}J_{\tau}(\pi_{\theta_k},\phi_{\psi_k})\|^2$ in \eqref{thm:pg_constanttau:eq7},
\begin{align*}
    3\delta^{\pi}_{k+1}+\delta^{\phi}_{k+1}&\leq 3\delta^{\pi}_{k}+(1-\frac{3(1-\gamma)\beta_k\tau\rho_{\min}^2 c^2}{8|\Scal|})\delta^{\phi}_{k}-\frac{\alpha_k}{8}\|\nabla_{\theta} J_{\tau}(\pi_{\theta_{k}},\phi_{\tau}(\pi_{\theta_k}))\|^2\notag\\
    &\hspace{20pt}+\frac{19\alpha_k}{2}\|\nabla_{\theta} J_{\tau}(\pi_{\theta_{k}},\phi_{\tau}(\pi_{\theta_k}))-\nabla_{\theta}J_{\tau}(\pi_{\theta_k},\phi_{\psi_k})\|^2\notag\\
    &\leq 3\delta^{\pi}_{k}+(1\hspace{-1pt}-\hspace{-1pt}\frac{3(1-\gamma)\beta_k\tau\rho_{\min}^2 c^2}{8|\Scal|})\delta^{\phi}_{k}\hspace{-1pt}-\hspace{-1pt}\frac{3\alpha_k(1-\gamma)\tau\rho_{\min}^2 c^2}{32|\Scal|}\delta_k^{\pi}\hspace{-1pt}+\hspace{-1pt}\frac{19\log(2) L^2 \alpha_k}{\tau\rho_{\min}}\delta_k^{\phi}\notag\\
    &=3(1-\frac{\alpha_k (1-\gamma) \tau \rho_{\min}^2 c^2}{32|\Scal|})\delta_k^{\pi}+(1-\frac{3(1-\gamma)\beta_k\tau\rho_{\min}^2 c^2}{8|\Scal|}+\frac{19\log(2) L^2 \alpha_k}{\tau\rho_{\min}})\delta^{\phi}_{k}.
\end{align*}

With the step sizes $\alpha_k=\alpha$, $\beta_k=\beta$ such that $\frac{\alpha}{\beta}\leq\min\{\frac{(1-\gamma)\tau^2\rho_{\min}^3 c^2}{152|\Scal|\log(2)L^2},8\}$, we can simplify the inequality above
\begin{align*}
    3\delta^{\pi}_{k+1}+\delta^{\phi}_{k+1}&\leq 3(1-\frac{\alpha_k (1-\gamma) \tau \rho_{\min}^2 c^2}{32|\Scal|})\delta_k^{\pi}+(1-\frac{3(1-\gamma)\beta_k\tau\rho_{\min}^2 c^2}{8|\Scal|}+\frac{19\log(2) L^2 \alpha_k}{\tau\rho_{\min}})\delta^{\phi}_{k}\\
    &\leq 3(1-\frac{\alpha (1-\gamma) \tau \rho_{\min}^2 c^2}{32|\Scal|})\delta_k^{\pi}+(1-\frac{(1-\gamma)\beta\tau\rho_{\min}^2 c^2}{4|\Scal|})\delta^{\phi}_{k}\\
    &\leq (1-\frac{\alpha (1-\gamma) \tau \rho_{\min}^2 c^2}{32|\Scal|})(3\delta_k^{\pi}+\delta_k^{\phi})\\
    &\leq (1-\frac{\alpha (1-\gamma) \tau \rho_{\min}^2 c^2}{32|\Scal|})^{k+1}(3\delta_0^{\pi}+\delta_0^{\phi}).
\end{align*}

\qed

\subsection{Proof of Corollary \ref{cor:GDA_piecewiseconstanttau}}

As a result of Lemma \ref{lem:V_tau_diff}, it is easy to verify
\begin{align}
    &(3\delta_{t+1,0}^{\pi}+\delta_{t+1,0}^{\phi})-(3\delta_{t,K_t}^{\pi}+\delta_{t,K_t}^{\phi})\notag\\
    &=(3J_{\tau_{t+1}}(\pi_{\tau_{t+1}}^{\star},\phi_{\tau_{t+1}}^{\star}) - 3J_{\tau_{t+1}}(\pi_{\theta_{t+1,0}},\phi_{\tau_{t+1}}(\pi_{\theta_{t+1,0}}))\notag\\
    &\hspace{20pt}+J_{\tau_{t+1}}(\pi_{\theta_{t+1,0}},\phi_{\psi_{t+1,0}})-J_{\tau_{t+1}}(\pi_{\theta_{t+1,0}},\phi_{\tau_{t+1}}(\pi_{\theta_{t+1,0}})))\notag\\
    &\hspace{20pt}-(3J_{\tau_{t}}(\pi_{\tau_{t}}^{\star},\phi_{\tau_t}^{\star}) - 3J_{\tau_{t}}(\pi_{\theta_{t,K_t}},\phi_{\tau_{t}}(\pi_{\theta_{t,K_t}}))\notag\\
    &\hspace{20pt}+J_{\tau_{t}}(\pi_{\theta_{t,K_t}},\phi_{\psi_{t,K_t}})-J_{\tau_{t}}(\pi_{\theta_{t,K_t}},\phi_{\tau_{t}}(\pi_{\theta_{t,K_t}}))\notag\\
    &=(3J_{\tau_{t+1}}(\pi_{\tau_{t+1}}^{\star},\phi_{\tau_{t+1}}^{\star}) - 3J_{\tau_{t+1}}(\pi_{\theta_{t+1,0}},\phi_{\tau_{t+1}}(\pi_{\theta_{t+1,0}}))\notag\\
    &\hspace{20pt}+J_{\tau_{t+1}}(\pi_{\theta_{t+1,0}},\phi_{\psi_{t+1,0}})-J_{\tau_{t+1}}(\pi_{\theta_{t+1,0}},\phi_{\tau_{t+1}}(\pi_{\theta_{t+1,0}})))\notag\\
    &\hspace{20pt}-(3J_{\tau_{t}}(\pi_{\tau_{t}}^{\star},\phi_{\tau_t}^{\star}) - 3J_{\tau_{t}}(\pi_{\theta_{t+1,0}},\phi_{\tau_{t}}(\pi_{\theta_{t+1,0}}))\notag\\
    &\hspace{20pt}+J_{\tau_{t}}(\pi_{\theta_{t+1,0}},\phi_{\psi_{t+1,0}})-J_{\tau_{t}}(\pi_{\theta_{t+1,0}},\phi_{\tau_{t}}(\pi_{\theta_{t+1,0}})))\notag\\
    &=3(J_{\tau_{t+1}}(\pi_{\tau_{t+1}}^{\star},\phi_{\tau_{t+1}}^{\star})-J_{\tau_{t}}(\pi_{\tau_{t}}^{\star},\phi_{\tau_t}^{\star}))\notag\\
    &\hspace{20pt}-4(J_{\tau_{t+1}}(\pi_{\theta_{t+1,0}},\phi_{\tau_{t+1}}(\pi_{\theta_{t+1,0}}))-J_{\tau_{t}}(\pi_{\theta_{t+1,0}},\phi_{\tau_{t}}(\pi_{\theta_{t+1,0}})))\notag\\
    &\hspace{20pt}+(J_{\tau_{t+1}}(\pi_{\theta_{t+1,0}},\phi_{\psi_{t+1,0}})-J_{\tau_{t}}(\pi_{\theta_{t+1,0}},\phi_{\psi_{t+1,0}}))\notag\\
    &\leq L_{\delta}(\tau_t-\tau_{t+1}).
    \label{cor:GDA_piecewiseconstanttau:proof_eq1}
\end{align}

We can choose $\tau_0$ large enough that
\begin{align*}
    3\delta^{\pi}_{0,0}+\delta^{\phi}_{0,0}\leq C_1\tau_0
\end{align*}
holds.
For any $t\geq 0$, if we run the inner loop for $K_t$ iterations such that 
\[3\delta_{t,K_t}^{\pi}+\delta_{t,K_t}^{\phi}\leq\frac{1}{2}(3\delta_{t,0}^{\pi}+\delta_{t,0}^{\phi})\leq\frac{C_1\tau_t}{2},\]
then we have
\begin{align*}
    3\delta_{t+1,0}^{\pi}+\delta_{t+1,0}^{\phi}&\leq 3\delta_{t,K_t}^{\pi}+\delta_{t,K_t}^{\phi}+L_{\delta}(\tau_t-\tau_{t+1})\leq \frac{C_1\tau_t}{2}+L_{\delta}(\tau_t-\tau_{t+1})\\
    &=\frac{(C_1+L_{\delta})C_1}{C_1+2L_{\delta}}\tau_{t+1}+\frac{C_1 L_{\delta}}{C_1+2L_{\delta}}\tau_{t+1}= C_1\tau_{t+1},
\end{align*}
where the first equality plugs in $\tau_{t}=\frac{2C_1+2L_{\delta}}{C_1+2L_{\delta}}\tau_{t+1}$. This means that the initial condition \eqref{thm:main_constanttau:eq1} is observed at the beginning of the every outer loop iteration.

Applying the inequality recursively,
\begin{align*}
    3\delta_{T,0}^{\pi}+\delta_{T,0}^{\phi}&\leq C_1\tau_T.
\end{align*}

With an argument similar to the one in \eqref{cor:GDA_piecewiseconstanttau:proof_eq1}, we can show
\begin{align*}
    &(3(J(\pi^{\star},\phi^{\star}) - J(\pi_{\theta_{T,0}},\phi_{0}(\pi_{\theta_{T,0}})))\\
    &\hspace{100pt}+(J(\pi_{\theta_{T,0}},\phi_{\psi_{T,0}})-J(\pi_{\theta_{T,0}},\phi_{0}(\pi_{\theta_{T,0}}))))-(3\delta_{T,0}^{\pi}+\delta_{T,0}^{\phi})\leq L_{\delta}\tau_T.
\end{align*}
In order to achieve \eqref{cor:GDA_piecewiseconstanttau:eq1}, it suffices to guarantee $3\delta_{T,0}^{\pi}+\delta_{T,0}^{\phi}+L_{\delta}\tau_T\leq \epsilon$, or $(C_1+L_{\delta})\tau_T\leq\epsilon$.
This implies that we need $\tau_T=\Ocal(\epsilon)$, or equivalently, $T=\Ocal(\log(\epsilon^{-1}))$ since $\tau_T=\left(\frac{C_1+2L_{\delta}}{2C_1+2L_{\delta}}\right)^T \tau_0$.

Ultimately we are interested in bounding $\sum_{t=0}^{T}K_t$. Note that $K_t$ needs to be at most
\begin{align*}
    K_t\leq\lceil\frac{\log(\frac{1}{2})}{\log(1-\frac{\alpha_t (1-\gamma) \tau_t \rho_{\min}^2 c^2}{32|\Scal|})}\rceil.
\end{align*}

To apply Theorem \ref{thm:main_constanttau}, we need to select the step sizes that satisfy the required condition. Since $\{\tau_t\}$ is a decaying sequence, the smoothness constant $L=3L_{\Hcal}\max\{\tau_0,1\}$ is valid across all outer loop iterations $t$.

We use $L_t=3 L_{\Hcal}\max\{\tau_t,1\}$ to denote the smoothness constant of the regularized value function in outer loop iteration $t$ and use $T_1$ to denote the index of the outer loop iteration such that $\tau_{T_1}\geq 1$ and $\tau_{T_1+1}<1$. Note that $T_1$ is an absolute constant that only depends on the structure of the Markov game. From iterations $t=0$ to $t=T_1$, the smoothness constant is proportional to regularization weight $L_t=3L_{\Hcal}\max\{\tau_t,1\}=3L_{\Hcal}\tau_t$. We need to choose $\alpha_t,\beta_t$ such that
\begin{align*}
    &\beta_t\leq\frac{1}{L_t}=\frac{1}{3L_{\Hcal}\tau_t},\quad\frac{\alpha_t}{\beta_t}\leq\min\{\frac{ (1-\gamma)\rho_{\min}^3 c^2 \tau_t^2}{152\log(2)|\Scal|L_t^2},8\}=\min\{\frac{ (1-\gamma)\rho_{\min}^3 c^2}{1368\log(2)|\Scal|L_{\Hcal}^2},8\},\\
    &\alpha_t\leq\min\{(L_t+\frac{2\sqrt{|\Scal|}L_t^2}{\sqrt{(1-\gamma)\rho_{\min}}\tau_t c})^{-1},\frac{16|\Scal|}{(1-\gamma)\rho_{\min}^2 c^2 \tau_t}\}\\
    &\hspace{20pt}=\min\{(3L_{\Hcal}\tau_t+\frac{18\sqrt{|\Scal|}L_{\Hcal}^2\tau_t}{\sqrt{(1-\gamma)\rho_{\min}} c})^{-1},\frac{16|\Scal|}{(1-\gamma)\rho_{\min}^2 c^2 \tau_t}\}.
\end{align*}
Then it is obvious that we can choose $\alpha_t=\Ocal(\tau_t^{-1})$, implying $\alpha_t\tau_t=\Ocal(1)$. Therefore, for all $0\leq t\leq T_1$,
\begin{align}
    K_t\leq \lceil\frac{\log(\frac{1}{2})}{\log(1-\frac{\alpha_t (1-\gamma) \tau_t \rho_{\min}^2 c^2}{32|\Scal|})}\rceil=\Ocal(1).
    \label{cor:GDA_piecewiseconstanttau:proof_eq2}
\end{align}

From iterations $t=T_1$ until $t=T$, the smoothness constant is $L_t=3L_{\Hcal}\max\{\tau_t,1\}=3L_{\Hcal}$.
Note that there is an upper and lower bound on $\beta_t$. In order for the upper bound to be no smaller than the lower bound, we need
\begin{align*}
    \frac{152\log(2)|\Scal|L^2 \alpha_t}{(1-\gamma)\rho_{\min}^3 c^2 \tau_t^2}\leq \frac{1}{L}.
\end{align*}
This means that we should choose $\alpha_t=\Ocal(\tau_t^2)$, implying $\alpha_t\tau_t=\Ocal(\tau_t^3)$. Plugging it in \eqref{cor:GDA_piecewiseconstanttau:proof_eq2},
\begin{align*}
    K_t=\lceil\frac{\log(\frac{1}{2})}{\log(1-\frac{\alpha_t (1-\gamma) \tau_t \rho_{\min}^2 c^2}{32|\Scal|})}\rceil=\Ocal(\frac{1}{\log(1-\tau_t^3)})\leq \Ocal(\tau_t^{-3}),
\end{align*}
where the last inequality follows from the fact that $1+x\leq \exp(x)$ for any scalar $x$.

Since $\tau_t=\tau_T(\frac{2C_1+2L_{\delta}}{C_1+2L_{\delta}})^{T-t}$,
\begin{align*}
    \sum_{t=0}^{T}K_t&=\sum_{t=0}^{T_1}K_t+\sum_{t=T_1}^{T}K_t\leq\sum_{t=0}^{T}\Ocal(\tau_t^{-3})=\Ocal(1)+\sum_{t=T_1}^{T}\Ocal(\tau_T^{-3}(\frac{2C_1+2L_{\delta}}{C_1+2L_{\delta}})^{-3(T-t)})\notag\\
    &\leq\Ocal(\tau_T^{-3}\sum_{t=0}^{T}(\frac{C_1+2L_{\delta}}{2C_1+2L_{\delta}})^{3(T-t)})=\Ocal(\tau_T^{-3}\sum_{t=0}^{T}(\frac{C_1+2L_{\delta}}{2C_1+2L_{\delta}})^{3t})\notag\\
    &\leq \Ocal(\tau_T^{-3}\frac{1}{1-(\frac{C_1+2L_{\delta}}{2C_1+2L_{\delta}})^3})=\Ocal(\tau_T^{-3}).
\end{align*}

Since $\tau_T=\Ocal(\epsilon)$, 
\begin{align*}
    \sum_{t=0}^{T}K_t\leq\Ocal(\tau_T^{-3})=\Ocal(\epsilon^{-3}).
\end{align*}

\qed

\subsection{Proof of Theorem \ref{thm:main}}

Define $L_0=L_{\Hcal}(2\tau_0+1)$. The exact conditions on the initial step sizes, regularization weight, and $h$ are
\begin{gather}
    \delta^{\pi}_{0}+\delta^{\phi}_{0}\leq \frac{C_1\tau_0}{h^{\frac{1}{3}}},\label{eq:stepsize_exist_eq1}\\
    \alpha_0=\frac{65536\log(2)(\log|\Acal|+\log|\Bcal|)+96(1-\gamma)\rho_{\min}c^2}{3(1-\gamma)^2 \rho_{\min}^3 c^4 \tau_0},\label{eq:stepsize_exist_eq2}\\
    \frac{\alpha_0}{h^{\frac{2}{3}}}\leq(2L_{\Hcal}+4L_{\Hcal}^2 C_2)\frac{\tau_0}{h^{\frac{1}{3}}}+(L_{\Hcal}+4L_{\Hcal}^2C_2)+\frac{L_{\Hcal}^2 C_2 h^{\frac{1}{3}}}{\tau_0},\label{eq:stepsize_exist_eq3}\\
    \beta_0\leq\frac{1}{L_0},\quad\frac{\alpha_0}{\beta_0}\leq\min\{\frac{ (1-\gamma)\tau_0^2\rho_{\min}^3 c^2}{152\log(2)|\Scal|L_0^2},1\}.
\end{gather}
In Remark~\ref{rem:stepsize_exist} at the end of this section, we show that there always exist $\alpha_0$, $\beta_0$, $\tau_0$, and $h$ that observe the conditions.

\eqref{eq:Lipschitz_Jtau} implies that for any $\theta_1,\theta_2\in\mathbb{R}^{|\Scal|\times|\Acal|}$, $\psi_1,\psi_2\in\mathbb{R}^{|\Scal|\times|\Bcal|}$, and $k\geq0$,
\begin{align}
    \|\nabla_{\theta}J_{\tau_k}(\pi_{\theta_1},\phi_{\psi_1})-\nabla_{\theta}J_{\tau_k}(\pi_{\theta_2},\phi_{\psi_2})\|
    &\leq (L_V+2\tau_k L_{\Hcal})(\|\theta_1-\theta_2\|+\|\psi_1-\psi_2\|)\notag\\
    &\leq L_0(\|\theta_1-\theta_2\|+\|\psi_1-\psi_2\|),\label{thm:main:proof_eq2}
\end{align}
where the last inequality follows from $\tau_k\leq\tau_0$.

\textbf{Convergence of $3\delta_k^{\pi}+\delta_k^{\phi}$:}

We will first use an induction argument to prove 
\[3\delta_k^{\pi}+\delta_k^{\phi}\leq\frac{\rho_{\min}\tau_0 c^2}{64\log(2)(k+h)^{1/3}},\quad \forall k\geq0.\]
The base case is $3\delta_0^{\pi}+\delta_0^{\phi}\leq\frac{\rho_{\min} c^2 \tau_0}{64\log(2)h^{\frac{1}{3}}}$, which holds by the initial condition.
Now, suppose
\begin{align}
    3\delta_k^{\pi}+\delta_k^{\phi}\leq\frac{\rho_{\min}\tau_0 c^2}{64\log(2)(k+h)^{1/3}}\label{thm:main:proof_eq1}
\end{align}
holds. We aim to show
\begin{align*}
    3\delta_{k+1}^{\pi}+\delta_{k+1}^{\phi}\leq\frac{\rho_{\min}\tau_0 c^2}{64\log(2)(k+1+h)^{1/3}}.
\end{align*}

We introduce the following technical lemmas.

\begin{lem}\label{lem:min_policy_iterates}
Suppose \eqref{thm:main:proof_eq1} holds. Then, we have
\begin{align}
    -\left(\min_{s,a}\pi_{\theta_k}(a\mid s)\right)^2 &\leq-\frac{3c^2}{8},
    \label{lem:min_policy_iterates:eq1}\\
    -\left(\min_{s,b}\phi_{\psi_k}(b\mid s)\right)^2&\leq -\frac{3c^2}{8}. \label{lem:min_policy_iterates:eq2}
\end{align}
\end{lem}

\begin{lem}\label{lem:g_smooth}
Suppose \eqref{thm:main:proof_eq1} holds.
Under Assumption \ref{assump:positive_rho} and \ref{assump:NE_completelymixed} and the step sizes of Theorem \ref{thm:main}, we have
\begin{align*}
    &g_{\tau_{k}}(\theta_k)-g_{\tau_{k}}(\theta_{k+1})\notag\\
    &=J_{\tau_{k}}(\pi_{\theta_k},\phi_{\tau_k}(\pi_{\theta_k}))-J_{\tau_{k}}(\pi_{\theta_{k+1}},\phi_{\tau_k}(\pi_{\theta_{k+1}}))\notag\\
    &\leq \frac{\alpha_k}{2}\left(\|\nabla_{\theta} J_{\tau_k}(\pi_{\theta_{k}},\phi_{\tau_k}(\pi_{\theta_k}))-\nabla_{\theta}J_{\tau_k}(\pi_{\theta_k},\phi_{\psi_k})\|^2-\|\nabla_{\theta} J_{\tau_k}(\pi_{\theta_{k}},\phi_{\tau_k}(\pi_{\theta_k}))\|^2\right).
\end{align*}
\end{lem}

We perform the following decomposition
\begin{align}
    &\delta^{\pi}_{k+1}-\delta^{\pi}_{k}\notag\\
    &=J_{\tau_{k}}(\pi_{\theta_k},\phi_{\tau_k}(\pi_{\theta_k}))-J_{\tau_{k+1}}(\pi_{\theta_{k+1}},\phi_{\tau_{k+1}}(\pi_{\theta_{k+1}}))+J_{\tau_{k+1}}(\pi_{\tau_{k+1}}^{\star},\phi_{\tau_{k+1}}^{\star})-J_{\tau_k}(\pi_{\tau_k}^{\star},\phi_{\tau_k}^{\star})\notag\\
    &=J_{\tau_{k}}(\pi_{\theta_k},\phi_{\tau_k}(\pi_{\theta_k}))-J_{\tau_{k}}(\pi_{\theta_{k+1}},\phi_{\tau_k}(\pi_{\theta_{k+1}}))\notag\\
    &\hspace{20pt}+J_{\tau_{k}}(\pi_{\theta_{k+1}},\phi_{\tau_k}(\pi_{\theta_{k+1}}))-J_{\tau_{k}}(\pi_{\theta_{k+1}},\phi_{\tau_{k+1}}(\pi_{\theta_{k+1}}))\notag\\
    &\hspace{20pt}+J_{\tau_{k}}(\pi_{\theta_{k+1}},\phi_{\tau_{k+1}}(\pi_{\theta_{k+1}}))-J_{\tau_{k+1}}(\pi_{\theta_{k+1}},\phi_{\tau_{k+1}}(\pi_{\theta_{k+1}}))\notag\\
    &\hspace{20pt}+J_{\tau_{k+1}}(\pi_{\tau_{k+1}}^{\star},\phi_{\tau_{k+1}}^{\star})-J_{\tau_k}(\pi_{\tau_k}^{\star},\phi_{\tau_k}^{\star})\notag\\
    &\leq J_{\tau_{k}}(\pi_{\theta_k},\phi_{\tau_k}(\pi_{\theta_k}))-J_{\tau_{k}}(\pi_{\theta_{k+1}},\phi_{\tau_k}(\pi_{\theta_{k+1}}))+\frac{\tau_k-\tau_{k+1}}{1-\gamma}\log|\Acal|+(\tau_k-\tau_{k+1})\log|\Bcal|\notag\\
    &\leq  \frac{\alpha_k}{2}\left(\|\nabla_{\theta} J_{\tau_k}(\pi_{\theta_{k}},\phi_{\tau_k}(\pi_{\theta_k}))-\nabla_{\theta}J_{\tau_k}(\pi_{\theta_k},\phi_{\psi_k})\|^2-\|\nabla_{\theta} J_{\tau_k}(\pi_{\theta_{k}},\phi_{\tau_k}(\pi_{\theta_k}))\|^2\right)\notag\\
    &\hspace{20pt}+\frac{\tau_k-\tau_{k+1}}{1-\gamma}(\log|\Acal|+\log|\Bcal|)\label{thm:main:proof_eq3}
\end{align}
where the first inequality comes from $J_{\tau_{k}}(\pi_{\theta_{k+1}},\phi_{\tau_k}(\pi_{\theta_{k+1}}))-J_{\tau_{k}}(\pi_{\theta_{k+1}},\phi_{\tau_{k+1}}(\pi_{\theta_{k+1}}))\leq 0$ by the definition of $\phi_{\tau}(\cdot)$ and the bound on $J_{\tau_{k}}(\pi_{\theta_{k+1}},\phi_{\tau_{k+1}}(\pi_{\theta_{k+1}}))-J_{\tau_{k+1}}(\pi_{\theta_{k+1}},\phi_{\tau_{k+1}}(\pi_{\theta_{k+1}}))$ and $J_{\tau_{k+1}}(\pi_{\tau_{k+1}}^{\star},\phi_{\tau_{k+1}}^{\star})-J_{\tau_k}(\pi_{\tau_k}^{\star},\phi_{\tau_k}^{\star})$ from Lemma \ref{lem:V_tau_diff} Eqs.~\eqref{lem:V_tau_diff:eq1} and \eqref{lem:V_tau_diff:eq3}. The second inequality uses Lemma \ref{lem:g_smooth}.

Similarly, we consider the decay of $\delta^{\phi}_{k}$.
\begin{align}
    \delta^{\phi}_{k+1}-\delta^{\phi}_{k}&=J_{\tau_{k+1}}(\pi_{\theta_{k+1}},\phi_{\psi_{k+1}})-g_{\tau_{k+1}}(\pi_{\theta_{k+1}})-J_{\tau_k}(\pi_{\theta_k},\phi_{\psi_k})+g_{\tau_{k}}(\pi_{\theta_k})\notag\\
    &=\left(J_{\tau_{k+1}}(\pi_{\theta_{k+1}},\phi_{\psi_{k+1}})\hspace{-2pt}-\hspace{-2pt}J_{\tau_{k}}(\pi_{\theta_{k+1}},\phi_{\psi_{k+1}})\right)\hspace{-2pt}+\hspace{-2pt}\left(J_{\tau_{k}}(\pi_{\theta_{k+1}},\phi_{\psi_{k+1}})\hspace{-2pt}-\hspace{-2pt}J_{\tau_k}(\pi_{\theta_{k+1}},\phi_{\psi_{k}})\right)\notag\\
    &\hspace{20pt}+\left(J_{\tau_k}(\pi_{\theta_{k+1}},\phi_{\psi_{k}})-J_{\tau_k}(\pi_{\theta_k},\phi_{\psi_k})\right)+\left(g_{\tau_{k}}(\pi_{\theta_k})-g_{\tau_{k+1}}(\pi_{\theta_{k+1}})\right).
    \label{thm:main:proof_eq4}
\end{align}

By Lemma \ref{lem:V_tau_diff} Eq.~\eqref{lem:V_tau_diff:eq1},
\begin{align}
    J_{\tau_{k+1}}(\pi_{\theta_{k+1}},\phi_{\psi_{k+1}})-J_{\tau_{k}}(\pi_{\theta_{k+1}},\phi_{\psi_{k+1}})\leq \frac{\tau_k-\tau_{k+1}}{1-\gamma}\log|\Bcal|.
    \label{thm:main:proof_eq5}
\end{align}

Using the $L_0$-smoothness of the value function derived in \eqref{thm:main:proof_eq2}
\begin{align*}
    &J_{\tau_k}(\pi_{\theta_{k+1}},\phi_{\psi_{k+1}})-J_{\tau_k}(\pi_{\theta_{k+1}},\phi_{\psi_{k}})\notag\\
    &\leq\langle \nabla_{\psi} J_{\tau_k}(\pi_{\theta_{k+1}},\phi_{\psi_{k}}),\psi_{k+1}-\psi_{k}\rangle+\frac{L_0}{2}\|\psi_{k+1}-\psi_k\|^2\notag\\
    &= -\beta_k\|\nabla_{\psi} J_{\tau_k}(\pi_{\theta_{k+1}},\phi_{\psi_{k}})\|^2+\frac{L_0\beta_k^2}{2}\|\nabla_{\psi} J_{\tau_k}(\pi_{\theta_{k+1}},\phi_{\psi_{k}})\|^2\notag\\
    &\leq -\frac{\beta_k}{2}\|\nabla_{\psi} J_{\tau_k}(\pi_{\theta_{k+1}},\phi_{\psi_{k}})\|^2\notag\\
    &\leq -\frac{(1-\gamma)\beta_k\tau_k\rho_{\min}^2}{|\Scal|}\left(\min_{s,b}\phi_{\psi_k}(b\mid s)\right)^2\left(J_{\tau_k}(\pi_{\theta_k},\phi_{\psi_k})-J_{\tau_k}(\pi_{\theta_k},\phi_{\tau_k}(\pi_{\theta_k}))\right)\notag\\
    &= -\frac{(1-\gamma)\beta_k\tau_k\rho_{\min}^2}{|\Scal|}\left(\min_{s,b}\phi_{\psi_k}(b\mid s)\right)^2\delta_k^{\phi},
\end{align*}
where the second inequality uses $\beta_k\leq\frac{1}{L_0}$ and the third inequality follows from Lemma \ref{lem:nonuniform_PL}.

Using Eq.~\eqref{lem:min_policy_iterates:eq2} of Lemma \ref{lem:min_policy_iterates} to further simplify this inequality,
\begin{align}
    J_{\tau_k}(\pi_{\theta_{k+1}},\phi_{\psi_{k+1}})-J_{\tau_k}(\pi_{\theta_{k+1}},\phi_{\psi_{k}})
    &\leq -\frac{3(1-\gamma)\beta_k\tau_k\rho_{\min}^2 c^2}{8|\Scal|}\delta_k^{\phi}
    \label{thm:main:proof_eq6}.
\end{align}

For the third term of \eqref{thm:main:proof_eq4}, we have from the $L_0$-smoothness of the value function derived in \eqref{thm:main:proof_eq2}
\begin{align}
    J_{\tau_k}(\pi_{\theta_{k+1}},\phi_{\psi_{k}})-J_{\tau_k}(\pi_{\theta_{k}},\phi_{\psi_{k}})&\leq\langle \nabla_{\theta} J_{\tau_k}(\pi_{\theta_{k}},\phi_{\psi_{k}}),\theta_{k+1}-\theta_k\rangle+\frac{L_0}{2}\|\theta_{k+1}-\theta_k\|^2\notag\\
    &= \alpha_k\|\nabla_{\theta} J_{\tau_k}(\pi_{\theta_{k}},\phi_{\psi_{k}})\|^2+\frac{L_0\alpha_k^2}{2}\|\nabla_{\theta} J_{\tau_k}(\pi_{\theta_{k}},\phi_{\psi_{k}})\|^2\notag\\
    &\leq \frac{3\alpha_k}{2}\|\nabla_{\theta} J_{\tau_k}(\pi_{\theta_{k}},\phi_{\psi_{k}})\|^2,
    \label{thm:main:proof_eq7}
\end{align}
where in the last inequality we use $\alpha_k L_0\leq 1$.

Using Lemma \ref{lem:g_smooth} and Lemma \ref{lem:V_tau_diff} \eqref{lem:V_tau_diff:eq4}, we bound the last term of \eqref{thm:main:proof_eq4}
\begin{align}
    &g_{\tau_k}(\pi_{\theta_k})-g_{\tau_{k+1}}(\pi_{\theta_{k+1}})\notag\\
    &=g_{\tau_k}(\pi_{\theta_k})-g_{\tau_{k}}(\pi_{\theta_{k+1}})+g_{\tau_k}(\pi_{\theta_{k+1}})-g_{\tau_{k+1}}(\pi_{\theta_{k+1}})\notag\\
    &\leq \frac{\alpha_k}{2}\left(\|\nabla_{\theta} J_{\tau_k}(\pi_{\theta_{k}},\phi_{\tau_k}(\pi_{\theta_k}))-\nabla_{\theta}J_{\tau_k}(\pi_{\theta_k},\phi_{\psi_k})\|^2-\|\nabla_{\theta} J_{\tau_k}(\pi_{\theta_{k}},\phi_{\tau_k}(\pi_{\theta_k}))\|^2\right)\notag\\
    &\hspace{20pt}+(\tau_k-\tau_{k+1})\log|\Acal|
    \label{thm:main:proof_eq8}
\end{align}

Using \eqref{thm:main:proof_eq5}-\eqref{thm:main:proof_eq8} in \eqref{thm:main:proof_eq4}, we have
\begin{align}
    \delta^{\phi}_{k+1} &=\delta^{\phi}_{k}+\left(J_{\tau_{k+1}}(\pi_{\theta_{k+1}},\phi_{\psi_{k+1}})\hspace{-2pt}-\hspace{-2pt}J_{\tau_k}(\pi_{\theta_{k+1}},\phi_{\psi_{k+1}})\right)+\left(J_{\tau_{k}}(\pi_{\theta_{k+1}},\phi_{\psi_{k+1}})\hspace{-2pt}-\hspace{-2pt}J_{\tau_k}(\pi_{\theta_{k+1}},\phi_{\psi_{k}})\right)\notag\\
    &\hspace{20pt}+\left(J_{\tau_k}(\pi_{\theta_{k+1}},\phi_{\psi_{k}})-J_{\tau_k}(\pi_{\theta_k},\phi_{\psi_k})\right)+\left(g_{\tau_{k}}(\pi_{\theta_k})-g_{\tau_{k+1}}(\pi_{\theta_{k+1}})\right)\notag\\
    &\leq\delta_k^{\phi}+\frac{\tau_k-\tau_{k+1}}{1-\gamma}\log|\Bcal|-\frac{3(1-\gamma)\beta_k\tau_k\rho_{\min}^2 c^2}{8|\Scal|}\delta_k^{\phi}+\frac{3\alpha_k}{2}\|\nabla_{\theta} J_{\tau_k}(\pi_{\theta_{k}},\phi_{\psi_{k}})\|^2\notag\\
    &\hspace{20pt}+\frac{\alpha_k}{2}\left(\|\nabla_{\theta} J_{\tau_k}(\pi_{\theta_{k}},\phi_{\tau_k}(\pi_{\theta_k}))-\nabla_{\theta}J_{\tau_k}(\pi_{\theta_k},\phi_{\psi_k})\|^2-\|\nabla_{\theta} J_{\tau_k}(\pi_{\theta_{k}},\phi_{\tau_k}(\pi_{\theta_k}))\|^2\right)\notag\\
    &\hspace{20pt}+(\tau_k-\tau_{k+1})\log|\Acal|\notag\\
    &\leq (1-\frac{3(1-\gamma)\beta_k\tau_k \rho_{\min}^2 c^2}{8|\Scal|})\delta^{\phi}_{k}+\frac{3\alpha_k}{2}\|\nabla_{\theta} J_{\tau_k}(\pi_{\theta_{k}},\phi_{\psi_{k}})\|^2\notag\\
    &\hspace{20pt}+\frac{\alpha_k}{2}\left(\|\nabla_{\theta} J_{\tau_k}(\pi_{\theta_{k}},\phi_{\tau_k}(\pi_{\theta_k}))-\nabla_{\theta}J_{\tau_k}(\pi_{\theta_k},\phi_{\psi_k})\|^2-\|\nabla_{\theta} J_{\tau_k}(\pi_{\theta_{k}},\phi_{\tau_k}(\pi_{\theta_k}))\|^2\right)\notag\\
    &\hspace{20pt}+\frac{\tau_k-\tau_{k+1}}{1-\gamma}(\log|\Acal|+\log|\Bcal|).
    \label{thm:main:proof_eq9}
\end{align}

Combining \eqref{thm:main:proof_eq3} and \eqref{thm:main:proof_eq9},
\begin{align*}
    &3\delta^{\pi}_{k+1}+\delta^{\phi}_{k+1}\\
    &\leq 3\delta^{\pi}_{k}+\frac{3\alpha_k}{2}\left(\|\nabla_{\theta} J_{\tau_k}(\pi_{\theta_{k}},\phi_{\tau_k}(\pi_{\theta_k}))-\nabla_{\theta}J_{\tau_k}(\pi_{\theta_k},\phi_{\psi_k})\|^2-\|\nabla_{\theta} J_{\tau_k}(\pi_{\theta_{k}},\phi_{\tau_k}(\pi_{\theta_k}))\|^2\right)\notag\\
    &\hspace{20pt}+\frac{3(\tau_k-\tau_{k+1})}{1-\gamma}(\log|\Acal|+\log|\Bcal|)+(1-\frac{3(1-\gamma)\beta_k\tau_k \rho_{\min}^2 c^2}{8|\Scal|})\delta^{\phi}_{k}\notag\\
    &\hspace{20pt}+\frac{\alpha_k}{2}\left(\|\nabla_{\theta} J_{\tau_k}(\pi_{\theta_{k}},\phi_{\tau_k}(\pi_{\theta_k}))-\nabla_{\theta}J_{\tau_k}(\pi_{\theta_k},\phi_{\psi_k})\|^2-\|\nabla_{\theta} J_{\tau_k}(\pi_{\theta_{k}},\phi_{\tau_k}(\pi_{\theta_k}))\|^2\right)\notag\\
    &\hspace{20pt}+\frac{3\alpha_k}{2}\|\nabla_{\theta} J_{\tau_k}(\pi_{\theta_{k}},\phi_{\psi_{k}})\|^2+\frac{\tau_k-\tau_{k+1}}{1-\gamma}(\log|\Acal|+\log|\Bcal|)\notag\\
    &\leq 3\delta^{\pi}_{k}+(1-\frac{3(1-\gamma)\beta_k\tau_k \rho_{\min}^2 c^2}{8|\Scal|})\delta^{\phi}_{k}+\frac{3\alpha_k}{2}\|\nabla_{\theta} J_{\tau_k}(\pi_{\theta_{k}},\phi_{\psi_{k}})\|^2\notag\\
    &\hspace{20pt}+2\alpha_k\|\nabla_{\theta} J_{\tau_k}(\pi_{\theta_{k}},\phi_{\tau_k}(\pi_{\theta_k}))-\nabla_{\theta}J_{\tau_k}(\pi_{\theta_k},\phi_{\psi_k})\|^2-2\alpha_k\|\nabla_{\theta} J_{\tau_k}(\pi_{\theta_{k}},\phi_{\tau_k}(\pi_{\theta_k}))\|^2\notag\\
    &\hspace{20pt}+\frac{4(\tau_k-\tau_{k+1})}{1-\gamma}(\log|\Acal|+\log|\Bcal|).
\end{align*}
Simplifying this inequality with 
\begin{align*}
    \|\nabla_{\theta} \hspace{-1pt}J_{\tau_k}\hspace{-1pt}(\pi_{\theta_{k}},\phi_{\psi_{k}})\|^2\hspace{-2pt}&=\hspace{-2pt}\|\nabla_{\theta} J_{\tau_k}(\pi_{\theta_{k}},\phi_{\tau_k}(\pi_{\theta_k}))-\left(\nabla_{\theta} J_{\tau_k}(\pi_{\theta_{k}},\phi_{\tau_k}(\pi_{\theta_k}))-\nabla_{\theta}J_{\tau_k}(\pi_{\theta_k},\phi_{\psi_k})\right)\|^2\\
    &\leq\hspace{-2pt} \|\nabla_{\theta} J_{\tau_k}(\pi_{\theta_{k}},\phi_{\tau_k}(\pi_{\theta_k}))\|^2+\|\nabla_{\theta} J_{\tau_k}(\pi_{\theta_{k}},\phi_{\tau_k}(\pi_{\theta_k}))-\nabla_{\theta}J_{\tau_k}(\pi_{\theta_k},\phi_{\psi_k})\|^2\\
    &\hspace{20pt}+2\langle\nabla_{\theta} J_{\tau_k}(\pi_{\theta_{k}},\phi_{\tau_k}(\pi_{\theta_k})),\nabla_{\theta} J_{\tau_k}(\pi_{\theta_{k}},\phi_{\tau_k}(\pi_{\theta_k}))-\nabla_{\theta}J_{\tau_k}(\pi_{\theta_k},\phi_{\psi_k})\rangle\\
    &\leq\hspace{-2pt} \frac{5}{4}\|\nabla_{\theta}\hspace{-1pt} J_{\tau_k}(\pi_{\theta_{k}},\phi_{\tau_k}(\pi_{\theta_k}))\|^2\hspace{-2pt}+\hspace{-2pt}5\|\nabla_{\theta} J_{\tau_k}(\pi_{\theta_{k}},\phi_{\tau_k}(\pi_{\theta_k}))\hspace{-2pt}-\hspace{-2pt}\nabla_{\theta}\hspace{-1pt}J_{\tau_k}(\pi_{\theta_k},\phi_{\psi_k})\|^2
\end{align*}
we have
\begin{align}
    3\delta^{\pi}_{k+1}+\delta^{\phi}_{k+1}&\leq 3\delta^{\pi}_{k}+(1-\frac{3(1-\gamma)\beta_k\tau_k \rho_{\min}^2 c^2}{8|\Scal|})\delta^{\phi}_{k}-\frac{\alpha_k}{8}\|\nabla_{\theta} J_{\tau_k}(\pi_{\theta_{k}},\phi_{\tau_k}(\pi_{\theta_k}))\|^2\notag\\
    &\hspace{20pt}+\frac{19\alpha_k}{2}\|\nabla_{\theta} J_{\tau_k}(\pi_{\theta_{k}},\phi_{\tau_k}(\pi_{\theta_k}))-\nabla_{\theta}J_{\tau_k}(\pi_{\theta_k},\phi_{\psi_k})\|^2\notag\\
    &\hspace{20pt}+\frac{4(\tau_k-\tau_{k+1})}{1-\gamma}(\log|\Acal|+\log|\Bcal|)
    \label{thm:main:proof_eq10}
\end{align}

Using Lemma \ref{lem:nonuniform_PL} to bound $-\|\nabla_{\theta} J_{\tau_k}(\pi_{\theta_{k}},\phi_{\tau_k}(\pi_{\theta_k}))\|^2$, 
\begin{align}
    &-\|\nabla_{\theta} J_{\tau_k}(\pi_{\theta_{k}},\phi_{\tau_k}(\pi_{\theta_k}))\|^2\notag  \\
    &\leq-\frac{2(1-\gamma)\tau_k\rho_{\min}^2}{|\Scal|}\left(\min_{s,a}\pi_{\theta_k}(a\mid s)\right)^2\left(J_{\tau_k}(\pi_{\tau_k}(\phi_{\tau_k}(\pi_{\theta_k})),\phi_{\tau_k}(\pi_{\theta_k}))-J_{\tau_k}(\pi_{\theta_k},\phi_{\tau_k}(\pi_{\theta_k}))\right)\notag\\
    &\leq-\frac{2(1-\gamma)\tau_k\rho_{\min}^2}{|\Scal|}\left(\min_{s,a}\pi_{\theta_k}(a\mid s)\right)^2\left(J_{\tau_k}(\pi_{\tau_k}^{\star},\phi_{\tau_k}^{\star})-J_{\tau_k}(\pi_{\theta_k},\phi_{\tau_k}(\pi_{\theta_k}))\right),
    \label{thm:main:proof_eq11}
\end{align}
where the second inequality follows from
\begin{align*}
    J_{\tau_k}(\pi_{\tau_k}(\phi_{\tau_k}(\pi_{\theta_k})),\phi_{\tau_k}(\pi_{\theta_k}))=\max_{\pi}J_{\tau_k}(\pi,\phi_{\tau_k}(\pi_{\theta_k}))\geq \max_{\pi}\min_{\phi}J_{\tau_k}(\pi,\phi)=J_{\tau_k}(\pi_{\tau_k}^{\star},\phi_{\tau_k}^{\star}).
\end{align*}

From Lemma \ref{lem:min_policy_iterates} Eq.~\eqref{lem:min_policy_iterates:eq1}, $-\left(\min_{s,a}\pi_{\theta_k}(a\mid s)\right)^2\leq-\frac{3c^2}{8}$, which further simplifies \eqref{thm:main:proof_eq11}
\begin{align}
    &-\|\nabla_{\theta} J_{\tau_k}(\pi_{\theta_{k}},\phi_{\tau_k}(\pi_{\theta_k}))\|^2\notag\\
    &\leq-\frac{2(1-\gamma)\tau_k\rho_{\min}^2}{|\Scal|}\left(\min_{s,a}\pi_{\theta_k}(a\mid s)\right)^2\left(J_{\tau_k}(\pi_{\tau_k}^{\star},\phi_{\tau_k}^{\star})-J_{\tau_k}(\pi_{\theta_k},\phi_{\tau_k}(\pi_{\theta_k}))\right)\notag\\
    &=-\frac{2(1-\gamma)\tau_k\rho_{\min}^2}{|\Scal|}\left(\min_{s,a}\pi_{\theta_k}(a\mid s)\right)^2\delta_k^{\pi}\leq-\frac{3(1-\gamma)\tau_k\rho_{\min}^2 c^2}{4|\Scal|}\delta_k^{\pi}.
    \label{thm:main:proof_eq12.1}
\end{align}

For $\|\nabla_{\theta} J_{\tau_k}(\pi_{\theta_{k}},\phi_{\tau_k}(\pi_{\theta_k}))-\nabla_{\theta}J_{\tau_k}(\pi_{\theta_k},\phi_{\psi_k})\|^2$, we have from the $L_0$-smoothness of the value function derived in \eqref{thm:main:proof_eq2}
\begin{align}
    \|\nabla_{\theta} J_{\tau_k}(\pi_{\theta_{k}},\phi_{\tau_k}(\pi_{\theta_k}))\hspace{-2pt}-\hspace{-2pt}\nabla_{\theta}J_{\tau_k}(\pi_{\theta_k},\phi_{\psi_k})\|^2&\leq L_0^2\|\phi_{\tau_k}(\pi_{\theta_k})-\phi_{\psi_k}\|^2\notag\\
    &\leq \frac{2\log(2) L_0^2}{\tau_k\rho_{\min}}\left(J_{\tau_k}(\pi_{\theta_k},\phi_{\psi_k})\hspace{-2pt}-\hspace{-2pt}J_{\tau_k}(\pi_{\theta_{k}},\phi_{\tau_k}(\pi_{\theta_k}))\right)\notag\\
    &=\frac{2\log(2) L_0^2}{\tau_k\rho_{\min}}\delta_k^{\phi},\label{thm:main:proof_eq12.2}
\end{align}
where the second inequality follows from Lemma \ref{lem:quadratic_growth} Eq.~\eqref{lem:quadratic_growth:eq2}.

Using \eqref{thm:main:proof_eq12.1} and \eqref{thm:main:proof_eq12.2} in \eqref{thm:main:proof_eq10},
\begin{align}
    &3\delta^{\pi}_{k+1}+\delta^{\phi}_{k+1}\notag\\
    &\leq 3\delta^{\pi}_{k}+(1-\frac{3(1-\gamma)\beta_k\tau_k\rho_{\min}^2 c^2}{8|\Scal|})\delta^{\phi}_{k}-\frac{\alpha_k}{8}\|\nabla_{\theta} J_{\tau_k}(\pi_{\theta_{k}},\phi_{\tau_k}(\pi_{\theta_k}))\|^2\notag\\
    &\hspace{20pt}+\frac{19\alpha_k}{2}\|\nabla_{\theta} J_{\tau_k}(\pi_{\theta_{k}},\phi_{\tau_k}(\pi_{\theta_k}))-\nabla_{\theta}J_{\tau_k}(\pi_{\theta_k},\phi_{\psi_k})\|^2+\frac{4(\tau_k-\tau_{k+1})}{1-\gamma}(\log|\Acal|+\log|\Bcal|)\notag\\
    &\leq 3\delta^{\pi}_{k}+(1-\frac{3(1-\gamma)\beta_k\tau_k\rho_{\min}^2 c^2}{8|\Scal|})\delta^{\phi}_{k}-\frac{3\alpha_k(1-\gamma)\tau_k\rho_{\min}^2 c^2}{32|\Scal|}\delta_k^{\pi}\notag\\
    &\hspace{20pt}+\frac{19\log(2) L_0^2 \alpha_k}{\tau_k\rho_{\min}}\delta_k^{\phi}+\frac{4(\tau_k-\tau_{k+1}))}{1-\gamma}(\log|\Acal|+\log|\Bcal|)\notag\\
    &=3(1-\frac{(1-\gamma) \alpha_k \tau_k \rho_{\min}^2 c^2}{32|\Scal|})\delta_k^{\pi}+(1-\frac{3(1-\gamma)\beta_k\tau_k\rho_{\min}^2 c^2}{8|\Scal|}+\frac{19\log(2) L_0^2 \alpha_k}{\tau_k\rho_{\min}})\delta^{\phi}_{k}\notag\\
    &\hspace{20pt}+\frac{4(\tau_k-\tau_{k+1})}{1-\gamma}(\log|\Acal|+\log|\Bcal|).\label{thm:main:proof_eq13}
\end{align}



With the step size rule $\frac{\alpha_0}{\beta_0}\leq\min\{\frac{(1-\gamma)\tau_0^2\rho_{\min}^3 c^2}{152\log(2)L_0^2|\Scal|},1\}$, we can simplify \eqref{thm:main:proof_eq13},
\begin{align*}
    3\delta^{\pi}_{k+1}+\delta^{\phi}_{k+1}&\leq 3(1-\frac{(1-\gamma) \alpha_k \tau_k \rho_{\min}^2 c^2}{32|\Scal|})\delta_k^{\pi}+(1\hspace{-2pt}-\hspace{-2pt}\frac{3(1-\gamma)\beta_k\tau_k\rho_{\min}^2 c^2}{8|\Scal|}+\frac{19\log(2) L_0^2 \alpha_k}{\tau_k\rho_{\min}})\delta^{\phi}_{k}\notag\\
    &\hspace{20pt}+\frac{4(\tau_k-\tau_{k+1})}{1-\gamma}(\log|\Acal|+\log|\Bcal|)\\
    &\leq 3(1-\frac{(1-\gamma) \alpha_k \tau_k \rho_{\min}^2 c^2}{32|\Scal|})\delta_k^{\pi}\notag\\
    &\hspace{20pt}+(1-\frac{3(1-\gamma)\beta_k\tau_k\rho_{\min}^2 c^2}{8|\Scal|}+\frac{19\log(2) L_0^2 }{\tau_k\rho_{\min}}\frac{(1-\gamma)\rho_{\min}^3 c^2\tau_k^2\beta_k}{152\log(2)L_0^2|\Scal|})\delta^{\phi}_{k}\notag\\
    &\hspace{20pt}+\frac{4(\tau_k-\tau_{k+1})}{1-\gamma}(\log|\Acal|+\log|\Bcal|)\\
    &\leq 3(1-\frac{(1-\gamma) \alpha_k \tau_k \rho_{\min}^2 c^2}{32|\Scal|})\delta_k^{\pi}+(1-\frac{(1-\gamma)\beta_k\tau_k\rho_{\min}^2 c^2}{4|\Scal|})\delta^{\phi}_{k}\notag\\
    &\hspace{20pt}+\frac{4(\tau_k-\tau_{k+1})}{1-\gamma}(\log|\Acal|+\log|\Bcal|)\\
    &\leq (1-\frac{(1-\gamma) \alpha_k \tau_k \rho_{\min}^2 c^2}{32|\Scal|})(3\delta_k^{\pi}+\delta^{\phi}_{k})+\frac{4(\tau_k-\tau_{k+1})}{1-\gamma}(\log|\Acal|+\log|\Bcal|)\\
    &\leq (1-\frac{(1-\gamma) \rho_{\min}^2 c^2  \alpha_0 \tau_0}{32|\Scal|(k+h)})\frac{C_1}{(k+h)^{1/3}}+\frac{32\tau_0}{3(1-\gamma)(k+h)^{4/3}}(\log|\Acal|+\log|\Bcal|),
\end{align*}
where the last inequality follows from \eqref{thm:main:proof_eq1} and Lemma \ref{lem:stepsize_diff}.

Letting $D_1=\frac{(1-\gamma)\rho_{\min}^2 c^2}{32|\Scal|}$ and $D_2=\frac{32}{3(1-\gamma)}(\log|\Acal|+\log|\Bcal|)$,
\begin{align*}
    3\delta^{\pi}_{k+1}+\delta^{\phi}_{k+1}&\leq \left(1-\frac{D_1 \alpha_0\tau_0}{k+h}\right)\frac{C_1\tau_0}{(k+h)^{1/3}}+\frac{D_2\tau_0}{(k+1)^{4/3}}\\
    &=\left(k+h-D_1 \alpha_0\tau_0+\frac{D_2}{C_1}\right) \frac{C_1\tau_0}{(k+h)^{4/3}}.
\end{align*}
By requiring
\begin{align*}
    \tau_0=\frac{65536\log(2)(\log|\Acal|+\log|\Bcal|)+96(1-\gamma)\rho_{\min}c^2}{3(1-\gamma)^2 \rho_{\min}^3 c^4\alpha_0}=\frac{1}{D_1\alpha_0}(1+\frac{D_2}{C_1}),
\end{align*}
we have
\begin{align*}
    3\delta^{\pi}_{k+1}+\delta^{\phi}_{k+1}&\leq \left(k+h-D_1 \alpha_0\tau_0+\frac{D_2}{C_1}\right)\cdot \frac{C_1\tau_0}{(k+h)^{4/3}}\\
    &= \left(k+h-(1+\frac{D_2}{C_1})+\frac{D_2}{C_1}\right)\cdot \frac{C_1\tau_0}{(k+h)^{4/3}}\\
    &= \frac{C_1\tau_0 (k-1+h)}{(k+h)^{4/3}},
\end{align*}

Since 
$(k-1+h)^3(k+1+h)\leq(k+h)^4$
for all $k\geq0$ and $h\geq1$, we have
\[\frac{k-1+h}{(k+h)^{4/3}}=\frac{(k-1+h)(k+1+h)^{1/3}}{(k+1)^{4/3}(k+1+h)^{1/3}}\leq\frac{(k+h)^{4/3}}{(k+h)^{4/3}(k+1+h)^{1/3}}=\frac{1}{(k+1+h)^{1/3}},\]
which leads to
\begin{align*}
    3\delta^{\pi}_{k+1}+\delta^{\phi}_{k+1}\leq \frac{C_1\tau_0 (k-1+h)}{(k+h)^{4/3}}\leq \frac{C_1\tau_0}{(k+1+h)^{1/3}}=\frac{\rho_{\min}\tau_0 c^2}{64\log(2)(k+1+h)^{1/3}}.
\end{align*}

This finishes our induction and implies that for all $k\geq0$
\begin{align*}
    J_{\tau_k}(\pi_{\tau_k}^{\star},\phi_{\tau_k}^{\star}) - J_{\tau_k}(\pi_{\theta_k},\phi_{\tau_k}(\pi_{\theta_k}))&\leq \frac{C_1\tau_0}{3(k+h)^{1/3}},\\
    J_{\tau_k}(\pi_{\theta_k},\phi_{\psi_k})- J_{\tau_k}(\pi_{\theta_k},\phi_{\tau_k}(\pi_{\theta_k}))&\leq \frac{C_1\tau_0}{(k+h)^{1/3}}.
\end{align*}

\textbf{Bounding the difference between value functions with and without the regularization:}

Ultimately, we are interested in $J(\pi^{\star},\phi^{\star}) - J(\pi_{\theta_k},\phi_{0}(\pi_{\theta_k}))$ and $J(\pi_{\theta_k},\phi_{\psi_k})- J(\pi_{\theta_k},\phi_{0}(\pi_{\theta_k}))$, which measure the performance of $\pi_{\theta_k}$ and $\phi_{\psi_k}$ in the original un-regularized Markov game.

By Lemma \ref{lem:V_tau_diff} Eq.~\eqref{lem:V_tau_diff:eq3}, \eqref{lem:V_tau_diff:eq4}, and \eqref{lem:V_tau_diff:eq1},
\begin{align*}
     J_{\tau_k}(\pi_{\tau_k}^{\star},\phi_{\tau_k}^{\star}) - J(\pi^{\star},\phi^{\star}) &\geq -\tau_k\log|\Bcal|\\
     J_{\tau_k}(\pi_{\theta_k},\phi_{\tau_k}(\pi_{\theta_k}))-J(\pi_{\theta_k},\phi_{0}(\pi_{\theta_k}))&\leq\tau_k\log|\Acal|\\
     J_{\tau_k}(\pi_{\theta_k},\phi_{\psi_k})-J(\pi_{\theta_k},\phi_{\psi_k})&\geq-\frac{\tau_k}{1-\gamma}\log|\Bcal|.
\end{align*}

Therefore,
\begin{align*}
    J(\pi^{\star},\phi^{\star}) - J(\pi_{\theta_k},\phi_{0}(\pi_{\theta_k}))&=J(\pi^{\star},\phi^{\star})-J_{\tau_k}(\pi_{\tau_k}^{\star},\phi_{\tau_k}^{\star})+J_{\tau_k}(\pi_{\tau_k}^{\star},\phi_{\tau_k}^{\star})-J_{\tau_k}(\pi_{\theta_k},\phi_{\tau_k}(\pi_{\theta_k}))\\
    &\hspace{20pt}+J_{\tau_k}(\pi_{\theta_k},\phi_{\tau_k}(\pi_{\theta_k})) - J(\pi_{\theta_k},\phi_{0}(\pi_{\theta_k}))\\
    &\leq \tau_k\log|\Bcal|+\frac{C_1\tau_0}{3(k+h)^{1/3}}+\tau_k\log|\Acal|\\
    &=\frac{C_1\tau_0+3(\log|\Acal|+\log|\Bcal|)\tau_0}{3{(k+h)^{1/3}}},
\end{align*}
and
\begin{align*}
    J(\pi_{\theta_k},\phi_{\psi_k})- J(\pi_{\theta_k},\phi_{0}(\pi_{\theta_k}))&=J(\pi_{\theta_k},\phi_{\psi_k})- J_{\tau_k}(\pi_{\theta_k},\phi_{\psi_k})\\
    &\hspace{20pt}+ J_{\tau_k}(\pi_{\theta_k},\phi_{\psi_k})-J_{\tau_k}(\pi_{\theta_k},\phi_{\tau_k}(\pi_{\theta_k}))\\
    &\hspace{20pt}+J_{\tau_k}(\pi_{\theta_k},\phi_{\tau_k}(\pi_{\theta_k}))- J(\pi_{\theta_k},\phi_{0}(\pi_{\theta_k}))\\
    &\leq \frac{\tau_k}{1-\gamma}\log|\Bcal|+\frac{C_1\tau_0}{(k+h)^{1/3}}+\tau_k\log|\Acal|\\
    &\leq \frac{(1-\gamma)C_1\tau_0+(\log|\Acal|+\log|\Bcal|)\tau_0}{(1-\gamma){(k+h)^{1/3}}}.
\end{align*}

\begin{remark}\label{rem:stepsize_exist}
To select $\alpha_0$, $\beta_0$, $\tau_0$, and $h$, we first make $\tau_0=\lambda h^{1/3}$ for some $\lambda>0$ large enough. This choice guarantees the validity of \eqref{eq:stepsize_exist_eq1} (we just need $\delta^{\pi}_{0}+\delta^{\phi}_{0}\leq C_1\lambda$). Viewing \eqref{eq:stepsize_exist_eq2}, it means
\[\alpha_0=\frac{65536\log(2)(\log|\Acal|+\log|\Bcal|)+96(1-\gamma)\rho_{\min}c^2}{3(1-\gamma)^2 \rho_{\min}^3 c^4 \lambda h^{\frac{1}{3}}}.\]
Now that $\lambda$ is fixed, to ensure \eqref{eq:stepsize_exist_eq3}, we choose $h$ large enough to observe
\begin{align*}
&\frac{65536\log(2)(\log|\Acal|+\log|\Bcal|)+96(1-\gamma)\rho_{\min}c^2}{3(1-\gamma)^2 \rho_{\min}^3 c^4 \lambda h}=\frac{\alpha_0}{h^{\frac{2}{3}}}\\
&\hspace{180pt}\leq(2L_{\Hcal}+4L_{\Hcal}^2 C_2)\lambda+(L_{\Hcal}+4L_{\Hcal}^2C_2)+\frac{L_{\Hcal}^2 C_2}{\lambda}.
\end{align*}
Once $\lambda$ and $h$ are chosen, $\alpha_0$, $\tau_0$, and $h$ are determined.
Finally, since $\frac{ (1-\gamma)\tau_0^2\rho_{\min}^3 c^2}{152\log(2)|\Scal|L_0^2}\leq 1$, we just need to select $\beta_0\in[\frac{152\log(2)|\Scal|L_0^2\alpha_0}{(1-\gamma)\tau_0^2\rho_{\min}^3 c^2},\frac{1}{L_0}]$.
Recall that $L_0=L_{\Hcal}(2\tau_0+1)$, it can be easily seen that the lower bound $\frac{152\log(2)|\Scal|L_0^2\alpha_0}{(1-\gamma)\tau_0^2\rho_{\min}^3 c^2}=\Ocal(\frac{1}{\lambda^3 h^{1/3}})$, which is much smaller than the upper bound $\frac{1}{L_0}=\Ocal(\frac{1}{\tau_0})=\Ocal(\frac{1}{\lambda h^{1/3}})$ since $\lambda$ was large enough.
\end{remark}

\qed

%% file: Proof_Lemma.tex
\section{Proof of Lemmas}

\subsection{Proof of Lemma \ref{lem:quadratic_growth}}

For a given $\phi$, let $\hat{\pi}\in\pi_{\tau}(\phi)$ (which is a possibly non-unique maximizer).

According to \citet{mei2020global}[Lemma 26],
\begin{align*}
    J_{\tau}(\hat{\pi},\phi)-J_{\tau}(\pi,\phi)=\frac{\tau}{1-\gamma}\sum_{s\in\Scal} d_{\rho}^{\pi,\phi}(s) D_{KL}(\pi(\cdot\mid s)||\hat{\pi}(\cdot\mid s)).
\end{align*}

The Pinsker's inequality states that for any two probability distributions $p_1$ and $p_2$
\begin{align*}
    D_{KL}(p_1||p_2)\geq\frac{1}{2\log(2)}\|p_1-p_2\|_1^2.
\end{align*}

Using this inequality, 
\begin{align*}
    J_{\tau}(\hat{\pi},\phi)-J_{\tau}(\pi,\phi)&=\frac{\tau}{1-\gamma}\sum_{s\in\Scal} d_{\rho}^{\pi,\phi}(s) D_{KL}(\pi(\cdot\mid s)||\hat{\pi}(\cdot\mid s))\notag\\
    &\geq \frac{\tau}{2\log(2)(1-\gamma)}\sum_{s\in\Scal} d_{\rho}^{\pi,\phi}(s) \|\pi(\cdot\mid s)-\hat{\pi}(\cdot\mid s)\|_1^2\notag\\
    &\geq \frac{\tau}{2\log(2)(1-\gamma)}\sum_{s\in\Scal} (1-\gamma)\rho(s) \|\pi(\cdot\mid s)-\hat{\pi}(\cdot\mid s)\|_1^2\notag\\
    &\geq \frac{\tau\min_{s\in\Scal}\rho(s)}{2\log(2)}\sum_{s\in\Scal}  \|\pi(\cdot\mid s)-\hat{\pi}(\cdot\mid s)\|_1^2\notag\\
    &\geq \frac{\tau\min_{s\in\Scal}\rho(s)}{2\log(2)}\|\pi-\hat{\pi}\|^2,
\end{align*}
where the second inequality follows from the fact that $d_{\rho}^{\pi,\hat{\phi}}(s)\geq(1-\gamma)\rho(s)$ entry-wise. This inequality means that $\hat{\pi}\in\pi_{\tau}(\phi)$ has to be unique, as no other policy can achieve the same value function.

The same argument can be used to show Eq.~\eqref{lem:quadratic_growth:eq2}.

\qed

\subsection{Proof of Lemma \ref{lem:minimax_regMG}}

Let $(\pi_1,\phi_1)$, $(\pi_2,\phi_2)$ be optimal solution pairs to the maximin and minimax problem, respectively, 
\begin{align}
    (\pi_1,\phi_1)\in\argmax_{\pi\in\Delta_{\Acal}^{\Scal}}\argmin_{\phi\in\Delta_{\Bcal}^{\Scal}}J_{\tau}(\pi,\phi) \quad \text { and } \quad(\pi_2,\phi_2)\in\argmin_{\phi\in\Delta_{\Bcal}^{\Scal}}\argmax_{\pi\in\Delta_{\Acal}^{\Scal}}J_{\tau}(\pi,\phi).
    \label{lem:minimax_regMG:eq1}
\end{align}
Since the policy simplex is a compact set, $(\pi_1,\phi_1)$ and $(\pi_2,\phi_2)$ exist and are well-defined. The following minimax inequality always holds
\begin{align}
    J_{\tau}(\pi_1,\phi_1) =\max_{\pi\in\Delta_{\Acal}^{\Scal}}\min_{\phi\in\Delta_{\Bcal}^{\Scal}}J_{\tau}(\pi,\phi) \leq \min_{\phi\in\Delta_{\Bcal}^{\Scal}}\max_{\pi\in\Delta_{\Acal}^{\Scal}}J_{\tau}(\pi,\phi)=J_{\tau}(\pi_2,\phi_2).
    \label{lem:minimax_regMG:ineq_minimax}
\end{align}
We first want to show that $\pi_1=\pi_{\tau}(\phi_1)$ and $\phi_1=\phi_{\tau}(\pi_1)$. Since
\begin{align*}
    J_{\tau}(\pi_1,\phi_1)=\max_{\pi\in\Delta_{\Acal}^{\Scal}}\min_{\phi\in\Delta_{\Bcal}^{\Scal}}J_{\tau}(\pi,\phi)=\min_{\phi\in\Delta_{\Bcal}^{\Scal}}J_{\tau}(\pi_1,\phi)=J_{\tau}(\pi_1,\phi_{\tau}(\pi_1)),
\end{align*}
we have $\phi_1\in\phi_{\tau}(\pi_1)$, and Lemma~\ref{lem:quadratic_growth} further implies $\phi_1=\phi_{\tau}(\pi_1)$ is unique. 
In addition, we know that $\pi_1$ is the optimizer of $g_{\tau}$ defined in \eqref{eq:def_g}. Let $\theta_1$ be an softmax parameter for $\pi_1$ (e.g. $\theta_1(s,a)=\log\pi(a\mid s)$ for all $s,a$). Since $\pi_1$ is an optimizer of $g_{\tau}$ in policy space, $\theta_1$ must also be an (not necessarily unique) optimizer of $\tilde{g}_{\tau}(\theta)=\min_{\phi}J_{\tau}(\pi_{\theta},\phi)$ in the parameter space. Therefore, we have $\forall \theta\in\mathbb{R}^{\Scal\times\Acal}$
\begin{align}
    0\geq\langle\nabla_{\theta}g_{\tau}(\pi_{\theta_1}),\theta-\theta_1\rangle
    =\langle\nabla_{\theta}J_{\tau}(\pi_{\theta_1},\phi_1),\theta-\theta_1\rangle,\label{lem:minimax_regMG:eq3}
\end{align}
where the first equality follows from Danskin's Theorem in \eqref{eq:Danskin}. 
Since $\theta$ is not constrained, \eqref{lem:minimax_regMG:eq3} means that
\begin{align*}
    \nabla_{\theta} J_{\tau}(\pi_{\theta_1},\phi_1)=0,
\end{align*}
implying that $\theta_1$ is a stationary point of
\begin{align*}
    \max_{\theta}J_{\tau}(\pi_{\theta},\phi_1).
\end{align*}
By Lemma~\ref{lem:nonuniform_PL}, every stationary point is also globally optimal. Therefore, we have $\pi_1=\pi_{\theta_1}=\pi_{\tau}(\phi_1)$.




A consequence of $\pi_1=\pi_{\tau}(\phi_1)$ and $\phi_1=\phi_{\tau}(\pi_1)$ is that $(\pi_1,\phi_1)$ is the unique optimal solution pair to the maximin problem, i.e. there does not exist $(\widehat{\pi}_1,\widehat{\phi}_1)\neq(\pi_1,\phi_1)$ such that $(\widehat{\pi}_1,\widehat{\phi}_1)\in\argmax_{\pi\in\Delta_{\Acal}^{\Scal}}\argmin_{\phi\in\Delta_{\Bcal}^{\Scal}}J_{\tau}(\pi,\phi)$. To see this, let us suppose that such a pair $(\widehat{\pi}_1,\widehat{\phi}_1)$ does exist. Then, the only possibility is $\widehat{\pi}_1\neq\pi_1$ and $\widehat{\phi}_1\neq\phi_1$ by Lemma \ref{lem:quadratic_growth}. Since $\widehat{\pi}_1\neq\pi_{\tau}(\phi_1)$ and $\phi_1\neq\phi_{\tau}(\widehat{\pi}_1)$, we have
\begin{align*}
    J_{\tau}(\widehat{\pi}_1,\phi_1) < J_{\tau}(\pi_1,\phi_1)=J_{\tau}(\widehat{\pi}_1,\widehat{\phi}_1) < J_{\tau}(\widehat{\pi}_1,\phi_1),
\end{align*}
which creates a contradiction.

Similarly, it can be shown that 
\begin{align*}
    \pi_2=\phi_{\tau}(\phi_2),\quad\text{and}\quad\phi_2=\phi_{\tau}(\pi_2),
\end{align*}
and that $(\pi_2,\phi_2)$ is the unique optimal solution pair to the minimax problem.

We now aim prove that $(\pi_1,\phi_1)=(\pi_2,\phi_2)$, i.e. the minimax and maximin problem have the same solution. Suppose $(\pi_1,\phi_1)\neq(\pi_2,\phi_2)$, which means that $\pi_1\neq\pi_2$ and $\phi_1\neq\phi_2$ have to hold due to Lemma \ref{lem:quadratic_growth}. Since $\pi_2\neq\pi_{\tau}(\phi_1)$ and $\phi_1\neq\phi_{\tau}(\pi_2)$, we have from \eqref{lem:minimax_regMG:ineq_minimax}
\begin{align*}
    J_{\tau}(\pi_2,\phi_1) < J_{\tau}(\pi_1,\phi_1)\leq J_{\tau}(\pi_2,\phi_2) < J_{\tau}(\pi_2,\phi_1).
\end{align*}
This is again a contradiction. Therefore, $(\pi_1,\phi_1)=(\pi_2,\phi_2)$ has to be true. Then, \eqref{lem:minimax_regMG:ineq_minimax} leads to
\begin{align*}
    \max_{\pi\in\Delta_{\Acal}^{\Scal}}\min_{\phi\in\Delta_{\Bcal}^{\Scal}}J_{\tau}(\pi,\phi)= \max_{\pi\in\Delta_{\Acal}^{\Scal}}\min_{\phi\in\Delta_{\Bcal}^{\Scal}}J_{\tau}(\pi,\phi).
\end{align*}

We also know that the Nash equilibrium has to be unique in this case, as the maximin and minimax problems both have a unique solution pair that agrees with each other.

\qed

\subsection{Proof of Lemma \ref{lem:V_tau_diff}}

By the definition of the value function,
\begin{align*}
    &J_{\tau}(\pi,\phi)-J_{\tau'}(\pi,\phi)\notag\\
    &=\mathbb{E}\left[\sum_{k=0}^{\infty} \gamma^k \Big(r\left(s_k, a_k, b_k\right)-\tau\log\pi(a_k\mid s_k)+\tau\log\phi(b_k\mid s_k)\Big) \mid s_0\sim\rho\right]\notag\\
    &\hspace{20pt}-\mathbb{E}\left[\sum_{k=0}^{\infty} \gamma^k \Big(r\left(s_k, a_k, b_k\right)-\tau'\log\pi(a_k\mid s_k)+\tau'\log\phi(b_k\mid s_k)\Big) \mid s_0\sim\rho\right]\notag\\
    &=\mathbb{E}\left[\sum_{k=0}^{\infty} \gamma^k \Big((\tau-\tau')\log\pi(a_k\mid s_k)+(\tau-\tau')\log\phi(b_k\mid s_k)\Big) \mid s_0\sim\rho\right]\notag\\
    &=\frac{\tau-\tau'}{1-\gamma}\mathbb{E}_{s'\sim d_{\rho}^{\pi,\phi},a\sim \pi(\cdot\mid s'), b\sim \phi(\cdot\mid s')}\left[-\log\pi(a\mid s')+\log\phi(b\mid s')\right]\notag\\
    &=\frac{\tau-\tau'}{1-\gamma}\mathbb{E}_{s'\sim d_{\rho}^{\pi,\phi}}[H(\pi(\cdot\mid s'))-H(\phi(\cdot\mid s'))],
\end{align*}
where $H$ denotes the entropy and is defined in \eqref{eq:H_def}.

We have the following upper and lower bound on the entropy
\begin{align*}
    0 \leq H(\pi(\cdot\mid s'))\leq \log|\Acal|, \quad 0 \leq H(\phi(\cdot\mid s'))\leq \log|\Bcal|.
\end{align*}

Therefore, if $\tau\geq\tau'\geq 0$,
\begin{align*}
    -\frac{\tau-\tau'}{1-\gamma}\log|\Bcal|\leq J_{\tau}(\pi,\phi)-J_{\tau'}(\pi,\phi)\leq \frac{\tau-\tau'}{1-\gamma}\log|\Acal|.
\end{align*}

For any $\tau\geq\tau'\geq0$,
\begin{align*}
    &J_{\tau}(\pi_{\tau}^{\star},\phi_{\tau}^{\star})-J_{\tau'}(\pi_{\tau'}^{\star},\phi_{\tau'}^{\star})\notag\\
    &=\max_{\pi}\min_{\phi}J_{\tau}(\pi,\phi)-\min_{\phi}J_{\tau'}(\pi_{\tau'}^{\star},\phi)\\
    &\geq \min_{\phi}J_{\tau}(\pi_{\tau'}^{\star},\phi)-\min_{\phi}J_{\tau'}(\pi_{\tau'}^{\star},\phi)\\
    &= \min_{\phi}\Big(J_{\tau'}(\pi_{\tau'}^{\star},\phi)+(\tau-\tau')\Hcal_{\pi}(\rho,\pi_{\tau'}^{\star},\phi)-(\tau-\tau')\Hcal_{\phi}(\rho,\pi_{\tau'}^{\star},\phi)\Big)-\min_{\phi}J_{\tau'}(\pi_{\tau'}^{\star},\phi)\\
    &\geq \min_{\phi}J_{\tau'}(\pi_{\tau'}^{\star},\phi)+(\tau\hspace{-2pt}-\hspace{-2pt}\tau')\min_{\phi}\Hcal_{\pi}(\rho,\pi_{\tau'}^{\star},\phi)+(\tau-\tau')\min_{\phi}-\Hcal_{\phi}(\rho,\pi_{\tau'}^{\star},\phi)\hspace{-2pt}-\hspace{-2pt}\min_{\phi}J_{\tau'}(\pi_{\tau'}^{\star},\phi)\\
    &=(\tau-\tau')\left(\min_{\phi}\Hcal_{\pi}(\rho,\pi_{\tau'}^{\star},\phi)-\max_{\phi}\Hcal_{\phi}(\rho,\pi_{\tau'}^{\star},\phi)\right)\\
    &\geq(\tau-\tau')(0-\log|\Bcal|)\\
    &=-(\tau-\tau')\log|\Bcal|,
\end{align*}
where the second inequality comes from the fact that $\min_{x}f_1(x)+f_2(x)\geq \min_{x}f_1(x)+\min_{x}f_2(x)$ for any functions $f_1,f_2$ of the same domain.

It can be shown by a similar argument
\[J_{\tau}(\pi_{\tau}^{\star},\phi_{\tau}^{\star})-J_{\tau'}(\pi_{\tau'}^{\star},\phi_{\tau'}^{\star}) \leq (\tau-\tau')\log|\Acal|.\]

In addition, for any $\tau\geq\tau'\geq0$ and any policy $\pi$,
\begin{align*}
    &J_{\tau}(\pi,\phi_{\tau}(\pi))-J_{\tau'}(\pi,\phi_{0}(\pi))\notag\\
    &=\min_{\phi}J_{\tau}(\pi,\phi)-\min_{\phi}J_{\tau'}(\pi,\phi)\\
    &= \min_{\phi}\left(J_{\tau'}(\pi,\phi)+(\tau-\tau')\Hcal_{\pi}(\rho,\pi,\phi)-(\tau-\tau')\Hcal_{\phi}(\rho,\pi,\phi)\right)-\min_{\phi}J_{\tau'}(\pi,\phi)\\
    &\leq \left(\min_{\phi}J_{\tau'}(\pi,\phi)+(\tau-\tau')\max_{\phi}\Hcal_{\pi}(\rho,\pi,\phi)+(\tau-\tau')\max_{\phi}(-\Hcal_{\phi}(\rho,\pi,\phi))\right)-\min_{\phi}J_{\tau'}(\pi,\phi)\\
    &=(\tau-\tau')\left(\max_{\phi}\Hcal_{\pi}(\rho,\pi,\phi)-\min_{\phi}\Hcal_{\phi}(\rho,\pi,\phi)\right)\\
    &\leq(\tau-\tau')\log|\Acal|.
\end{align*}

It can be shown by a similar argument
\[J_{\tau}(\pi,\phi_{\tau}(\pi))-J_{\tau'}(\pi,\phi_{0}(\pi))\geq -(\tau-\tau')\log|\Bcal|.\]

\qed

\subsection{Proof of Lemma \ref{lem:nonuniform_PL}}

Adapting \citet{mei2020global}[Lemma 15], we have for any $\theta\in\mathbb{R}^{\Scal\times\Acal}$ and $\psi\in\mathbb{R}^{\Scal\times\Bcal}$
\begin{align*}
    &\|\nabla_{\theta} J_{\tau}(\pi_{\theta},\phi_{\psi})\|^2\geq\frac{2\tau\rho_{\min}}{|\Scal|}\left(\min_{s,a}\pi_{\theta}(a\mid s)\right)^2\hspace{-2pt}\left\|\frac{d_{\rho}^{\pi_{\tau}(\phi_{\psi}),\phi_{\psi}}}{d_{\rho}^{\pi_{\theta},\phi_{\psi}}}\right\|_{\infty}^{-1}\left(J_{\tau}(\pi_{\tau}(\phi_{\psi}),\phi_{\psi})\hspace{-2pt}-\hspace{-2pt}J_{\tau}(\pi_{\theta},\phi_{\psi})\right),\\
    &\|\nabla_{\psi} J_{\tau}(\pi_{\theta},\phi_{\psi})\|^2\geq\frac{2\tau\rho_{\min}}{|\Scal|}\left(\min_{s,b}\phi_{\psi}(b\mid s)\right)^2\left\|\frac{d_{\rho}^{\pi_{\theta},\phi_{\tau}(\pi_{\theta})}}{d_{\rho}^{\pi_{\theta},\phi_{\psi}}}\right\|_{\infty}^{-1}\left(J_{\tau}(\pi_{\theta},\phi_{\psi})-J_{\tau}(\pi_{\theta},\phi_{\tau}(\pi_{\theta}))\right).
\end{align*}

Then, the first inequality follows from $d_{\rho}^{\pi_{\tau}(\phi_{\psi}),\phi_{\psi}}(s)\leq 1$ and $d_{\rho}^{\pi_{\theta},\phi_{\psi}}(s)\geq(1-\gamma)\rho(s)\geq(1-\gamma)\rho_{\min}$ for all $s\in\Scal$, and the second inequality from $d_{\rho}^{\pi_{\theta},\phi_{\tau}(\pi_{\theta})}\leq1$ and $d_{\rho}^{\pi_{\theta},\phi_{\psi}}\geq(1-\gamma)\rho_{\min}$ for all $s\in\Scal$.

\qed

\subsection{Proof of Lemma \ref{lem:value_LipschitzGrad}}

\citet{mei2020global}[Lemma 7, Lemma 14] establishes the smoothness condition of the value function and the regularization entropy with respect to one player's policy, i.e.
\begin{align*}
    \|\nabla_{\theta}J(\pi_{\theta_1},\phi_{\psi_1})-\nabla_{\theta}J(\pi_{\theta_2},\phi_{\psi_1})\|&\leq L_V\|\theta_1-\theta_2\|,\\
    \|\nabla_{\psi}J(\pi_{\theta_1},\phi_{\psi_1})-\nabla_{\psi}J(\pi_{\theta_1},\phi_{\psi_2})\|&\leq L_V\|\psi_1-\psi_2\|.
\end{align*}

Therefore, we only need to show
\begin{align*}
    \|\nabla_{\theta}J(\pi_{\theta_1},\phi_{\psi_1})-\nabla_{\theta}J(\pi_{\theta_1},\phi_{\psi_2})\|&\leq L_V\|\psi_1-\psi_2\|,\\
    \|\nabla_{\psi}J(\pi_{\theta_1},\phi_{\psi})-\nabla_{\psi}J(\pi_{\theta_2},\phi_{\psi})\|&\leq L_V\|\theta_1-\theta_2\|.
\end{align*}

Given a fixed $\theta$ and $\psi$, with
arbitrary vectors $u$ and $v$ such that $\|u\|_2=\|v\|_2=1$, we define the shorthand notation
\[\pi_{\alpha,u}=\pi_{\theta+\alpha u},\quad\phi_{\beta,v}=\pi_{\psi+\beta v}.\]


According to \citet{zeng2021decentralized}[Lemma B.5], 
\begin{align*}
    &\sum_a\left|\frac{d \pi_{\alpha,u}(a\mid s)}{d\alpha}\right|\leq 2,\quad\sum_b\left|\frac{d \phi_{\beta,v}(b\mid s)}{d\beta}\right|\leq 2,\\
    &\sum_{a,b}\left|\frac{d \pi_{\alpha}(a \mid s)}{d \alpha}\frac{d \phi_{\beta,v}(b \mid s)}{d \beta}\right|\leq\left(\sum_{a}\left|\frac{d \pi_{\alpha}(a \mid s)}{d \alpha}\right|\right)\left(\sum_{b}\left|\frac{d \phi_{\beta}(b \mid s)}{d \beta}\right|\right)\leq 4.
\end{align*}

Let $P(\alpha,\beta,u,v)\in\mathbb{R}^{|\Scal||\Acal||\Bcal|\times|\Scal||\Acal||\Bcal|}$ denote the state-action transition matrix induced by the policy pair $(\pi_{\alpha,u},\phi_{\beta,v})$
\[P(\alpha,\beta,u,v)_{(s,a,b)\rightarrow(s',a',b')}=\Pcal(s'\mid s,a,b)\pi_{\alpha,u}(a'\mid s')\phi_{\beta,v}(b'\mid s').\]

Differentiating with respect to $\alpha$ and $\beta$,
\begin{align*}
    \left[\frac{d^2 P(\alpha,\beta,u,v)}{d \alpha d \beta}\right]_{(s, a,b)\rightarrow(s',a',b')}=\frac{d \pi_{\alpha,u}(a' \mid s')}{d \alpha} \frac{d \phi_{\beta,v}(b' \mid s')}{d \beta} \Pcal(s' \mid s, a, b),
\end{align*}
which implies for any vector $x$
\begin{align*}
    \left[\frac{d^2 P(\alpha,\beta,u,v)}{d \alpha d \beta}x\right]_{s, a,b}=\sum_{s',a',b'}\frac{d \pi_{\alpha}(a' \mid s')}{d \alpha}\frac{d \phi_{\beta,v}(b' \mid s')}{d \beta}\Pcal(s' \mid s, a, b)x_{s',a',b'}.
\end{align*}

The $\ell_{\infty}$ norm of this quantity can be upper bounded
\begin{align}
    &\max_{\|u\|_2=\|v\|_2=1}\left\|\frac{d^2 P(\alpha,\beta,u,v)}{d \alpha d \beta}x\right\|_{\infty}\notag\\
    &=\max_{s,a,b}\max_{\|u\|_2=\|v\|_2=1}\left|\left[\frac{d^2 P(\alpha,\beta,u,v)}{d \alpha d \beta}x\right]_{s, a,b}\right|\notag\\
    &=\max_{s,a,b}\max_{\|u\|_2=\|v\|_2=1}\left|\sum_{s',a',b'}\frac{d \pi_{\alpha}(a' \mid s')}{d \alpha} \frac{d \phi_{\beta,v}(b' \mid s')}{d \beta}\Pcal(s' \mid s, a, b)x_{s',a',b'}\right|\notag\\
    &\leq \max_{s,a,b}\sum_{s'}\Pcal(s' \mid s, a, b)\|x\|_{\infty}\max_{\|u\|_2=\|v\|_2=1}\sum_{a',b'}\left|\frac{d \pi_{\alpha}(a' \mid s')}{d \alpha}\frac{d \phi_{\beta,v}(b' \mid s')}{d \beta}\right|\notag\\
    &\leq 4\|x\|_{\infty}.
    \label{lem:LipschitzGrad:proof_eq1}
\end{align}

Using an identical argument, we can show that 
\begin{align}
    &\max_{\|u\|_2=\|v\|_2=1}\left\|\frac{d P(\alpha,\beta,u,v)}{d \alpha}x\right\|_{\infty}\leq\sum_a\left|\frac{d \pi_{\alpha,u}(a\mid s)}{d\alpha}\right|\|x\|_{\infty}\leq 2\|x\|_{\infty},\\
    &\max_{\|u\|_2=\|v\|_2=1}\left\|\frac{d P(\alpha,\beta,u,v)}{d \beta}x\right\|_{\infty}\leq\sum_b\left|\frac{d \pi_{\beta,v}(b\mid s)}{d\beta}\right|\|x\|_{\infty}\leq 2\|x\|_{\infty}.\label{lem:LipschitzGrad:proof_eq2}
\end{align}

With $M(\alpha,\beta,u,v)=(I-\gamma P(\alpha,\beta,u,v))^{-1}$ and $r=[r(s_0,a_0,b_0),\cdots, r(s_{|\Scal|},a_{|\Acal|},b_{|\Bcal|})]$,
\[Q^{\pi_{\alpha,u},\phi_{\beta,v}}(s,a,b)=e_{s,a,b}^{\top}M(\alpha,\beta,u,v)r.\]

Taking the derivatives,
\begin{align*}
    &\frac{d Q^{\pi_{\alpha,u},\phi_{\beta,v}}(s,a,b)}{d\alpha}=\gamma e_{s,a,b}^{\top}M(\alpha,\beta,u,v)\frac{d P(\alpha,\beta,u,v)}{d \alpha}M(\alpha,\beta,u,v)r,\\
    &\frac{d Q^{\pi_{\alpha,u},\phi_{\beta,v}}(s,a,b)}{d\beta}=\gamma e_{s,a,b}^{\top}M(\alpha,\beta,u,v)\frac{d P(\alpha,\beta,u,v)}{d \beta}M(\alpha,\beta,u,v)r.
\end{align*}
Taking the second-order derivative,
\begin{align*}
    &\frac{d^2 Q^{\pi_{\alpha,u},\phi_{\beta,v}}(s,a,b)}{d\alpha d\beta}\notag\\
    &=\gamma^2 e_{s,a,b}^{\top}M(\alpha,\beta,u,v)\frac{d P(\alpha,\beta,u,v)}{d \alpha}M(\alpha,\beta,u,v)\frac{d P(\alpha,\beta,u,v)}{d \beta}M(\alpha,\beta,u,v)r\\
    &\hspace{20pt}+\gamma^2 e_{s,a,b}^{\top}M(\alpha,\beta,u,v)\frac{d P(\alpha,\beta,u,v)}{d \beta}M(\alpha,\beta,u,v)\frac{d P(\alpha,\beta,u,v)}{d \alpha}M(\alpha,\beta,u,v)r\\
    &\hspace{20pt}+\gamma e_{s,a,b}^{\top}M(\alpha,\beta,u,v)\frac{d^2 P(\alpha,\beta,u,v)}{d \alpha d \beta}M(\alpha,\beta,u,v)r
\end{align*}

Using $M(\alpha,\beta,u,v)\1=(I-\gamma P(\alpha,\beta,u,v))^{-1}\1=\frac{1}{1-\gamma}\1$ and inequalities \eqref{lem:LipschitzGrad:proof_eq1} and \eqref{lem:LipschitzGrad:proof_eq2}, we have
\begin{align*}
    &\max_{\|u\|_2=\|v\|_2=1}\hspace{-2pt}\left|\frac{d Q^{\pi_{\alpha,u},\phi_{\beta,v}}(s,a,b)}{d\alpha}\right|\hspace{-2pt}\leq \hspace{-2pt}\|\gamma M(\alpha,\beta,u,v)\frac{d P(\alpha,\beta,u,v)}{d \alpha}M(\alpha,\beta,u,v)r\|_{\infty}\hspace{-2pt}\leq \hspace{-2pt}\frac{2\gamma }{(1-\gamma)^2},\\
    &\max_{\|u\|_2=\|v\|_2=1}\hspace{-2pt}\left|\frac{d Q^{\pi_{\alpha,u},\phi_{\beta,v}}(s,a,b)}{d\beta}\right|\hspace{-2pt}\leq\hspace{-2pt} \|\gamma M(\alpha,\beta,u,v)\frac{d P(\alpha,\beta,u,v)}{d \beta}M(\alpha,\beta,u,v)r\|_{\infty}\hspace{-2pt}\leq\hspace{-2pt} \frac{2\gamma}{(1-\gamma)^2},
\end{align*}
and
\begin{align*}
    &\max_{\|u\|_2=\|v\|_2=1}\left|\frac{d^2 Q^{\pi_{\alpha,u},\phi_{\beta,v}}(s,a,b)}{d\alpha d\beta}\right|\\
    &\leq \|\gamma^2 M(\alpha,\beta,u,v)\frac{d P(\alpha,\beta,u,v)}{d \alpha}M(\alpha,\beta,u,v)\frac{d P(\alpha,\beta,u,v)}{d \beta}M(\alpha,\beta,u,v)r\|_{\infty}\\
    &\hspace{20pt}+\|\gamma^2 M(\alpha,\beta,u,v)\frac{d P(\alpha,\beta,u,v)}{d \beta}M(\alpha,\beta,u,v)\frac{d P(\alpha,\beta,u,v)}{d \alpha}M(\alpha,\beta,u,v)r\|_{\infty}\\
    &\hspace{20pt}+\|\gamma M(\alpha,\beta,u,v)\frac{d^2 P(\alpha,\beta,u,v)}{d \alpha d \beta}M(\alpha,\beta,u,v)r\|_{\infty}\\
    &\leq\frac{2\gamma^2}{(1-\gamma)^3}+\frac{4\gamma}{(1-\gamma)^2}.
\end{align*}

Since $V^{\pi_{\alpha,u},\phi_{\beta,v}}(s)=\sum_{a,b}\pi_{\alpha,u}(a\mid s)\phi_{\beta,v}(b\mid s)Q^{\pi_{\alpha,u},\phi_{\beta,v}}(s,a,b)$,
\begin{align*}
    \frac{d^2 V^{\pi_{\alpha,u},\phi_{\beta,v}}(s)}{d\alpha d\beta}&=\sum_{a,b}\frac{d\pi_{\alpha,u}(a\mid s)}{d\alpha}\frac{d\phi_{\beta,v}(b\mid s)}{d\beta}Q^{\pi_{\alpha,u},\phi_{\beta,v}}(s,a,b)\\
    &\hspace{20pt}+\sum_{a,b}\pi_{\alpha,u}(a\mid s)\phi_{\beta,v}(b\mid s)\frac{d^2 Q^{\pi_{\alpha,u},\phi_{\beta,v}}(s,a,b)}{d\alpha d\beta}\\
    &\hspace{20pt}+\sum_{a,b}\frac{d\pi_{\alpha,u}(a\mid s)}{d\alpha}\phi_{\beta,v}(b\mid s)\frac{d Q^{\pi_{\alpha,u},\phi_{\beta,v}}(s,a,b)}{d\beta}\\
    &\hspace{20pt}+\sum_{a,b}\pi_{\alpha,u}(a\mid s)\frac{d \phi_{\beta,v}(b\mid s)}{d\beta}\frac{d Q^{\pi_{\alpha,u},\phi_{\beta,v}}(s,a,b)}{d\alpha}.
\end{align*}

Therefore,
\begin{align*}
    \max_{\|u\|_2=\|v\|_2=1}\left|\frac{d V^{\pi_{\alpha,u},\phi_{\beta,v}}(s)}{d\alpha d\beta}\right|\leq\frac{4}{1-\gamma}+\left(\frac{2\gamma^2}{(1-\gamma)^3}+\frac{4\gamma}{(1-\gamma)^2}\right)+2\frac{4\gamma}{(1-\gamma)^2}\leq\frac{8}{(1-\gamma)^3},
\end{align*}
which implies
\begin{align*}
    \|\nabla_{\theta}J(\pi_{\theta},\phi_{\psi_1})-\nabla_{\theta}J(\pi_{\theta},\phi_{\psi_2})\|\leq \frac{8}{(1-\gamma)^3}\|\psi_1-\psi_2\|.
\end{align*}
Similarly, it follows by the same argument that
\begin{align*}
    \|\nabla_{\psi}J(\pi_{\theta_1},\phi_{\psi})-\nabla_{\psi}J(\pi_{\theta_2},\phi_{\psi})\|\leq \frac{8}{(1-\gamma)^3}\|\theta_1-\theta_2\|.
\end{align*}

\citet{zeng2021decentralized}[Lemma B.5] implies
\begin{align}
    \|J(\pi_{\theta_1},\phi_{\psi_1})-J(\pi_{\theta_2},\phi_{\psi_2})\|\leq \frac{2}{(1-\gamma)^2}(\|\theta_1-\theta_2\|+\|\psi_1-\psi_2\|),\label{lem:LipschitzGrad:proof_eq3}
\end{align}
and we simply use $\frac{2}{(1-\gamma)^2}\leq L_V$.

\qed

\subsection{Proof of Lemma \ref{lem:reg_LipschitzGrad}}

We will prove the first two inequalities on the Lipschitz gradient of $\Hcal_{\pi}$. The next two inequalities are completely symmetric and can be derived using an identical argument.

\citet{mei2020global}[Lemma 14] implies
\begin{align*}
    \|\nabla_{\theta}\Hcal_{\pi}(s,\pi_{\theta_1},\phi_{\psi_1})-\nabla_{\theta}\Hcal_{\pi}(s,\pi_{\theta_2},\phi_{\psi_1})\|\leq L_{\Hcal}\|\theta_1-\theta_2\|,
\end{align*}
so we just need to show
\begin{align}
    \begin{aligned}
    \|\nabla_{\theta}\Hcal_{\pi}(s,\pi_{\theta_1},\phi_{\psi_1})-\nabla_{\theta}\Hcal_{\pi}(s,\pi_{\theta_1},\phi_{\psi_2})\|\leq L_{\Hcal}\|\psi_1-\psi_2\|,\\
    \|\nabla_{\psi}\Hcal_{\pi}(s,\pi_{\theta_1},\phi_{\psi_1})-\nabla_{\psi}\Hcal_{\pi}(s,\pi_{\theta_2},\phi_{\psi_1})\|\leq L_{\Hcal}\|\theta_1-\theta_2\|,\\
    \|\nabla_{\psi}\Hcal_{\pi}(s,\pi_{\theta_1},\phi_{\psi_1})-\nabla_{\psi}\Hcal_{\pi}(s,\pi_{\theta_1},\phi_{\psi_2})\|\leq L_{\Hcal}\|\psi_1-\psi_2\|.
    \end{aligned}
    \label{lem:reg_LipschitzGrad:proof_eq1}
\end{align}

Given a fixed $\theta$ and $\psi$, with
arbitrary vectors $u$ and $v$ such that $\|u\|_2=\|v\|_2=1$, we define the shorthand notation
\[\pi_{\alpha,u}=\pi_{\theta+\alpha u},\quad\phi_{\beta,v}=\pi_{\psi+\beta v}.\]

Note that to show \eqref{lem:reg_LipschitzGrad:proof_eq1}, it suffices to show for any $u,v$
\begin{align*}
    \left|\frac{d^2\Hcal_\pi(s,\pi_{\alpha,u},\phi_{\beta,v})}{d\alpha d\beta}\right|\leq L_{\Hcal},\quad\left|\frac{d^2\Hcal_\pi(s,\pi_{\alpha,u},\phi_{\beta,v})}{d\beta^2}\right|\leq L_{\Hcal}.
\end{align*}

We define the state transition matrix $P\in\mathbb{R}^{|\Scal|\times|\Scal|}$ such that
\[P(\alpha,\beta,u,v)_{s\rightarrow s'}=\sum_{a,b}\Pcal(s'\mid s,a,b)\pi_{\alpha,u}(a\mid s)\phi_{\beta,v}(b\mid s).\]
Let $M(\alpha,\beta,u,v)=(I-\gamma P(\alpha,\beta,u,v))^{-1}$. Then, we can re-write $\Hcal_\pi(s,\pi,\phi)$ in the matrix form
\begin{align*}
    \Hcal_\pi(s,\pi,\phi)=e_{s}^{\top}M(\alpha,\beta,u,v)h_{\alpha,u},
\end{align*}
where $h_{\alpha,u}=[h_{\alpha,u}(s_0),\cdots,h_{\alpha,u}(s_{|\Scal|})]\in\mathbb{R}^{|\Scal|}$ is a vector with
\[h_{\alpha,u}(s)=-\sum_{a}\pi_{\alpha,u}(a\mid s) \log\pi_{\alpha,u}(a\mid s).\]

According to \citet{mei2020global}[Lemma 14],
\[\left\|\frac{d h_{\alpha,u}}{d\alpha}\right\|_{\infty}\leq2\log|\Acal|\|u\|_2
=2\log|\Acal|.\]

Taking the derivatives of $\Hcal_\pi(s,\pi,\phi)$,
\begin{align*}
    &\frac{d\Hcal_\pi(s,\pi_{\alpha,u},\phi_{\beta,v})}{d\alpha}\notag\\
    &=\gamma e_{s}^{\top} M(\alpha,\beta,u,v)\frac{d P(\alpha,\beta,u,v)}{d\alpha}M(\alpha,\beta,u,v)h_{\alpha,u}+e_{s}^{\top}M(\alpha,\beta,u,v)\frac{d h_{\alpha,u}}{d\alpha},
\end{align*}
and taking second order derivative
\begin{align*}
    &\frac{d^2\Hcal_\pi(s,\pi_{\alpha,u},\phi_{\beta,v})}{d\alpha d\beta}\notag\\
    &=\gamma^2 e_{s}^{\top} M(\alpha,\beta,u,v)\frac{d P(\alpha,\beta,u,v)}{d\alpha}M(\alpha,\beta,u,v)\frac{d P(\alpha,\beta,u,v)}{d\beta}M(\alpha,\beta,u,v)h_{\alpha,u}\\
    &\hspace{20pt}+\gamma^2 e_{s}^{\top} M(\alpha,\beta,u,v)\frac{d P(\alpha,\beta,u,v)}{d\beta}M(\alpha,\beta,u,v)\frac{d P(\alpha,\beta,u,v)}{d\alpha}M(\alpha,\beta,u,v)h_{\alpha,u}\\
    &\hspace{20pt}+\gamma e_{s}^{\top} M(\alpha,\beta,u,v)\frac{d^2 P(\alpha,\beta,u,v)}{d\alpha d\beta}M(\alpha,\beta,u,v)h_{\alpha,u}\\
    &\hspace{20pt}+\gamma e_{s}^{\top} M(\alpha,\beta,u,v)\frac{d P(\alpha,\beta,u,v)}{d\beta}M(\alpha,\beta,u,v)\frac{d h_{\alpha,u}}{d\alpha}.
\end{align*}

Using a similar line of argument to \citet{mei2020global}[Eq.~(192)-(195)] and analysis in Lemma \ref{lem:value_LipschitzGrad} of our work, we can show that for any vector $x$
\begin{align*}
    \left\|\frac{d P(\alpha,\beta,u,v)}{d\alpha}x\right\|_{\infty}\leq2\|x\|_{\infty},\,\,\left\|\frac{d P(\alpha,\beta,u,v)}{d\beta}\right\|_{\infty}\leq2\|x\|_{\infty},\,\,\left\|\frac{d^2 P(\alpha,\beta,u,v)}{d\alpha d\beta}\right\|_{\infty}\leq4\|x\|_{\infty}.
\end{align*}

From the fact that $\|M(\alpha,\beta,u,v)x\|_{\infty}\leq\frac{1}{1-\gamma}\|x\|_{\infty}$, we have for any vectors $u,v$
\begin{align*}
    &\left|\frac{d^2\Hcal_\pi(s,\pi_{\alpha,u},\phi_{\beta,v})}{d\alpha d\beta}\right|\\
    &\leq\gamma^2 \left\| M(\alpha,\beta,u,v)\frac{d P(\alpha,\beta,u,v)}{d\alpha}M(\alpha,\beta,u,v)\frac{d P(\alpha,\beta,u,v)}{d\beta}M(\alpha,\beta,u,v)h_{\alpha,u}\right\|\\
    &\hspace{20pt}+\gamma^2 \left\| M(\alpha,\beta,u,v)\frac{d P(\alpha,\beta,u,v)}{d\beta}M(\alpha,\beta,u,v)\frac{d P(\alpha,\beta,u,v)}{d\alpha}M(\alpha,\beta,u,v)h_{\alpha,u}\right\|\\
    &\hspace{20pt}+\gamma \left\| M(\alpha,\beta,u,v)\frac{d^2 P(\alpha,\beta,u,v)}{d\alpha d\beta}M(\alpha,\beta,u,v)h_{\alpha,u}\right\|\\
    &\hspace{20pt}+\gamma \left\| M(\alpha,\beta,u,v)\frac{d P(\alpha,\beta,u,v)}{d\beta}M(\alpha,\beta,u,v)\frac{d h_{\alpha,u}}{d\alpha}\right\|\\
    &\leq \frac{4\gamma^2 \log|\Acal|}{(1-\gamma)^3}+\frac{4\gamma^2 \log|\Acal|}{(1-\gamma)^3}+\frac{4\gamma\log|\Acal|}{(1-\gamma)^2}+\frac{2\gamma}{(1-\gamma)^2}\cdot2\log|\Acal|\\
    &\leq\frac{8\log|\Acal|}{(1-\gamma)^3}.
\end{align*}

Now it remains to be shown
\[\left|\frac{d^2\Hcal_\pi(s,\pi_{\alpha,u},\phi_{\beta,v})}{d\beta^2}\right|\leq L_{\Hcal}.\]
From the eye of the second player, $\Hcal_{\pi}(s,\pi_{\theta},\phi_{\psi})$ is simply the value function of a regular MDP with itself as the only agent (the first player's policy combines with $\Pcal$) with the reward function $r(s,b)=-\sum_{a\in\Acal}\pi_{\theta}(a\mid s)\log\pi_{\theta}(a\mid s)\in[0,\log|\Acal|]$. Therefore, by Lemma \ref{lem:value_LipschitzGrad} which is derived with reward bounded between 0 and 1, we know \[\left|\frac{d^2\Hcal_\pi(s,\pi_{\alpha,u},\phi_{\beta,v})}{d\beta^2}\right|\leq \log|\Acal|L_V\leq L_{\Hcal}.\]

To show the Lipschitz continuity, we note that
\begin{align*}
    &\left|\frac{d\Hcal_\pi(s,\pi_{\alpha,u},\phi_{\beta,v})}{d\alpha}\right|\\
    &=\left|\gamma e_{s}^{\top} M(\alpha,\beta,u,v)\frac{d P(\alpha,\beta,u,v)}{d\alpha}M(\alpha,\beta,u,v)h_{\alpha,u}+e_{s}^{\top}M(\alpha,\beta,u,v)\frac{d h_{\alpha,u}}{d\alpha}\right|\\
    &\leq \gamma\|M(\alpha,\beta,u,v)\frac{d P(\alpha,\beta,u,v)}{d\alpha}M(\alpha,\beta,u,v)h_{\alpha,u}\|+\|M(\alpha,\beta,u,v)\frac{d h_{\alpha,u}}{d\alpha}\|\\
    &\leq \frac{4\gamma\log|\Acal|}{(1-\gamma)^2}+\frac{2\log|\Acal|}{1-\gamma}\leq L_{\Hcal}.
\end{align*}

To show the Lipschitz continuity of $\Hcal_{\pi}$ with respect to $\psi$, we use the same argument as above and note that from the eye of the second player, $\Hcal_{\pi}(s,\pi_{\theta},\phi_{\psi})$ is simply the value function of a regular MDP with itself as the only agent (the first player's policy combines with $\Pcal$) with the reward function $r(s,b)=-\sum_{a\in\Acal}\pi_{\theta}(a\mid s)\log\pi_{\theta}(a\mid s)\in[0,\log|\Acal|]$.
Adapting \eqref{lem:LipschitzGrad:proof_eq3}, we have
\begin{align*}
    \left|\frac{d\Hcal_\pi(s,\pi_{\alpha,u},\phi_{\beta,v})}{d\beta}\right|\leq\frac{2}{(1-\gamma)^2}\cdot\log|\Acal|\leq L_{\Hcal}.
\end{align*}

\qed

\subsection{Proof of Lemma \ref{lem:stepsize_diff}}

We first show that for any $\tilde{k}>0$, we have $\frac{1}{\tilde{k}^a}-\frac{1}{(\tilde{k}+1)^a}\leq\frac{8}{3(\tilde{k}+1)^{a+1}}$.

Since the integer $\tilde{k}$ is positive, it can be lower bound by $\frac{\tilde{k}+1}{2}$.
\begin{align*}
    &\frac{1}{\tilde{k}^a}-\frac{1}{(\tilde{k}+1)^a}\\
    &=\frac{(\tilde{k}+1)^a-\tilde{k}^a}{\tilde{k}^a (\tilde{k}+1)^a}\leq \frac{2((\tilde{k}+1)^a-\tilde{k}^a)}{(\tilde{k}+1)^{2a}}=\frac{2((\tilde{k}+1)^a-\tilde{k}^a)\left((\tilde{k}+1)^{1-a}+\tilde{k}^{1-a}\right)}{(\tilde{k}+1)^{2a}\left((\tilde{k}+1)^{1-a}+\tilde{k}^{1-a}\right)}\\
    &\leq \frac{2((\tilde{k}+1)^a-\tilde{k}^a)\left((\tilde{k}+1)^{1-a}+\tilde{k}^{1-a}\right)}{(\tilde{k}+1)^{2a}\left((\tilde{k}+1)^{1-a}+\frac{1}{2}(\tilde{k}+1)^{1-a}\right)}= \frac{4((\tilde{k}+1)^a-\tilde{k}^a)\left((\tilde{k}+1)^{1-a}+\tilde{k}^{1-a}\right)}{3(\tilde{k}+1)^{a+1}}\\
    &=\frac{4\left((\tilde{k}+1)-\tilde{k}^a(\tilde{k}+1)^{1-a}+\tilde{k}^{1-a}(\tilde{k}+1)^a-\tilde{k}\right)}{3(\tilde{k}+1)^{a+1}}\\
    &=\frac{4\left(1-\tilde{k}^a(\tilde{k}+1)^{1-a}+\tilde{k}^{1-a}(\tilde{k}+1)^a\right)}{3(\tilde{k}+1)^{a+1}} \leq \frac{8}{3(\tilde{k}+1)^{a+1}},
\end{align*}
where the last inequality follows from
\begin{align*}
    \tilde{k}^{1-a}(\tilde{k}+1)^a-\tilde{k}^a(\tilde{k}+1)^{1-a}\leq (\tilde{k}+1)^{1-a}(\tilde{k}+1)^a-\tilde{k}^a \tilde{k}^{1-a} = \tilde{k}+1-\tilde{k}=1.
\end{align*}

Choosing $\tilde{k}=k+h$ yields
\begin{align*}
    \frac{1}{(k+h)^a}-\frac{1}{(k+1+h)^a}\leq\frac{8}{3(k+1+h)^{a+1}}\leq\frac{8}{3(k+h)^{a+1}}.
\end{align*}

\qed

\subsection{Proof of Lemma \ref{lem:min_policy_iterates_constanttau}}

The property of the min and max function implies that
\begin{align*}
    \max_{s,a}(\pi_{\tau}^{\star}(a\mid s)-\pi_{\theta_k}(a\mid s))+\min_{s,a}\pi_{\theta_k}(a\mid s)\geq\min_{s,a}\pi_{\tau}^{\star}(a\mid s).
\end{align*}
Since the three terms are all non-negative, the inequality holds after taking the square
\begin{align*}
    (\min_{s,a}\pi_{\tau}^{\star}(a\mid s))^2&\leq(\max_{s,a}(\pi_{\tau}^{\star}(a\mid s)-\pi_{\theta_k}(a\mid s))+\min_{s,a}\pi_{\theta_k}(a\mid s))^2\notag\\
    &\leq \frac{4}{3}(\min_{s,a}\pi_{\theta_k}(a\mid s))^2+4(\max_{s,a}(\pi_{\tau}^{\star}(a\mid s)-\pi_{\theta_k}(a\mid s)))^2.
\end{align*}
Re-arranging the terms,
\begin{align*}
    -\left(\min_{s,a}\pi_{\theta_k}(a\mid s)\right)^2&\leq-\frac{3}{4}\left(\min_{s,a}\pi_{\tau}^{\star}(a\mid s)\right)^2+3\left(\max_{s,a}\pi_{\tau}^{\star}(a\mid s)-\pi_{\phi_k}(a\mid s)\right)^2\notag\\
    &\leq-\frac{3}{4}\left(\min_{s,a}\pi_{\tau}^{\star}(a\mid s)\right)^2+3\|\pi_{\tau}^{\star}-\pi_{\phi_k}\|^2
\end{align*}

From Lemma \ref{lem:quadratic_growth},
\begin{align}
    -\left(\min_{s,a}\pi_{\theta_k}(a\mid s)\right)^2&\leq-\frac{3}{4}\left(\min_{s,a}\pi_{\tau}^{\star}(a\mid s)\right)^2+3\|\pi_{\tau}^{\star}-\pi_{\phi_k}\|^2\notag\\
    &\leq-\frac{3}{4}\left(\min_{s,a}\pi_{\tau}^{\star}(a\mid s)\right)^2+\frac{6\log(2)}{\tau\rho_{\min}}(J_{\tau}(\pi_{\tau}^{\star},\phi_{\tau}^{\star})-J_{\tau}(\pi_{\theta_k},\phi_{\tau}^{\star}))\notag\\
    &\leq-\frac{3}{4}\left(\min_{s,a}\pi_{\tau}^{\star}(a\mid s)\right)^2+\frac{6\log(2)}{\tau\rho_{\min}}(J_{\tau}(\pi_{\tau}^{\star},\phi_{\tau}^{\star})-J_{\tau}(\pi_{\theta_k},\phi_{\tau}(\pi_{\theta_k})))\notag\\
    &=-\frac{3}{4}\left(\min_{s,a}\pi_{\tau}^{\star}(a\mid s)\right)^2+\frac{6\log(2)}{\tau\rho_{\min}}\delta_k^{\pi}\label{lem:min_policy_iterates_constanttau:proof_eq0.5}
\end{align}

Since $3\delta_k^{\pi}+\delta_k^{\phi}\leq (1-\frac{\alpha (1-\gamma) \tau \rho_{\min}^2 c^2}{32|\Scal|})^{k}(3\delta_0^{\pi}+\delta_0^{\phi})\leq3\delta_0^{\pi}+\delta_0^{\phi}\leq\frac{\rho_{\min}c^2}{64\log(2)}$, we have $\delta_k^{\pi}\leq\frac{\rho_{\min}c^2}{64\log(2)}$. Then, \eqref{lem:min_policy_iterates_constanttau:proof_eq0.5} implies
\begin{align*}
    -\left(\min_{s,a}\pi_{\theta_k}(a\mid s)\right)^2&\leq-\frac{3}{4}\left(\min_{s,a}\pi_{\tau}^{\star}(a\mid s)\right)^2+\frac{6\log(2)}{\tau\rho_{\min}}\delta_k^{\pi}\leq -\frac{3c^2}{4}+\frac{3c^2}{32}\leq-\frac{3c^2}{8}.
\end{align*}

Similarly, the property of the min and max function implies that
\begin{align*}
    \max_{s,b}(\phi_{\tau}^{\star}(b\mid s)-\phi_{\psi_k}(b\mid s))+\min_{s,b}\phi_{\psi_k}(b\mid s)\geq\min_{s,b}\phi_{\tau}^{\star}(b\mid s).
\end{align*}
Again, all three terms are non-negative, which means that the inequality is preserved after taking the square
\begin{align*}
    (\min_{s,b}\phi_{\tau}^{\star}(b\mid s))^2&\leq(\min_{s,b}\phi_{\psi_k}(b\mid s)+\max_{s,b}(\phi_{\tau}^{\star}(b\mid s)-\phi_{\psi_k}(b\mid s)))^2\notag\\
    &\leq \frac{4}{3}(\min_{s,b}\phi_{\psi_k}(b\mid s))^2+4(\max_{s,b}(\phi_{\tau}^{\star}(b\mid s)-\phi_{\psi_k}(b\mid s)))^2,
\end{align*}
which leads to
\begin{align}
    -(\min_{s,b}\phi_{\psi_k}(b\mid s))^2&\leq-\frac{3}{4}(\min_{s,b}\phi_{\tau}^{\star}(b\mid s))^2+3(\max_{s,b}(\phi_{\tau}^{\star}(b\mid s)-\phi_{\psi_k}(b\mid s)))^2\notag\\
    &\leq -\frac{3}{4}(\min_{s,b}\phi_{\tau}^{\star}(b\mid s))^2+3\|\phi_{\tau}^{\star}-\phi_{\psi_k}\|^2\notag\\
    &\leq -\frac{3}{4}(\min_{s,b}\phi_{\tau}^{\star}(b\mid s))^2+6\|\phi_{\tau}(\pi_{\theta_k})-\phi_{\psi_k}\|^2+6\|\phi_{\tau}^{\star}-\phi_{\tau}(\pi_{\theta_k})\|^2.\label{lem:min_policy_iterates_constanttau:proof_eq1}
\end{align}

From Lemma \ref{lem:quadratic_growth},
\begin{align}
    \|\phi_{\tau}(\pi_{\theta_k})-\phi_{\psi_k}\|^2\leq\frac{2\log(2)}{\tau\rho_{\min}}\left(J_{\tau}(\pi_{\theta_k},\phi_{\psi_k})-J_{\tau}(\pi_{\theta_k},\phi_{\tau}(\pi_{\theta_k}))\right)=\frac{2\log(2)}{\tau\rho_{\min}}\delta_k^{\phi},\label{lem:min_policy_iterates_constanttau:proof_eq2}
\end{align}
and
\begin{align}
    \|\phi_{\tau}^{\star}-\phi_{\tau}(\pi_{\theta_k})\|^2
    &\leq\frac{2\log(2)}{\tau\rho_{\min}}\left(J_{\tau}(\pi_{\theta_k},\phi_{\tau}^{\star})-J_{\tau}(\pi_{\theta_k},\phi_{\tau}(\pi_{\theta_k}))\right)\notag\\
    &\leq\frac{2\log(2)}{\tau\rho_{\min}}\left(J_{\tau}(\pi_{\tau}^{\star},\phi_{\tau}^{\star})-J_{\tau}(\pi_{\theta_k},\phi_{\tau}(\pi_{\theta_k}))\right)\notag\\
    &=\frac{2\log(2)}{\tau\rho_{\min}}\delta_k^{\pi},\label{lem:min_policy_iterates_constanttau:proof_eq3}
\end{align}

Using \eqref{lem:min_policy_iterates_constanttau:proof_eq2} and \eqref{lem:min_policy_iterates_constanttau:proof_eq3} in \eqref{lem:min_policy_iterates_constanttau:proof_eq1},
\begin{align*}
    -(\min_{s,b}\phi_{\psi_k}(b\mid s))^2&\leq -\frac{3}{4}(\min_{s,b}\phi_{\tau}^{\star}(b\mid s))^2+6\|\phi_{\tau}(\pi_{\theta_k})-\phi_{\psi_k}\|^2+6\|\phi_{\tau}^{\star}-\phi_{\tau}(\pi_{\theta_k})\|^2\\
    &\leq -\frac{3}{4}(\min_{s,b}\phi_{\tau}^{\star}(b\mid s))^2+\frac{12\log(2)}{\tau\rho_{\min}}\delta_k^{\phi}+\frac{12\log(2)}{\tau\rho_{\min}}\delta_k^{\pi}\\
    &=-\frac{3}{4}(\min_{s,b}\phi_{\tau}^{\star}(b\mid s))^2+\frac{12\log(2)}{\tau\rho_{\min}}(\delta_k^{\pi}+\delta_k^{\phi}).
\end{align*}

$3\delta_k^{\pi}+\delta_k^{\phi}\leq (1-\frac{\alpha (1-\gamma) \tau \rho_{\min}^2 c^2}{32|\Scal|})^{k}(3\delta_0^{\pi}+\delta_0^{\phi})\leq3\delta_0^{\pi}+\delta_0^{\phi}\leq\frac{\rho_{\min}c^2}{64\log(2)}$ guarantees $\delta_k^{\pi}+\delta_k^{\phi}\leq\frac{\rho_{\min}c^2}{32\log(2)}$. Using this in the inequality above, we have
\begin{align*}
    -(\min_{s,b}\phi_{\psi_k}(b\mid s))^2&\leq -\frac{3}{4}(\min_{s,b}\phi_{\tau}^{\star}(b\mid s))^2+\frac{12\log(2)}{\tau\rho_{\min}}(\delta_k^{\pi}+\delta_k^{\phi})\leq -\frac{3c^2}{4}+\frac{3c^2}{8}\leq-\frac{3c^2}{8}.
\end{align*}


\qed

\subsection{Proof of Lemma \ref{lem:g_smooth_constanttau}}

From Lemma \ref{lem:nonuniform_PL}, for any $\psi\in\mathbb{R}^{|\Scal|\times|\Bcal|}$
\begin{align*}
    J_{\tau}(\pi_{\theta_2},\phi_{\psi})\hspace{-2pt}-\hspace{-2pt}J_{\tau}(\pi_{\theta_2},\phi_{\tau}(\pi_{\theta_2}))\hspace{-2pt}&\leq\hspace{-2pt}\frac{|\Scal|}{2\tau\rho_{\min}\left(\min_{s,a}\phi_{\psi}(a\mid s)\right)^2}\left\|\frac{d_{\rho}^{\pi_{\theta_2},\phi_{\tau}(\pi_{\theta_2})}}{d_{\rho}^{\pi_{\theta_2},\phi_{\psi}}}\right\|_{\infty}\hspace{-3pt}\|\nabla_{\psi} J_{\tau}(\pi_{\theta_2},\phi_{\psi})\|^2\notag\\
    &\leq\hspace{-2pt} \frac{|\Scal|}{2\tau(1-\gamma)\left(\min_{s,a}\phi_{\psi}(a\mid s)\right)^2}\|\nabla_{\psi} J_{\tau}(\pi_{\theta_2},\phi_{\psi})\|^2,
\end{align*}
where the second inequality follows by an argument similar to \eqref{thm:main:proof_eq11}. Letting $\psi$ be the parameter that parameterizes $\phi_{\tau}(\pi_{\theta_1})$, we have
\begin{align*}
    &J_{\tau}(\pi_{\theta_2},\phi_{\tau}(\pi_{\theta_1}))-J_{\tau}(\pi_{\theta_2},\phi_{\tau}(\pi_{\theta_2}))\\
    &\leq \frac{|\Scal|}{2\tau(1-\gamma)\left(\min_{s,a}\phi_{\tau}(\pi_{\theta_1})(a\mid s)\right)^2}\|\nabla_{\psi} J_{\tau}(\pi_{\theta_2},\phi_{\tau}(\pi_{\theta_1}))\|^2\notag\\
    &=\frac{|\Scal|}{2\tau(1-\gamma)\left(\min_{s,a}\phi_{\tau}(\pi_{\theta_1})(a\mid s)\right)^2}\|\nabla_{\psi} J_{\tau}(\pi_{\theta_2},\psi_{\rho,\tau}^{\star}(\pi_{\theta_1}))-\nabla_{\psi} J_{\tau}(\pi_{\theta_1},\psi_{\rho,\tau}^{\star}(\pi_{\theta_1}))\|^2\notag\\
    &\leq \frac{L^2|\Scal|}{2\tau(1-\gamma)\left(\min_{s,a}\phi_{\tau}(\pi_{\theta_1})(a\mid s)\right)^2}\|\theta_1-\theta_2\|^2,
\end{align*}
where the last inequality follows from the fact that for any $\theta_1,\theta_2\in\mathbb{R}^{|\Scal|\times|\Acal|}$, $\psi_1,\psi_2\in\mathbb{R}^{|\Scal|\times|\Bcal|}$
\begin{align}
    \|\nabla_{\psi}J_{\tau}(\pi_{\theta_1},\phi_{\psi_1})-\nabla_{\psi}J_{\tau}(\pi_{\theta_2},\phi_{\psi_2})\|
    &\leq \|\nabla_{\psi}J(\pi_{\theta_1},\phi_{\psi_1})-\nabla_{\psi}J(\pi_{\theta_2},\phi_{\psi_2})\|\notag\\
    &\hspace{20pt}+\tau\|\nabla_{\psi}\Hcal_{\pi}(s,\pi_{\theta_1},\phi_{\psi_1})-\nabla_{\psi}\Hcal_{\pi}(s,\pi_{\theta_2},\phi_{\psi_2})\|\notag\\
    &\hspace{20pt}+\tau\|\nabla_{\psi}\Hcal_{\phi}(s,\pi_{\theta_1},\phi_{\psi_1})-\nabla_{\psi}\Hcal_{\phi}(s,\pi_{\theta_2},\phi_{\psi_2})\|\notag\\
    &\leq L(\|\theta_1-\theta_2\|+\|\psi_1-\psi_2\|),\label{lem:g_smooth_constanttau:proof_eq0}
\end{align}
which is a result of Lemmas~\ref{lem:value_LipschitzGrad} and \ref{lem:reg_LipschitzGrad}.

By Lemma \ref{lem:quadratic_growth}, we also have
\begin{align*}
    J_{\tau}(\pi_{\theta_2},\phi_{\tau}(\pi_{\theta_1}))-J_{\tau}(\pi_{\theta_2},\phi_{\tau}(\pi_{\theta_2}))&\geq\frac{\tau\rho_{\min}}{2\log(2)}\|\phi_{\tau}(\pi_{\theta_1})-\phi_{\tau}(\pi_{\theta_2})\|^2.
\end{align*}

Combining the two inequalities and re-arranging the terms, we have
\begin{align}
    \|\phi_{\tau}(\pi_{\theta_1})-\phi_{\tau}(\pi_{\theta_2})\|&\leq \frac{\sqrt{|\Scal|\log(2)}L}{\sqrt{(1-\gamma)\rho_{\min}}\tau\left(\min_{s,a}\phi_{\tau}(\pi_{\theta_1})(a\mid s)\right)}\|\theta_1-\theta_2\|.\label{lem:g_smooth_constanttau:proof_eq0.5}
\end{align}

Therefore, by \eqref{eq:Lipschitz_Jtau},
\begin{align*}
    &\|\nabla_{\theta} J_{\tau}(\pi_{\theta_{k}},\phi_{\tau}(\pi_{\theta_k}))-\nabla_{\theta} J_{\tau}(\pi_{\theta_{k+1}},\phi_{\tau}(\pi_{\theta_{k+1}}))\|\\
    &\leq L\|\theta_k-\theta_{k+1}\|+L\|\phi_{\tau}(\pi_{\theta_k})-\phi_{\tau}(\pi_{\theta_{k+1}})\|\\
    &\leq L\left(1+\frac{\sqrt{|\Scal|\log(2)}L}{\sqrt{(1-\gamma)\rho_{\min}}\tau\left(\min_{s,a}\phi_{\tau}(\pi_{\theta_k})(a\mid s)\right)}\right)\|\theta_k-\theta_{k+1}\|
\end{align*}

Due to the Danskin's Theorem \eqref{eq:Danskin}, this implies that we can perform the expansion
\begin{align}
    &J_{\tau}(\pi_{\theta_k},\phi_{\tau}(\pi_{\theta_k}))-J_{\tau}(\pi_{\theta_{k+1}},\phi_{\tau}(\pi_{\theta_{k+1}}))\notag\\
    &\leq-\langle\nabla_{\theta} J_{\tau}(\pi_{\theta_{k}},\phi_{\tau}(\pi_{\theta_k})),\theta_{k+1}-\theta_k\rangle\notag\\
    &\hspace{20pt}+\frac{L}{2}\left(1+\frac{\sqrt{|\Scal|\log(2)}L}{\sqrt{(1-\gamma)\rho_{\min}}\tau\left(\min_{s,a}\phi_{\tau}(\pi_{\theta_k})(a\mid s)\right)}\right)\|\theta_{k+1}-\theta_k\|^2\notag\\
    &\leq -\alpha_k\langle\nabla_{\theta} J_{\tau}(\pi_{\theta_{k}},\phi_{\tau}(\pi_{\theta_k})),\nabla_{\theta}J_{\tau}(\pi_{\theta_k},\phi_{\psi_k})\rangle\notag\\
    &\hspace{20pt}+\frac{L\alpha_k^2}{2}\left(1+\frac{\sqrt{|\Scal|\log(2)}L}{\sqrt{(1-\gamma)\rho_{\min}}\tau\left(\min_{s,a}\phi_{\tau}(\pi_{\theta_k})(a\mid s)\right)}\right)\|\nabla_{\theta}J_{\tau}(\pi_{\theta_k},\phi_{\psi_k})\|^2.
    \label{lem:g_smooth_constanttau:proof_eq1}
\end{align}

Note that by the property of the min function
\begin{align}
    \min_{s,a}\phi_{\tau}(\pi_{\theta_k})(a\mid s)&\geq\min_{s,a}\phi_{\tau}^{\star}(a\mid s)-\max_{s,a}(\phi_{\tau}^{\star}(a\mid s)-\phi_{\tau}(\pi_{\theta_k})(a\mid s))\notag\\
    &\geq \min_{s,a}\phi_{\tau}^{\star}(a\mid s)-\|\phi_{\tau}^{\star}-\phi_{\tau}(\pi_{\theta_k})\|\notag\\
    &\geq c-\sqrt{\frac{2\log(2)}{\tau\rho_{\min}}(\delta_k^{\pi}+\delta_k^{\phi})},\label{lem:g_smooth_constanttau:proof_eq1.5}
\end{align}
where the last inequality uses the same argument as in \eqref{lem:min_policy_iterates:proof_eq3}. 
Since \eqref{thm:main:proof_eq1} implies $\delta_k^{\pi}+\delta_k^{\phi}\leq\frac{\rho_{\min}c^2\tau}{64\log(2)(k+1)^{1/3}}$, we further have
\begin{align*}
    \min_{s,a}\phi_{\tau}(\pi_{\theta_k})(a\mid s)&\geq c-\sqrt{\frac{2\log(2)}{\tau\rho_{\min}}(\delta_k^{\pi}+\delta_k^{\phi})}\geq c(1-\sqrt{\frac{1}{32}})\geq \frac{c\sqrt{\log(2)}}{2}.
\end{align*}
Using this bound in \eqref{lem:g_smooth_constanttau:proof_eq1},
\begin{align}
    &J_{\tau}(\pi_{\theta_k},\phi_{\tau}(\pi_{\theta_k}))-J_{\tau}(\pi_{\theta_{k+1}},\phi_{\tau}(\pi_{\theta_{k+1}}))\notag\\
    &\leq -\alpha_k\langle\nabla_{\theta} J_{\tau}(\pi_{\theta_{k}},\phi_{\tau}(\pi_{\theta_k})),\nabla_{\theta}J_{\tau}(\pi_{\theta_k},\phi_{\psi_k})\rangle\notag\\
    &\hspace{20pt}+\frac{L\alpha_k^2}{2}\left(1+\frac{\sqrt{|\Scal|\log(2)}L}{\sqrt{(1-\gamma)\rho_{\min}}\tau\left(\min_{s,a}\phi_{\tau}(\pi_{\theta_1})(a\mid s)\right)}\right)\|\nabla_{\theta}J_{\tau}(\pi_{\theta_k},\phi_{\psi_k})\|^2\notag\\
    &\leq -\alpha_k\langle\nabla_{\theta} J_{\tau}(\pi_{\theta_{k}},\phi_{\tau}(\pi_{\theta_k})),\nabla_{\theta}J_{\tau}(\pi_{\theta_k},\phi_{\psi_k})\rangle\notag\\
    &\hspace{20pt}+\frac{L\alpha_k^2}{2}\left(1+\frac{2\sqrt{|\Scal|}L}{\sqrt{(1-\gamma)\rho_{\min}}\tau c}\right)\|\nabla_{\theta}J_{\tau}(\pi_{\theta_k},\phi_{\psi_k})\|^2, \label{lem:g_smooth_constanttau:proof_eq2}
\end{align}

With the step size choice $\alpha_k\leq\left(L+\frac{2\sqrt{|\Scal|}L^2}{\sqrt{(1-\gamma)\rho_{\min}}\tau c}\right)^{-1}$, we get
\begin{align*}
    &J_{\tau}(\pi_{\theta_k},\phi_{\tau}(\pi_{\theta_k}))-J_{\tau}(\pi_{\theta_{k+1}},\phi_{\tau}(\pi_{\theta_{k+1}}))\notag\\
    &\leq -\alpha_k\langle\nabla_{\theta} J_{\tau}(\pi_{\theta_{k}},\phi_{\tau}(\pi_{\theta_k})),\nabla_{\theta}J_{\tau}(\pi_{\theta_k},\phi_{\psi_k})\rangle\notag\\
    &\hspace{20pt}+\frac{L\alpha_k^2}{2}\left(1+\frac{2\sqrt{|\Scal|}L}{\sqrt{(1-\gamma)\rho_{\min}}\tau c}\right)\|\nabla_{\theta}J_{\tau}(\pi_{\theta_k},\phi_{\psi_k})\|^2\\
    &\leq -\alpha_k\langle\nabla_{\theta} J_{\tau}(\pi_{\theta_{k}},\phi_{\tau}(\pi_{\theta_k})),\nabla_{\theta}J_{\tau}(\pi_{\theta_k},\phi_{\psi_k})\rangle\notag\\
    &\hspace{20pt}+\frac{\alpha_k}{2}\|\nabla_{\theta}J_{\tau}(\pi_{\theta_k},\phi_{\psi_k})\|^2\\
    &= \frac{\alpha_k}{2}\|\nabla_{\theta} J_{\tau}(\pi_{\theta_{k}},\phi_{\tau}(\pi_{\theta_k}))-\nabla_{\theta}J_{\tau}(\pi_{\theta_k},\phi_{\psi_k})\|^2-\|\nabla_{\theta} J_{\tau}(\pi_{\theta_{k}},\phi_{\tau}(\pi_{\theta_k}))\|^2.
\end{align*}

\qed

\subsection{Proof of Lemma \ref{lem:min_policy_iterates}}

The property of the min and max function implies that
\begin{align*}
    \max_{s,a}(\pi_{\tau_k}^{\star}(a\mid s)-\pi_{\theta_k}(a\mid s))+\min_{s,a}\pi_{\theta_k}(a\mid s)\geq\min_{s,a}\pi_{\tau_k}^{\star}(a\mid s).
\end{align*}
Since the three terms are all non-negative, the inequality holds after taking the square
\begin{align*}
    (\min_{s,a}\pi_{\tau_k}^{\star}(a\mid s))^2&\leq(\max_{s,a}(\pi_{\tau_k}^{\star}(a\mid s)-\pi_{\theta_k}(a\mid s))+\min_{s,a}\pi_{\theta_k}(a\mid s))^2\notag\\
    &\leq \frac{4}{3}(\min_{s,a}\pi_{\theta_k}(a\mid s))^2+4(\max_{s,a}(\pi_{\tau_k}^{\star}(a\mid s)-\pi_{\theta_k}(a\mid s)))^2.
\end{align*}
Re-arranging the terms,
\begin{align*}
    -\left(\min_{s,a}\pi_{\theta_k}(a\mid s)\right)^2&\leq-\frac{3}{4}\left(\min_{s,a}\pi_{\tau_k}^{\star}(a\mid s)\right)^2+3\left(\max_{s,a}\pi_{\tau_k}^{\star}(a\mid s)-\pi_{\phi_k}(a\mid s)\right)^2\notag\\
    &\leq-\frac{3}{4}\left(\min_{s,a}\pi_{\tau_k}^{\star}(a\mid s)\right)^2+3\|\pi_{\tau_k}^{\star}-\pi_{\phi_k}\|^2
\end{align*}

From Lemma \ref{lem:quadratic_growth},
\begin{align}
    -\left(\min_{s,a}\pi_{\theta_k}(a\mid s)\right)^2&\leq-\frac{3}{4}\left(\min_{s,a}\pi_{\tau_k}^{\star}(a\mid s)\right)^2+3\|\pi_{\tau_k}^{\star}-\pi_{\phi_k}\|^2\notag\\
    &\leq-\frac{3}{4}\left(\min_{s,a}\pi_{\tau_k}^{\star}(a\mid s)\right)^2+\frac{6\log(2)}{\tau_k\rho_{\min}}(J_{\tau_k}(\pi_{\tau_k}^{\star},\phi_{\tau_k}^{\star})-J_{\tau_k}(\pi_{\theta_k},\phi_{\tau_k}^{\star}))\notag\\
    &\leq-\frac{3}{4}\left(\min_{s,a}\pi_{\tau_k}^{\star}(a\mid s)\right)^2+\frac{6\log(2)}{\tau_k\rho_{\min}}(J_{\tau_k}(\pi_{\tau_k}^{\star},\phi_{\tau_k}^{\star})-J_{\tau_k}(\pi_{\theta_k},\phi_{\tau_k}(\pi_{\theta_k})))\notag\\
    &=-\frac{3}{4}\left(\min_{s,a}\pi_{\tau_k}^{\star}(a\mid s)\right)^2+\frac{6\log(2)}{\tau_k\rho_{\min}}\delta_k^{\pi},\label{lem:min_policy_iterates:proof_eq0.5}
\end{align}

Since $3\delta_k^{\pi}+\delta_k^{\phi}\leq\frac{\rho\tau_k c^2}{64\log(2)}$, we have $\delta_k^{\pi}\leq\frac{\rho\tau_k c^2}{64\log(2)}$, which along with \eqref{lem:min_policy_iterates:proof_eq0.5} implies
\begin{align*}
    -\left(\min_{s,a}\pi_{\theta_k}(a\mid s)\right)^2&\leq-\frac{3}{4}\left(\min_{s,a}\pi_{\tau_k}^{\star}(a\mid s)\right)^2+\frac{6\log(2)}{\tau_k\rho_{\min}}\delta_k^{\pi}\leq-\frac{3c^2}{4}+\frac{3c^2}{32}\leq-\frac{3c^2}{8}.
\end{align*}

Similarly, the property of the min and max function implies that
\begin{align*}
    \max_{s,b}(\phi_{\tau_k}^{\star}(b\mid s)-\phi_{\psi_k}(b\mid s))+\min_{s,b}\phi_{\psi_k}(b\mid s)\geq\min_{s,b}\phi_{\tau_k}^{\star}(b\mid s).
\end{align*}
Again, all three terms are non-negative, which means that the inequality is preserved after taking the square
\begin{align*}
    (\min_{s,b}\phi_{\tau_k}^{\star}(b\mid s))^2&\leq(\min_{s,b}\phi_{\psi_k}(b\mid s)+\max_{s,b}(\phi_{\tau_k}^{\star}(b\mid s)-\phi_{\psi_k}(b\mid s)))^2\notag\\
    &\leq \frac{4}{3}(\min_{s,b}\phi_{\psi_k}(b\mid s))^2+4(\max_{s,b}(\phi_{\tau_k}^{\star}(b\mid s)-\phi_{\psi_k}(b\mid s)))^2,
\end{align*}
which leads to
\begin{align}
    -(\min_{s,b}\phi_{\psi_k}(b\mid s))^2&\leq-\frac{3}{4}(\min_{s,b}\phi_{\tau_k}^{\star}(b\mid s))^2+3(\max_{s,b}(\phi_{\tau_k}^{\star}(b\mid s)-\phi_{\psi_k}(b\mid s)))^2\notag\\
    &\leq -\frac{3}{4}(\min_{s,b}\phi_{\tau_k}^{\star}(b\mid s))^2+3\|\phi_{\tau_k}^{\star}-\phi_{\psi_k}\|^2\notag\\
    &\leq -\frac{3}{4}(\min_{s,b}\phi_{\tau_k}^{\star}(b\mid s))^2+6\|\phi_{\tau_k}(\pi_{\theta_k})-\phi_{\psi_k}\|^2+6\|\phi_{\tau_k}^{\star}-\phi_{\tau_k}(\pi_{\theta_k})\|^2.\label{lem:min_policy_iterates:proof_eq1}
\end{align}

From Lemma \ref{lem:quadratic_growth},
\begin{align}
    \|\phi_{\tau_k}(\pi_{\theta_k})-\phi_{\psi_k}\|^2\leq\frac{2\log(2)}{\tau_k\rho_{\min}}\left(J_{\tau_k}(\pi_{\theta_k},\phi_{\psi_k})-J_{\tau_k}(\pi_{\theta_k},\phi_{\tau_k}(\pi_{\theta_k}))\right)=\frac{2\log(2)}{\tau_k\rho_{\min}}\delta_k^{\phi},\label{lem:min_policy_iterates:proof_eq2}
\end{align}
and
\begin{align}
    &\|\phi_{\tau_k}^{\star}-\phi_{\tau_k}(\pi_{\theta_k})\|^2\notag\\
    &\leq\frac{2\log(2)}{\tau_k\rho_{\min}}\left(J_{\tau_k}(\pi_{\theta_k},\phi_{\tau_k}^{\star})-J_{\tau_k}(\pi_{\theta_k},\phi_{\tau_k}(\pi_{\theta_k}))\right)\notag\\
    &\leq \frac{2\log(2)}{\tau_k\rho_{\min}}\Big(\big(J_{\tau_k}(\pi_{\theta_k},\phi_{\tau_k}^{\star})-J_{\tau_k}(\pi_{\theta_k},\phi_{\psi_k})\big)+\underbrace{\big(J_{\tau_k}(\pi_{\theta_k},\phi_{\psi_k})-J_{\tau_k}(\pi_{\theta_k},\phi_{\tau_k}(\pi_{\theta_k}))\big)}_{\delta_k^{\phi}}\Big)\notag\\
    &=\frac{2\log(2)}{\tau_k\rho_{\min}}\left(\left(J_{\tau_k}(\pi_{\theta_k},\phi_{\tau_k}^{\star})\hspace{-2pt}-\hspace{-2pt}J_{\tau_k}(\pi_{\theta_k},\phi_{\tau_k}(\pi_{\theta_k}))\right)+\left(J_{\tau_k}(\pi_{\theta_k},\phi_{\tau_k}(\pi_{\theta_k}))\hspace{-2pt}-\hspace{-2pt}J_{\tau_k}(\pi_{\theta_k},\phi_{\psi_k})\right)+\delta_k^{\phi}\right)\notag\\
    &\leq\frac{2\log(2)}{\tau_k\rho_{\min}}\left(J_{\tau_k}(\pi_{\theta_k},\phi_{\tau_k}^{\star})-J_{\tau_k}(\pi_{\theta_k},\phi_{\tau_k}(\pi_{\theta_k}))+\delta_k^{\phi}\right)\notag\\
    &\leq\frac{2\log(2)}{\tau_k\rho_{\min}}\left(J_{\tau_k}(\pi_{\tau_k}^{\star},\phi_{\tau_k}^{\star})-J_{\tau_k}(\pi_{\theta_k},\phi_{\tau_k}(\pi_{\theta_k}))+\delta_k^{\phi}\right)\notag\\
    &=\frac{2\log(2)}{\tau_k\rho_{\min}}\left(\delta_k^{\pi}+\delta_k^{\phi}\right),\label{lem:min_policy_iterates:proof_eq3}
\end{align}
where the third inequality follows from $J_{\tau_k}(\pi_{\theta_k},\phi_{\tau_k}(\pi_{\theta_k}))-J_{\tau_k}(\pi_{\theta_k},\phi_{\psi_k})\leq0$.

Using \eqref{lem:min_policy_iterates:proof_eq2} and \eqref{lem:min_policy_iterates:proof_eq3} in \eqref{lem:min_policy_iterates:proof_eq1},
\begin{align*}
    -(\min_{s,b}\phi_{\psi_k}(b\mid s))^2&\leq -\frac{3}{4}(\min_{s,b}\phi_{\tau_k}^{\star}(b\mid s))^2+6\|\phi_{\tau_k}(\pi_{\theta_k})-\phi_{\psi_k}\|^2+6\|\phi_{\tau_k}^{\star}-\phi_{\tau_k}(\pi_{\theta_k})\|^2\\
    &\leq -\frac{3}{4}(\min_{s,b}\phi_{\tau_k}^{\star}(b\mid s))^2+\frac{12\log(2)}{\tau_k\rho_{\min}}\delta_k^{\phi}+\frac{12\log(2)}{\tau_k\rho_{\min}}(\delta_k^{\pi}+\delta_k^{\phi})\\
    &=-\frac{3}{4}(\min_{s,b}\phi_{\tau_k}^{\star}(b\mid s))^2+\frac{12\log(2)}{\tau_k\rho_{\min}}(\delta_k^{\pi}+2\delta_k^{\phi}).
\end{align*}

$3\delta_k^{\pi}+\delta_k^{\phi}\leq\frac{\rho\tau_k c^2}{64\log(2)}$ implies that $\delta_k^{\pi}+2\delta_k^{\phi}\leq\frac{\rho\tau_k c^2}{32\log(2)}$. Using this in the inequality above,
\begin{align*}
    -(\min_{s,b}\phi_{\psi_k}(b\mid s))^2&\leq-\frac{3}{4}(\min_{s,b}\phi_{\tau_k}^{\star}(b\mid s))^2+\frac{12\log(2)}{\tau_k\rho_{\min}}(\delta_k^{\pi}+2\delta_k^{\phi})\leq-\frac{3c^2}{4}+\frac{12c^2}{32}\leq-\frac{3c^2}{8}.
\end{align*}


\qed

\subsection{Proof of Lemma \ref{lem:g_smooth}}

From Lemma \ref{lem:nonuniform_PL}, for any $\psi\in\mathbb{R}^{|\Scal|\times|\Bcal|}$
\begin{align*}
    &J_{\tau_k}(\pi_{\theta_2},\phi_{\psi})-J_{\tau_k}(\pi_{\theta_2},\phi_{\tau_k}(\pi_{\theta_2}))\\
    &\leq\frac{|\Scal|}{2\tau_k\rho_{\min}\left(\min_{s,a}\phi_{\psi}(a\mid s)\right)^2}\left\|\frac{d_{\rho}^{\pi_{\theta_2},\phi_{\tau_k}(\pi_{\theta_2})}}{d_{\rho}^{\pi_{\theta_2},\phi_{\psi}}}\right\|_{\infty}\|\nabla_{\psi} J_{\tau_k}(\pi_{\theta_2},\phi_{\psi})\|^2\notag\\
    &\leq \frac{|\Scal|}{2\tau_k(1-\gamma)\left(\min_{s,a}\phi_{\psi}(a\mid s)\right)^2}\|\nabla_{\psi} J_{\tau_k}(\pi_{\theta_2},\phi_{\psi})\|^2,
\end{align*}
where the second inequality follows by an argument similar to \eqref{thm:main:proof_eq11}. Letting $\psi$ be the parameter that parameterizes $\phi_{\tau_k}(\pi_{\theta_1})$ and defining $L_k=L_{\Hcal}(2\tau_k+1)$, we have
\begin{align*}
    &J_{\tau_k}(\pi_{\theta_2},\phi_{\tau_k}(\pi_{\theta_1}))-J_{\tau_k}(\pi_{\theta_2},\phi_{\tau_k}(\pi_{\theta_2}))\\
    &\leq \frac{|\Scal|}{2\tau_k(1-\gamma)\left(\min_{s,a}\phi_{\tau_k}(\pi_{\theta_1})(a\mid s)\right)^2}\|\nabla_{\psi} J_{\tau_k}(\pi_{\theta_2},\phi_{\tau_k}(\pi_{\theta_1}))\|^2\notag\\
    &=\frac{|\Scal|}{2\tau_k(1-\gamma)\left(\min_{s,a}\phi_{\tau_k}(\pi_{\theta_1})(a\mid s)\right)^2}\|\nabla_{\psi} J_{\tau_k}(\pi_{\theta_2},\psi_{\rho,\tau_k}^{\star}(\pi_{\theta_1}))-\nabla_{\psi} J_{\tau_k}(\pi_{\theta_1},\psi_{\rho,\tau_k}^{\star}(\pi_{\theta_1}))\|^2\notag\\
    &\leq \frac{L_k^2|\Scal|}{2\tau_k(1-\gamma)\left(\min_{s,a}\phi_{\tau_k}(\pi_{\theta_1})(a\mid s)\right)^2}\|\theta_1-\theta_2\|^2,
\end{align*}
where the last inequality uses the same argument as \eqref{lem:g_smooth_constanttau:proof_eq0}.

By Lemma \ref{lem:quadratic_growth}, we also have
\begin{align*}
    J_{\tau_k}(\pi_{\theta_2},\phi_{\tau_k}(\pi_{\theta_1}))-J_{\tau_k}(\pi_{\theta_2},\phi_{\tau_k}(\pi_{\theta_2}))&\geq\frac{\tau_k\rho_{\min}}{2\log(2)}\|\phi_{\tau_k}(\pi_{\theta_1})-\phi_{\tau_k}(\pi_{\theta_2})\|^2.
\end{align*}

Combining the two inequalities and re-arranging the terms, we have
\begin{align}
    \|\phi_{\tau_k}(\pi_{\theta_1})-\phi_{\tau_k}(\pi_{\theta_2})\|&\leq \frac{\sqrt{|\Scal|\log(2)}L_k}{\sqrt{(1-\gamma)\rho_{\min}}\tau_k\left(\min_{s,a}\phi_{\tau_k}(\pi_{\theta_1})(a\mid s)\right)}\|\theta_1-\theta_2\|.\label{lem:g_smooth:proof_eq0.5}
\end{align}

Therefore, by \eqref{eq:Lipschitz_Jtau},
\begin{align*}
    &\|\nabla_{\theta} J_{\tau_k}(\pi_{\theta_{k}},\phi_{\tau_k}(\pi_{\theta_k}))-\nabla_{\theta} J_{\tau_k}(\pi_{\theta_{k+1}},\phi_{\tau_k}(\pi_{\theta_{k+1}}))\|\\
    &\leq L_k\|\theta_k-\theta_{k+1}\|+L_k\|\phi_{\tau_k}(\pi_{\theta_k})-\phi_{\tau_k}(\pi_{\theta_{k+1}})\|\\
    &\leq L_k\left(1+\frac{\sqrt{|\Scal|\log(2)}L_k}{\sqrt{(1-\gamma)\rho_{\min}}\tau_k\left(\min_{s,a}\phi_{\tau_k}(\pi_{\theta_k})(a\mid s)\right)}\right)\|\theta_k-\theta_{k+1}\|
\end{align*}

Due to the Danskin's Theorem \eqref{eq:Danskin}, this implies that we can perform the expansion
\begin{align}
    &J_{\tau_{k}}(\pi_{\theta_k},\phi_{\tau_k}(\pi_{\theta_k}))-J_{\tau_{k}}(\pi_{\theta_{k+1}},\phi_{\tau_k}(\pi_{\theta_{k+1}}))\notag\\
    &\leq-\langle\nabla_{\theta} J_{\tau_k}(\pi_{\theta_{k}},\phi_{\tau_k}(\pi_{\theta_k})),\theta_{k+1}-\theta_k\rangle\notag\\
    &\hspace{20pt}+\frac{L_k}{2}\left(1+\frac{\sqrt{|\Scal|\log(2)}L_k}{\sqrt{(1-\gamma)\rho_{\min}}\tau_k\left(\min_{s,a}\phi_{\tau_k}(\pi_{\theta_k})(a\mid s)\right)}\right)\|\theta_{k+1}-\theta_k\|^2\notag\\
    &\leq -\alpha_k\langle\nabla_{\theta} J_{\tau_k}(\pi_{\theta_{k}},\phi_{\tau_k}(\pi_{\theta_k})),\nabla_{\theta}J_{\tau_k}(\pi_{\theta_k},\phi_{\psi_k})\rangle\notag\\
    &\hspace{20pt}+\frac{L_k\alpha_k^2}{2}\left(1+\frac{\sqrt{|\Scal|\log(2)}L_k}{\sqrt{(1-\gamma)\rho_{\min}}\tau_k\left(\min_{s,a}\phi_{\tau_k}(\pi_{\theta_k})(a\mid s)\right)}\right)\|\nabla_{\theta}J_{\tau_k}(\pi_{\theta_k},\phi_{\psi_k})\|^2.
    \label{lem:g_smooth:proof_eq1}
\end{align}

Note that by the property of the min function
\begin{align}
    \min_{s,a}\phi_{\tau_k}(\pi_{\theta_k})(a\mid s)&\geq\min_{s,a}\phi_{\tau_k}^{\star}(a\mid s)-\max_{s,a}(\phi_{\tau_k}^{\star}(a\mid s)-\phi_{\tau_k}(\pi_{\theta_k})(a\mid s))\notag\\
    &\geq \min_{s,a}\phi_{\tau_k}^{\star}(a\mid s)-\|\phi_{\tau_k}^{\star}-\phi_{\tau_k}(\pi_{\theta_k})\|\notag\\
    &\geq c-\sqrt{\frac{2\log(2)}{\tau_k\rho_{\min}}(\delta_k^{\pi}+\delta_k^{\phi})},\label{lem:g_smooth:proof_eq1.5}
\end{align}
where the last inequality uses the same argument as in \eqref{lem:min_policy_iterates:proof_eq3}. 
Since \eqref{thm:main:proof_eq1} implies $\delta_k^{\pi}+\delta_k^{\phi}\leq\frac{\rho_{\min}c^2\tau_0}{64\log(2)(k+1)^{1/3}}$, we further have
\begin{align*}
    \min_{s,a}\phi_{\tau_k}(\pi_{\theta_k})(a\mid s)&\geq c-\sqrt{\frac{2\log(2)}{\tau_k\rho_{\min}}(\delta_k^{\pi}+\delta_k^{\phi})}\geq c(1-\sqrt{\frac{1}{32}})\geq \frac{c\sqrt{\log(2)}}{2}.
\end{align*}
Using this bound in \eqref{lem:g_smooth:proof_eq1},
\begin{align}
    &J_{\tau_{k}}(\pi_{\theta_k},\phi_{\tau_k}(\pi_{\theta_k}))-J_{\tau_{k}}^{\pi_{\theta_{k+1}},\phi_{\tau_k}(\pi_{\theta_{k+1}})}(\rho)\notag\\
    &\leq -\alpha_k\langle\nabla_{\theta} J_{\tau_k}(\pi_{\theta_{k}},\phi_{\tau_k}(\pi_{\theta_k})),\nabla_{\theta}J_{\tau_k}(\pi_{\theta_k},\phi_{\psi_k})\rangle\notag\\
    &\hspace{20pt}+\frac{L_k\alpha_k^2}{2}\left(1+\frac{\sqrt{|\Scal|\log(2)}L_k}{\sqrt{(1-\gamma)\rho_{\min}}\tau_k\left(\min_{s,a}\phi_{\tau_k}(\pi_{\theta_1})(a\mid s)\right)}\right)\|\nabla_{\theta}J_{\tau_k}(\pi_{\theta_k},\phi_{\psi_k})\|^2\notag\\
    &\leq -\alpha_k\langle\nabla_{\theta} J_{\tau_k}(\pi_{\theta_{k}},\phi_{\tau_k}(\pi_{\theta_k})),\nabla_{\theta}J_{\tau_k}(\pi_{\theta_k},\phi_{\psi_k})\rangle\notag\\
    &\hspace{20pt}+\frac{L_k\alpha_k^2}{2}\left(1+\frac{2\sqrt{|\Scal|}L_k}{\sqrt{(1-\gamma)\rho_{\min}}\tau_kc}\right)\|\nabla_{\theta}J_{\tau_k}(\pi_{\theta_k},\phi_{\psi_k})\|^2\notag\\
    &\leq-\alpha_k\langle\nabla_{\theta} J_{\tau_k}(\pi_{\theta_{k}},\phi_{\tau_k}(\pi_{\theta_k})),\nabla_{\theta}J_{\tau_k}(\pi_{\theta_k},\phi_{\psi_k})\rangle+\frac{\alpha_k^2}{2}\left(L_k+\frac{C_2 L_k^2}{\tau_k}\right)\|\nabla_{\theta}J_{\tau_k}(\pi_{\theta_k},\phi_{\psi_k})\|^2.\label{lem:g_smooth:proof_eq2}
\end{align}

The condition on $h$, which is $\frac{\alpha_0}{h^{2/3}}\leq(2L_{\Hcal}+4L_{\Hcal}^2 C_2)\frac{\tau_0}{h^{1/3}}+(L_{\Hcal}+4L_{\Hcal}^2C_2)+\frac{L_{\Hcal}^2 C_2 h^{1/3}}{\tau_0}$, can be equivalently expressed as $\alpha_0\left(L_0+\frac{C_2 L_0^2}{\tau_0}\right)\leq1$. Since $\alpha_k$ decays faster than $\tau_k$, this guarantees for all $k\geq 0$
\begin{align*}
    \alpha_k\left(L_k+\frac{C_2 L_k^2}{\tau_k}\right)\leq 1.
\end{align*}
Using this inequality in \eqref{lem:g_smooth:proof_eq2}, we get
\begin{align*}
    &J_{\tau_{k}}(\pi_{\theta_k},\phi_{\tau_k}(\pi_{\theta_k}))-J_{\tau_{k}}(\pi_{\theta_{k+1}},\phi_{\tau_k}(\pi_{\theta_{k+1}}))\notag\\
    &\leq -\alpha_k\langle\nabla_{\theta} J_{\tau_k}(\pi_{\theta_{k}},\phi_{\tau_k}(\pi_{\theta_k})),\nabla_{\theta}J_{\tau_k}(\pi_{\theta_k},\phi_{\psi_k})\rangle+\frac{\alpha_k}{2}\|\nabla_{\theta}J_{\tau_k}(\pi_{\theta_k},\phi_{\psi_k})\|^2\\
    &= \frac{\alpha_k}{2}\left(\|\nabla_{\theta} J_{\tau_k}(\pi_{\theta_{k}},\phi_{\tau_k}(\pi_{\theta_k}))-\nabla_{\theta}J_{\tau_k}(\pi_{\theta_k},\phi_{\psi_k})\|^2-\|\nabla_{\theta} J_{\tau_k}(\pi_{\theta_{k}},\phi_{\tau_k}(\pi_{\theta_k}))\|^2\right).
\end{align*}

\qed

%% file: Remark_InitialCondition.tex
\section{Discussion on the Initial Condition for Corollary~\ref{cor:GDA_piecewiseconstanttau} and Theorem~\ref{thm:main}}\label{remark:large_tau0}
We show that as $\tau_0\hspace{-2.5pt}\rightarrow\hspace{-2.5pt}\infty$, both $\delta^{\pi}_0\hspace{-1pt}/\hspace{-1pt}\tau_0$ and $\delta^{\phi}_0\hspace{-1pt}/\hspace{-1pt}\tau_0$ approach $0$. 
Decomposing $\delta^{\pi}_0\hspace{-1pt}/\hspace{-1pt}\tau_0$, we have
\begin{align}
    \frac{\delta^{\pi}_0}{\tau_0}\hspace{-2pt}&= 
    \frac{1}{\tau_0}\left(J(\pi_{\tau_0}^{\star},\phi_{\tau_0}^{\star}) - J(\pi_{\theta_0},\phi_{\tau_0}(\pi_{\theta_0}))\right)\label{remark:large_tau0:eq1}\\
    &+\left(\Hcal_{\pi}(\rho,\pi_{\tau_0}^{\star},\phi_{\tau_0}^{\star})\hspace{-1pt}-\hspace{-1pt}\Hcal_{\pi}(\rho,\pi_{\theta_0},\phi_{\tau_0}(\pi_{\theta_0}))\right)\hspace{-1pt}+\hspace{-1pt}\left(\Hcal_{\phi}(\rho,\pi_{\tau_0}^{\star},\phi_{\tau_0}^{\star})\hspace{-1pt}-\hspace{-1pt}\Hcal_{\phi}(\rho,\pi_{\theta_0},\phi_{\tau_0}(\pi_{\theta_0}))\right)\label{remark:large_tau0:eq2}
\end{align}

The original value functions are bounded within $[0,\frac{1}{1-\gamma}]$, which implies that the term \eqref{remark:large_tau0:eq1} decays inversely with $\tau_0$ in the worst case.
When $\tau_0\rightarrow\infty$, the Nash equilibrium policy pair $\pi_{\tau_0}^{\star}$ and $\phi_{\tau_0}^{\star}$ both approach the uniform distribution, and so does $\phi_{\tau_0}(\pi_{\theta_0})$. This means that \eqref{remark:large_tau0:eq2} approaches 0. Therefore, as the sum of \eqref{remark:large_tau0:eq1} and \eqref{remark:large_tau0:eq2}, $\delta_k^{\pi}/\tau_0$ decays to 0 as $\tau_0\rightarrow\infty$. A similar argument can be used for $\delta^{\phi}_0/\tau_0$.

%% file: ExperimentDetails.tex
\section{Experiment Details}\label{sec:experimentdeatils}
We first discuss the design of the completely mixed Markov game.
The dimension of state space is 2, and so is the dimension of the action spaces of both players. Using $s_1$, $s_2$ to denote the two states, 
we can essentially describe $\Pcal$ as a $2\times2\times2\times2$ tensor where $\Pcal(s'\mid s,\cdot,\cdot)$ is a $2\times2$ matrix for any $s,s'\in\Scal$ with rows corresponding to the action of the first player and columns corresponding to the second player
\begin{align*}
    \Pcal(s_1\mid s_1,\cdot,\cdot)=\left[\begin{array}{ll}
    0.2 & 0.5 \\
    0.5 & 0.1
    \end{array}\right],\quad \Pcal(s_2\mid s_1,\cdot,\cdot)=\left[\begin{array}{ll}
    0.8 & 0.5 \\
    0.5 & 0.9
    \end{array}\right],\notag\\
    \Pcal(s_1\mid s_2,\cdot,\cdot)=\left[\begin{array}{ll}
    0.3 & 0.2 \\
    0.6 & 0.2
    \end{array}\right],\quad \Pcal(s_2\mid s_2,\cdot,\cdot)=\left[\begin{array}{ll}
    0.7 & 0.8 \\
    0.4 & 0.8
    \end{array}\right].
\end{align*}

Similarly, the reward function can be described by a $2\times2\times2$ tensor where $r(s,\cdot,\cdot)$ is a $2\times2$ matrix for any $s\in\Scal$ with rows corresponding to the action of the first player and columns corresponding to the second player
\begin{align*}
    r(s_1,\cdot,\cdot)=\left[\begin{array}{ll}
    1 & 2 \\
    2 & 1
    \end{array}\right],\quad r(s_2,\cdot,\cdot)=\left[\begin{array}{ll}
    6 & 4 \\
    3 & 10
    \end{array}\right].
\end{align*}

Under the initial distribution $\rho=[0.5,0.5]^{\top}$ and discount factor $\gamma=0.9$, the (approximate) Nash equilibrium of this Markov game is
\begin{align*}
    \pi^{\star}(\cdot\mid s_1)=[0.812, 0.188],\quad \pi^{\star}(\cdot\mid s_2)=[0.837, 0.163],\\
    \phi^{\star}(\cdot\mid s_1)=[0.880, 0.120],\quad \phi^{\star}(\cdot\mid s_2)=[0.597, 0.403].
\end{align*}

To design the Markov game that does not observe Assumption~\ref{assump:NE_completelymixed}, we use the same transition probability matrices as in the completely mixed Markov game case. The reward function is
\begin{align*}
    r(s_1,\cdot,\cdot)=\left[\begin{array}{ll}
    1 & 2 \\
    3 & 4
    \end{array}\right],\quad r(s_2,\cdot,\cdot)=\left[\begin{array}{ll}
    1 & 2 \\
    3 & 4
    \end{array}\right].
\end{align*}

Under the initial distribution $\rho=[0.5,0.5]^{\top}$ and discount factor $\gamma=0.9$, it can be easily seen that the Nash equilibrium of this Markov game is unique and is
\begin{align*}
    \pi^{\star}(\cdot\mid s_1)=[0, 1],\quad \pi^{\star}(\cdot\mid s_2)=[0, 1],\\
    \phi^{\star}(\cdot\mid s_1)=[1, 0],\quad \phi^{\star}(\cdot\mid s_2)=[1, 0].
\end{align*}
Since the Nash equilibrium consists of a pair of deterministic policies, Assumption~\ref{assump:NE_completelymixed} is not satisfied in this case. 